\documentclass[11pt,a4paper]{amsart}
\usepackage[left=2cm, right=2cm, top=2.8cm, bottom=2.8cm]{geometry}

\allowdisplaybreaks

\usepackage{subcaption}
\usepackage[colorlinks=true,citecolor=red,linkcolor=blue,urlcolor=cyan,filecolor=magenta]{hyperref}
\numberwithin{equation}{section}

\usepackage{lmodern}			
\usepackage[T1]{fontenc}		
\usepackage[utf8]{inputenc}		
\usepackage{indentfirst}		
\usepackage{color,xcolor}		
\usepackage{graphicx}			
\usepackage{microtype} 			
 \usepackage{amsbsy}			
 \usepackage{amssymb}			
\usepackage{amsthm}			
 \usepackage{mathrsfs} 
 \usepackage{lineno}

 \usepackage{tikz}
\usepackage{tikz-cd}

 \usepackage{comment}
\usetikzlibrary{trees}

\newtheorem{theorem}{Theorem}[section]
\newtheorem{corollary}[theorem]{Corollary}
\newtheorem{lemma}[theorem]{Lemma}
\newtheorem{proposition}[theorem]{Proposition}
\newtheorem{remark}[theorem]{Remark}
\newtheorem{definition}[theorem]{Definition}

\newtheorem{question}{Question}

\title[Keller--Segel model on moving domains] 
{A Keller--Segel Chemotaxis Model on One-Dimensional Moving Domains}
\author[J. C. Huaccha Neyra]{Jackeline C. Huaccha Neyra}
\address{Instituto de Ci\^encia e Tecnologia, Universidade Federal de S\~ao Paulo,
Av. Cesare Mansueto Giulio Lattes, 1201, 12247-014 S\~ao Jos\'e dos Campos, SP, Brazil}
\email{jackeline.neyra@unifesp.br}
\thanks{Corresponding author: Jackeline C. Huaccha Neyra (jackeline.neyra@unifesp.br).}

\author[J. Matos Ascona]{John Matos Ascona}
\address{Universidade Estadual de Mato Grosso do Sul, Rodovia Dourados-Itahum (Rodovia MS/162), km 12, Caixa Postal 351, Cidade Universitária de Dourados, 79804-970, Dourados, MS, Brazil}
\email{john.ascona@uems.br}

\author[L. Santa Cruz Calder\'on]{Lizet Santa Cruz Calder\'on}
\address{Universidade Federal de Mato Grosso do Sul, Rua Oscar Trindade de Barros, 740, Bairro da Serraria, 79200-000, Aquidauana, MS, Brazil}
\email{lizet.calderon@ufms.br}

\author[H. L\'opez-L\'azaro]{Heraclio L\'opez-L\'azaro}
\address{Departamento de Matem\'atica, Universidade Federal de Santa Catarina,
Campus Universit\'ario Reitor Jo\~ao David Ferreira Lima, Trindade, 88040-900 Florian\'opolis, SC, Brazil}
\email{heraclio.lopez.lazaro@ufsc.br}

\begin{document}
\allowdisplaybreaks 

\begin{abstract}
We study a one-dimensional Keller--Segel system on an interval whose boundary moves according to a prescribed contracting law. The boundary motion makes the problem genuinely nonautonomous and destroys the mass conservation available on fixed domains. A time-dependent change of variables, combined with a time reparametrization, reduces the problem to a semilinear system on a fixed interval with time-independent principal part and nonautonomous transport terms. We prove global well-posedness and regularity of weak solutions, and we show that the contraction of the domain yields a uniform bound on the total mass. For the long-time dynamics, our main tool is a singular Gronwall type inequality with an exponentially damped kernel, for which we provide a self-contained version with explicit decay rate and constants, tailored to our setting; under a small gain condition on the boundary velocities, it yields dissipativity with an explicit exponential decay rate, which cannot be obtained from the classical inequality of Henry. With it we prove the existence of the minimal pullback attractor in the universe of families of bounded sets, on the moving domain and then on the fixed rescaled reference domains. Finally, numerical simulations on contracting domains illustrate the effect of the boundary motion on cell aggregation.

\vskip .1in
\noindent {\it 2020 Mathematics Subject Classification.} Primary 35B41, 35K57, 35R37, 37L30; Secondary 35Q92, 92C17, 35B40, 35A01.

\noindent {\it Key words and phrases:} Keller--Segel system, chemotaxis, moving domain, non-autonomous evolution process, pullback attractor, global well-posedness.
\end{abstract}

\maketitle

\tableofcontents

\section{Introduction} Chemotaxis is the process by which organisms or cell populations orient their movement in response to gradients of a chemical signal present in the
environment \cite{Eisenbach2004}. This mechanism plays a central role in
embryonic development \cite{DormannWeijer2003}, in wound healing and the
immune response \cite{GillitzerGoebeler2001}, in the bacterial colonization of injury sites \cite{ZhouSzymanskiBaylink2022}, and in tumor invasion and
metastasis \cite{RoussosCondeelisPatsialou2011}. In bacteria, the perception
of the stimulus is converted into changes in motility through signal
transduction pathways that typically regulate flagellar rotation
\cite{Karmakar2021,SourjikWingreen2012}. In mathematical
biology, chemotaxis has become a paradigm for the study of self-organization, collective migration, and pattern formation
\cite{BellomoTaoWinkler2015,HillenPainter2009}.

Among the models proposed, the Keller--Segel system, introduced in the early
1970s to describe the aggregation of the amoebae \textit{Dictyostelium
discoideum} \cite{KellerSegel1971a}, remains the most
influential formulation. It captures the competition between the random
diffusion of cells and their directed movement in response to a chemical
signal that they themselves produce, a competition that may lead to
aggregation, pattern formation, and singularities \cite{CieslakWinkler2008}. In the one-dimensional
parabolic--parabolic form adopted by Osaki and Yagi \cite{OY01}, the system on
a fixed, bounded, nondegenerate interval $I=(\alpha,\beta)\subset\mathbb{R}$
reads
\begin{equation}\tag{KS}\label{eq:KS}
\begin{cases}
\dfrac{\partial u}{\partial t}
   = a\,\dfrac{\partial^{2} u}{\partial x^{2}}
   -\kappa \dfrac{\partial}{\partial x}\!\left(u\,\dfrac{\partial \rho}{\partial x}\right),
   & (x,t)\in I\times(0,\infty),\\[3mm]
\dfrac{\partial \rho}{\partial t}
   = b\,\dfrac{\partial^{2} \rho}{\partial x^{2}} + c\,u - d\,\rho,
   & (x,t)\in I\times(0,\infty),\\[3mm]
\dfrac{\partial u}{\partial x}
   = \dfrac{\partial \rho}{\partial x} = 0,
   & x\in\{\alpha,\beta\},\ t\in(0,\infty),\\[3mm]
u(x,0)=u_{0}(x),\qquad \rho(x,0)=\rho_{0}(x), & x\in I,
\end{cases}
\end{equation}
where $a,b,c,d,\kappa>0$ are constants. The unknowns $u=u(x,t)\ge 0$ and
$\rho=\rho(x,t)\ge 0$ represent the cell density and the concentration of the
chemical signal, respectively. The term $-\kappa\partial_x(u\,\partial_x\rho)$
describes chemotactic transport, which drives the cells in the direction of
increasing signal; $c\,u$ models the production of the signal by the cells,
and $-d\,\rho$ its degradation. In general, the chemotactic flux takes the
form $u\,\partial_x\chi(\rho)$, with chemotactic potential $\chi$ and
sensitivity function $\chi'$. Besides the linear case $\chi(\rho)=\kappa\rho$
\cite{KellerSegel1971a}, notable examples are $\chi(\rho)=\kappa\log\rho$, associated
with the Weber--Fechner law \cite{KellerSegel1971b}, and
$\chi(\rho)=\kappa\rho/(K+\rho)$, corresponding to the so-called receptor law,
$\chi'(\rho)=\kappa K/(K+\rho)^2$, which accounts for receptor saturation at
high concentrations \cite{LapidusSchiller1976}. General conditions on $\chi$
covering these cases were analyzed in \cite{Schaaf85}; comprehensive surveys
can be found in \cite{HillenPainter2009,Horstmann2004}. In this
work we consider exclusively the linear sensitivity $\chi(\rho)=\kappa \rho$, which
is the most widely studied and the one for which the theory in the autonomous
case is available.

For fixed domains, the theory of \eqref{eq:KS} is by now well developed.
Local existence and uniqueness, as well as the behavior of the norms of
solutions, were established by Yagi \cite{Yagi97}. In higher dimensions, the
dichotomy between global existence and blow-up depends strongly on the
dimension and on the total mass e.g. \cite{DavilaDelPinoDolbeaultMussoWei2019}: in two dimensions, global existence holds below a mass threshold \cite{NSY97}, which is optimal and equal to $8\pi$ in the parabolic--elliptic case in $\mathbb{R}^2$ \cite{BDP2006}, whereas above it chemotactic collapse occurs, with or without symmetry \cite{HW01,NSS}; for developments in dimension three or higher, see \cite{BuseghinDavilaDelPinoMusso2026, Winkler2013} and the surveys \cite{BellomoTaoWinkler2015}. The asymptotic behavior of global solutions was studied in \cite{GZ98}, the set of
stationary solutions in \cite{Schaaf85,Wiebers85}, and an optimal control problem in \cite{ryu2001}.

In one dimension there is no collapse: solutions of \eqref{eq:KS} are global
and bounded for every admissible initial datum \cite{OY01,Yagi97}. It is
precisely this unconditional global well-posedness, unavailable in dimension
$n \ge 2$ for the same sensitivity, that makes the long-time dynamics a
well-defined object of study. Osaki and Yagi \cite{OY01} formulated
\eqref{eq:KS} as a semilinear evolution equation in a Hilbert space, verified
a squeezing property of the associated semigroup, and constructed an
exponential attractor in the sense of Eden, Foias, Nicolaenko, and Temam
\cite{EFNT94}, that is, a compact set of finite fractal dimension that
attracts every trajectory at an exponential rate; analogous results were
obtained for chemotaxis systems with growth \cite{OTYM2002}. One aspect of
that analysis is essential here: since the mass
$\|u(t)\|_{L^1}=\|u_0\|_{L^1}$ is conserved and each mass level contains a
spatially homogeneous equilibrium, there is no compact global attractor on the whole phase space, and the analysis is carried out on the
invariant sets of fixed mass.

This entire theory presupposes that the cells occupy a fixed region. In many
situations, however, the medium in which the population moves varies in
time, and it is convenient to distinguish between two classes of problems. In
\emph{free boundary} problems, the motion of the domain is driven by the
population itself, as in the expansion of a colony or a biofilm, and the
position of the boundary is an unknown coupled to the solution through an
additional law; for chemotaxis, existence of solutions has been established
for parabolic and parabolic--elliptic formulations
\cite{ChenLvWu2018, Wu2010}, mostly in radially
symmetric domains, where the problem reduces to a scalar equation for the
radius. In \emph{prescribed boundary} problems, the motion of the domain is
imposed by a mechanism external to the modeled population. This is the case
in embryonic development, where the tissue grows through its own
proliferation, while the cells migrating on it by chemotaxis do not determine
this growth: in this context, domain growth selects and stabilizes spatial
patterns \cite{CrampinGaffneyMaini1999,MadzvamuseGaffneyMaini2010}, and
chemotaxis models on growing domains reproduce, for instance, the formation
of stripes on the skin of growing fish \cite{PainterMainiOthmer1999}; see also
\cite{ESV2012}. The prescribed formulation is also natural when the boundary
evolves on a slow time scale relative to the cell dynamics, or when the
geometry is externally controlled, as in experimental devices with moving
walls. The present work belongs to this second class. In the works cited
above, the focus is on pattern formation and stability; to the best of our
knowledge, the long-time dynamics of Keller--Segel systems with prescribed
boundary, in the sense of pullback attractors, has not been addressed in the
literature.

More precisely, let $\alpha,\beta:\mathbb{R}\to\mathbb{R}$ with
$\alpha(t)<\beta(t)$, and set $\gamma(t):=\beta(t)-\alpha(t)$,
$\mathcal I_t:=(\alpha(t),\beta(t))$, and
$Q_{\tau,T}:=\bigcup_{t\in(\tau,T)}\mathcal I_t\times\{t\}$. We study the
system
\begin{equation}\tag{I-KS}\label{equ4}
\begin{cases}
\dfrac{\partial u}{\partial t}
 = a\,\dfrac{\partial^{2} u}{\partial x^{2}}
 - \kappa\,\dfrac{\partial}{\partial x}\!\left(u\,\dfrac{\partial \rho}{\partial x}\right),
 & (x,t)\in Q_{\tau,\infty},\\[3mm]
\dfrac{\partial \rho}{\partial t}
 = b\,\dfrac{\partial^{2} \rho}{\partial x^{2}} + c\,u - d\,\rho,
 & (x,t)\in Q_{\tau,\infty},\\[3mm]
\dfrac{\partial u}{\partial x}
 = \dfrac{\partial \rho}{\partial x} = 0,
 & x\in\{\alpha(t),\beta(t)\},\ t>\tau,\\[3mm]
u(x,\tau)=u_\tau(x),\qquad \rho(x,\tau)=\rho_\tau(x), & x\in\mathcal I_\tau,
\end{cases}
\end{equation}
where $a,b,c,d,\kappa>0$ are constants and the initial time
$\tau\in\mathbb{R}$ is arbitrary, under the following hypotheses:
\begin{enumerate}
    \item[{\bf (H1)}] $\alpha,\beta\in C^{1}(\mathbb{R})$ and  $0<\gamma_{0}<\gamma(t)<\gamma_{1}$ for all $t\in\mathbb{R}$.
    \item[{\bf (H2)}] $\alpha'(t)\geq 0$ and $\beta'(t)\leq 0$ for all $t\in \mathbb{R}$.
\end{enumerate}
For the asymptotic analysis in the pullback sense the following hypothesis will be necessary:
\begin{enumerate}
\item [{\bf (H3)}] $\alpha',\beta' \in L^{\infty}(\mathbb{R})$.
\end{enumerate}
And, for the exponential decay estimate underlying the dissipativity of the process (see
Lemma \ref{lemmaR11}), the following small gain condition will be necessary:
\begin{enumerate}
\item [{\bf(H4)}] For some specific positive constant $C>0$ and some $\theta \in (0,1)$ $$C \, \gamma_{1}\, \big(\|\alpha'\|_{L^{\infty}(\mathbb{R})} + \|\gamma'\|_{L^{\infty}(\mathbb{R})} \big)\,\Gamma(1-\theta)<1$$ where $\displaystyle\Gamma(1-\theta)=\int_0^\infty\sigma^{-\theta}e^{-\sigma}\,d\sigma$ is the Euler gamma function.
\end{enumerate}

The mobility of the boundary changes the structure of the problem in an
essential way. Integrating the first equation of \eqref{equ4} and using the homogeneous Neumann boundary conditions, one obtains the balance law
\begin{equation}\label{eq:balance-n}
\frac{d}{dt}\int_{\alpha(t)}^{\beta(t)} u(x,t)\,dx
 = \beta'(t)\,u(\beta(t),t) - \alpha'(t)\,u(\alpha(t),t).
\end{equation}
We impose hypothesis (H2) to ensure that the mass is conserved and that the level sets underlying the analysis of \cite{OY01} remain invariant; see Section \ref{sec3}. Moreover, the map
\begin{equation}\label{kappa}
r(x,t):=\frac{x-\alpha(t)}{\gamma(t)},
\end{equation}
which carries $\mathcal I_t$ onto the reference interval
$\mathcal I:=(0,1)$, converts \eqref{equ4} into a system with diffusion
coefficients proportional to $1/\gamma^2(t)$ and transport terms proportional
to $\alpha'$ and $\gamma'$, which makes the problem genuinely nonautonomous,
even though the constants $a,b,c,d,\kappa$ are fixed. Its asymptotic dynamics
is then described not by a global attractor, but by a \emph{pullback}
attractor, that is, a time-dependent family of compact sets
\cite{CLR2006,KloedenRasmussen2011,CLR}; for the alternative approach via
uniform attractors, see \cite{ChepyzhovVishik2002}. On noncylindrical
domains, this theory was developed and applied to the semilinear heat
equation in \cite{PPJ-2008,KLOEDEN20094702}. In our case, however, problem
\eqref{equ4} presents four additional difficulties. First, the cross-diffusion
term $\kappa\,\partial_x(u\,\partial_x\rho)$ requires $\rho$ to have one more
derivative of regularity than $u$, which forces us to work in product phase
spaces of the type $L^2\times H^1$, as in \cite{OY01}; consequently, obtaining
regularity requires rather delicate uniform estimates with respect to the
initial data. Second, after a first change of variables leading to an
equivalent problem on a fixed domain, the time reparametrization that removes
the nonautonomous factor $1/\gamma^2(t)$ from the principal part transfers the
time dependence to the nonlinearity, which then contains transport and
reaction terms with nonautonomous coefficients. Third, by
\eqref{eq:balance-n}, the mass is neither conserved nor dissipated, so that
neither the restriction to mass levels used in \cite{OY01} nor a standard
dissipative estimate is available. To overcome this difficulty, we adopt a
contraction hypothesis on the moving domain, namely (H2), which requires that each endpoint of $\mathcal I_t$ moves inwards; in particular, the length $\gamma$ of the domain is nonincreasing.
Fourth, the transport terms generated by the boundary motion must be controlled. On bounded time intervals this is done by means of the classical singular Gronwall inequality of Henry \cite{Henry1981}, which, under (H1)--(H2) alone, shows that the solutions do not blow up in finite time. To obtain the dissipativity of the associated evolution processes, with estimates uniform in time, we need in addition the boundedness (H3) and the small gain condition (H4) on the boundary velocities, which allow us to apply a singular Gronwall-type inequality with an exponentially damped kernel (see Lemma \ref{lem:main}). Inequalities of
this type are closely related to the classical theory of Volterra integral equations with $L^1$ kernels in exponentially weighted spaces (see \cite{GLS1990}). Since our arguments require an explicit decay rate and
explicit constants, we provide a self-contained version tailored to our setting, which combines the iteration technique of Henry \cite{Henry1981} with resolvent arguments in the spirit of Gripenberg, Londen, and Staffans
\cite{GLS1990}. Having overcome these four difficulties,
we obtain the existence of a pullback attractor in the universe of bounded sets by applying the recent results of \cite{HRHJ-2024}, since the analysis
is carried out on time-dependent families of phase spaces
$\mathbb{H}^M_{1,t}$, indexed by a parameter $M>0$ that controls the mass. In
turn, the control of the mass along the flow follows from the control of the growth of the moving domain established in (H2).

Our main contributions are the following.
\begin{enumerate}
\item By combining the change of variables \eqref{kappa} with a time
reparametrization, we reduce \eqref{equ4} to an abstract semilinear problem
on the fixed interval $\mathcal I$ whose principal part is time-independent.
For this problem, we prove the existence, uniqueness, and regularity of weak
solutions, and we show that they define continuous evolution processes.

\item We show that the contraction of the domain yields a uniform bound on the total mass, and we establish a priori estimates uniform in time by means of a singular Gronwall inequality with an exponentially damped kernel (see Lemma \ref{lem:main}), for which we give a self-contained
version with explicit decay rate and constants. Under a small-gain condition, this inequality yields an explicit exponential decay rate, which cannot be obtained from the classical version of Henry \cite{Henry1981}, whose kernel is not integrable on $(0,\infty)$.

\item We prove the existence of the minimal pullback attractor in the universe of families of bounded sets, first for the problem on the reference domain and then by transferring this result to the moving domain.

\item We complement the analysis with numerical simulations on contracting domains, illustrating the effect of the boundary motion on cell aggregation.
\end{enumerate}

The paper is organized as follows. In Section \ref{sec2} we perform the change of variables and the time reparametrization, and we collect the results on time-dependent Bochner spaces, nonautonomous dynamical systems, and the singular Gronwall inequality. In Section \ref{sec3} we prove the well-posedness and regularity of weak solutions. In Section \ref{sec4} we construct the pullback attractor. In Section \ref{sec5} we present the numerical simulations.
\section{Preliminaries}\label{sec2}

\subsection{Change of variables}
 
To transform \eqref{equ4}, posed on the moving domain $\mathcal I_t$, into an equivalent problem on the fixed domain $\mathcal I$, we use the map $r(x,t)$ from \eqref{kappa} and look for a solution of the form
\begin{equation}\label{recuperar}
u(x,t)=v\big(r(x,t),t\big),\qquad
\rho(x,t)=\psi\big(r(x,t),t\big),\qquad
v,\psi\in L^1_{loc}\big(\mathcal I\times(\tau,T)\big).
\end{equation}
 
Writing $y:=r(x,t)$, the chain rule gives
\[
\frac{\partial u}{\partial t}=\frac{\partial v}{\partial y}\frac{\partial r}{\partial t}+\frac{\partial v}{\partial t},
\qquad
\frac{\partial u}{\partial x}=\frac{\partial v}{\partial y}\frac{\partial r}{\partial x}=\frac1{\gamma(t)}\frac{\partial v}{\partial y}.
\]
A direct computation using $x=r^{-1}(y,t)=y\gamma(t)+\alpha(t)$ gives
\[
\frac{\partial r}{\partial t}\big(r^{-1}(y,t),t\big)
=\frac{-\alpha'(t)\gamma(t)-\big[y\gamma(t)+\alpha(t)\big]\gamma'(t)+\alpha(t)\gamma'(t)}{\gamma^2(t)}
=-\frac{\alpha'(t)+y\,\gamma'(t)}{\gamma(t)}=:b(y,t),
\]
so that
\begin{equation*}
\frac{\partial u}{\partial t}=\frac{\partial v}{\partial t}+b(y,t)\frac{\partial v}{\partial y} .
\end{equation*}
Differentiating $\partial u/\partial x$ once more in $x$,
\begin{equation*}
\frac{\partial^2u}{\partial x^2}=\frac1{\gamma^2(t)}\frac{\partial^2v}{\partial y^2} .
\end{equation*}
The same computation applied to $\rho(x,t)=\psi(r(x,t),t)$ gives
\begin{equation*}
\frac{\partial\rho}{\partial t}=\frac{\partial\psi}{\partial t}+b(y,t)\frac{\partial\psi}{\partial y},
\qquad
\frac{\partial^2\rho}{\partial x^2}=\frac1{\gamma^2(t)}\frac{\partial^2\psi}{\partial y^2} .
\end{equation*}
Finally, for the coupling term,
\[
\frac{\partial}{\partial x}\Big(u\frac{\partial\rho}{\partial x}\Big)
=\frac{\partial u}{\partial x}\frac{\partial\rho}{\partial x}+u\frac{\partial^2\rho}{\partial x^2}
=\frac1{\gamma^2(t)}\Big[\frac{\partial v}{\partial y}\frac{\partial\psi}{\partial y}+v\frac{\partial^2\psi}{\partial y^2}\Big]
=\frac1{\gamma^2(t)}\frac{\partial}{\partial y}\Big(v\frac{\partial\psi}{\partial y}\Big) .
\]
 
\paragraph{The problem on the fixed domain.}
Substituting these identities into \eqref{equ4}, the system becomes,
on $\mathcal I\times(\tau,\infty)$,
\begin{equation}\tag{II-KS}\label{equ5}
\begin{cases}
\displaystyle\frac{\partial v}{\partial t}=\frac{a}{\gamma^2(t)}\frac{\partial^2v}{\partial y^2}
-\frac{\kappa}{\gamma^2(t)}\frac{\partial}{\partial y}\Big(v\frac{\partial\psi}{\partial y}\Big)-b(y,t)\frac{\partial v}{\partial y},\quad (y,t)\in \mathcal{I}\times(\tau,\infty)\\[8pt]
\displaystyle\frac{\partial\psi}{\partial t}=\frac{b}{\gamma^2(t)}\frac{\partial^2\psi}{\partial y^2}-b(y,t)\frac{\partial\psi}{\partial y}+cv-d\psi,\hspace{1.6cm} (y,t)\in \mathcal{I}\times(\tau,\infty)\\[8pt]
\displaystyle \frac{\partial v}{\partial y}=\frac{\partial \psi}{\partial y}=0, \hspace{6cm} y\in \{0,1\},\quad \forall\,t\geq\tau,\\[4pt]
v(y,\tau)=v_\tau(y),\quad\psi(y,\tau)=\psi_\tau(y), \hspace{2.7cm} y\in\mathcal I,
\end{cases}
\end{equation}
where $v(y,t)=u(r^{-1}(y,t),t)$ and $b(y,t)=-(\alpha'(t)+\gamma'(t)y)/\gamma(t)$.
 
\begin{remark}
Here $b>0$ (without arguments) denotes the constant diffusion
coefficient in the second equation, distinct from the drift function $b(y,t)$ defined above; the two never appear together without their respective arguments, so no ambiguity arises.
\end{remark}
 
\paragraph{\bf Time rescaling.} System \eqref{equ5} still carries the time-dependent factor
$1/\gamma^2(t)$ in front of the diffusive terms. Multiplying both
equations by $\gamma^2(t)$,
\[
\gamma^2(t)\frac{\partial v}{\partial t}=a\frac{\partial^2v}{\partial y^2}
-\kappa\frac{\partial}{\partial y}\Big(v\frac{\partial\psi}{\partial y}\Big)-\gamma^2(t)b(y,t)\frac{\partial v}{\partial y},
\]
\[
\gamma^2(t)\frac{\partial\psi}{\partial t}=b\frac{\partial^2\psi}{\partial y^2}-\gamma^2(t)b(y,t)\frac{\partial\psi}{\partial y}+c\gamma^2(t)v-d\gamma^2(t)\psi,
\]
and introduce the new time variable
\begin{equation}\label{rescaling}
s(t):=\tau+\int_\tau^t\frac{dr}{\gamma^2(r)}, \qquad t\geq\tau,
\end{equation}
with inverse $t(\cdot):=s^{-1}(\cdot)$. Setting $w(y,s):=v(y,t(s))$ and
$\phi(y,s):=\psi(y,t(s))$, the chain rule
$\partial_sw=\gamma^2(t(s))\,\partial_tv$ turns the two equations above
into
\[
\frac{\partial w}{\partial s}=a\frac{\partial^2w}{\partial y^2}
-\kappa\frac{\partial}{\partial y}\Big(w\frac{\partial\phi}{\partial y}\Big)-\gamma^2(t(s))b(y,t(s))\frac{\partial w}{\partial y},
\]
\[
\frac{\partial\phi}{\partial s}=b\frac{\partial^2\phi}{\partial y^2}-\gamma^2(t(s))b(y,t(s))\frac{\partial\phi}{\partial y}+c\gamma^2(t(s))w-d\gamma^2(t(s))\phi .
\]
Writing
\[
g_1(s,w,\phi):=-\kappa\frac{\partial}{\partial y}\Big(w\frac{\partial\phi}{\partial y}\Big)-\gamma^2(t(s))b(y,t(s))\frac{\partial w}{\partial y},
\]
\[
g_2(s,w,\phi):=c\gamma^2(t(s))w-d\gamma^2(t(s))\phi-\gamma^2(t(s))b(y,t(s))\frac{\partial\phi}{\partial y},
\]
we arrive at the semilinear system with time-independent principal part
\begin{equation}\tag{III-KS}\label{equ4iii}
\begin{cases}
\displaystyle\frac{\partial w}{\partial s}-a\frac{\partial^2w}{\partial y^2}=g_1(s,w,\phi), \hspace{2cm} \text{on }\mathcal I\times(\tau,\infty),\\[8pt]
\displaystyle\frac{\partial\phi}{\partial s}-b\frac{\partial^2\phi}{\partial y^2}=g_2(s,w,\phi), \hspace{2.2cm}   \text{on }\mathcal I\times(\tau,\infty),\\[8pt]
\displaystyle \frac{\partial w}{\partial y}=\frac{\partial \phi}{\partial y}=0, \hspace{4cm} y\in \{0,1\},\quad \forall\,t\geq\tau,\\[4pt]
w(y,\tau)=w_\tau(y),\quad\phi(y,\tau)=\phi_\tau(y), \hspace{0.65cm} y\in\mathcal I,
\end{cases}
\end{equation}
 
Thanks to Propositions~\ref{Prop1} and~\ref{Prop7}, and using the
change of variables $r(\cdot,\cdot)$ from \eqref{kappa}, a solution of \eqref{equ4} on the moving domain is recovered directly from a
solution of \eqref{equ4iii}.

\subsection{Time-dependent Bochner spaces} The construction of the time-dependent spaces in \cite{KLOEDEN20094702} relies on the extension by
zero of functions in $H^1_0(\mathcal I_t)$ to the whole real line. Since in this work we deal with
homogeneous Neumann boundary conditions, this construction cannot be used directly: the extension by
zero of a function in $H^1(\mathcal I_t)$ which does not vanish at the endpoints has jumps, and thus it
does not belong to $H^1(\mathbb R)$. Instead, we define these spaces by transporting the moving interval $\mathcal I_t$ onto a fixed reference interval. This approach is in the spirit of the framework of evolving Bochner spaces developed in \cite{AlphonseElliottStinner2015}, which does not depend on the boundary conditions. Since in our one-dimensional setting the map $y\mapsto r^{-1}(y,t)=\alpha(t)+\gamma(t)y$ is affine, the results below can be proven by direct calculation.

\medskip
\subsubsection{\textbf{Sobolev spaces with Neumann boundary conditions.}}
For an open interval $J=(a,b)$ and an integer $k\geq0$, we consider the space $H^k(J)$ with the norm
$\|w\|_{H^k(J)}^2:=\sum_{j=0}^k\|w^{(j)}\|_{L^2(J)}^2$. Since $H^k(J)\subset C^{k-1}(\overline J)$ for
$k\geq1$, the boundary values of $w^{(j)}$, $0\leq j\leq k-1$, are well defined, and we set
\[
H^k_N(J):=\big\{w\in H^k(J)\,:\,w'(a)=w'(b)=0\big\},\qquad k\geq2,
\]
\[
H^k_{N^2}(J):=\big\{w\in H^k_N(J)\,:\,w'''(a)=w'''(b)=0\big\},\qquad k\geq4 .
\]
These are closed subspaces of $H^k(J)$, hence Hilbert spaces with the induced norm. For each $t\in\mathbb{R}$, we denote by $\|\cdot\|_{L^p(\mathcal{I}_t)}$ the norm of $L^{p}(\mathcal{I}_{t})$, by $\|\cdot\|_{L^2(\mathcal{I}_t)}$ and
$(\cdot,\cdot)_{L^2(\mathcal{I}_t)}$ the norm and the inner product of $L^{2}(\mathcal{I}_{t})$, by $\|\cdot\|_{H^1(\mathcal{I}_t)}$ the norm of $H^{1}(\mathcal{I}_{t})$, by $\|\cdot\|_{H^k(\mathcal{I}_t)}$ the norm of $H^{k}(\mathcal{I}_{t})$, and by
$\langle\cdot,\cdot\rangle_t$ the duality pairing between $H^1(\mathcal I_t)'$ and $H^1(\mathcal I_t)$.

The following lemma, which follows from the change of variables $x=r^{-1}(y,t)$, shows that the transformation to the reference interval preserves all these spaces, with equivalence constants which are independent of $t$, by (H1).

\begin{lemma}\label{lem:scaling}
Let $t\in\mathbb R$, $w\in L^1_{loc}(\mathcal I_t)$, and $\widetilde w:=w(r^{-1}(\cdot,t))$. Then:
\begin{itemize}
\item[(i)] For $p\in[1,\infty)$, $\|w\|_{L^p(\mathcal I_t)}=\gamma(t)^{1/p}\|\widetilde w\|_{L^p(\mathcal I)}$, and
$\|w\|_{L^\infty(\mathcal I_t)}=\|\widetilde w\|_{L^\infty(\mathcal I)}$.
\item[(ii)] For every integer $k\geq0$, $w\in H^k(\mathcal I_t)$ if and only if $\widetilde w\in H^k(\mathcal I)$, and in
this case $\partial_y^j\widetilde w=\gamma(t)^j\,(\partial_x^jw)\circ r^{-1}(\cdot,t)$ for $0\leq j\leq k$, so that
\begin{equation}\label{eq:scaling}
\Big\|\frac{\partial^j w}{\partial x^j}\Big\|_{L^2(\mathcal I_t)}^2=\gamma(t)^{1-2j}\,\|\partial_y^j\widetilde w\|_{L^2(\mathcal I)}^2,
\qquad 0\leq j\leq k .
\end{equation}
Consequently,
\begin{equation}\label{eq:equiv-Hk}
c_k\|\widetilde w\|_{H^k(\mathcal I)}\leq\|w\|_{H^k(\mathcal I_t)}\leq C_k\|\widetilde w\|_{H^k(\mathcal I)},
\end{equation}
where
$c_k^2:=\min_{0\leq j\leq k}\min\{\gamma_0^{1-2j},\gamma_1^{1-2j}\}$ and
$C_k^2:=\max_{0\leq j\leq k}\max\{\gamma_0^{1-2j},\gamma_1^{1-2j}\}$.
\item[(iii)] $w\in H^k_N(\mathcal I_t)$ if and only if $\widetilde w\in H^k_N(\mathcal I)$ $(k\geq2)$, and
$w\in H^k_{N^2}(\mathcal I_t)$ if and only if $\widetilde w\in H^k_{N^2}(\mathcal I)$ $(k\geq4)$.
\item[(iv)] For $F\in H^1(\mathcal I_t)'$, let $\widetilde F\in H^1(\mathcal I)'$ be defined by
$\langle\widetilde F,z\rangle:=\langle F,z(r(\cdot,t))\rangle_t$, $z\in H^1(\mathcal I)$. Then $F\mapsto\widetilde F$ is an
isomorphism from $H^1(\mathcal I_t)'$ onto $H^1(\mathcal I)'$ and
$C_1^{-1}\|\widetilde F\|_{H^1(\mathcal I)'}\leq\|F\|_{H^1(\mathcal I_t)'}\leq c_1^{-1}\|\widetilde F\|_{H^1(\mathcal I)'}$.
\end{itemize}
\end{lemma}

\begin{remark}\label{rem:dualH1}
For an open interval $J$, every $F\in H^1(J)'$ can be written as
\[
\langle F,\varphi\rangle=\int_Jf_0\,\varphi\,dx+\int_Jf_1\,\varphi'\,dx,\qquad\forall\,\varphi\in H^1(J),
\]
with $f_0,f_1\in L^2(J)$ and $\|f_0\|_{L^2(J)}^2+\|f_1\|_{L^2(J)}^2=\|F\|_{H^1(J)'}^2$. Indeed, by the Riesz
representation theorem in $H^1(J)$, there exists $g\in H^1(J)$ such that
$\langle F,\varphi\rangle=(g,\varphi)_{L^2(J)}+(g',\varphi')_{L^2(J)}$ and $\|g\|_{H^1(J)}=\|F\|_{H^1(J)'}$; it suffices to
take $f_0=g$ and $f_1=g'$. In particular, functionals involving boundary values, such as
$\varphi\mapsto\varphi(a)$, are represented by interior $L^2$-terms, and no explicit boundary contributions
appear in this representation. Note, however, that $H^1(J)'$ is strictly larger than $H^{-1}(J)$: two
elements of $H^1(J)'$ which coincide on $C^\infty_c(J)$ may differ by a combination of boundary values.
\end{remark}

\medskip
\subsubsection{\textbf{Bochner spaces on the moving interval.}}
In what follows, $X(J)$ denotes one of the spaces $L^p(J)$ $(1\leq p\leq\infty)$, $H^k(J)$, $H^k_N(J)$,
$H^k_{N^2}(J)$ or $H^1(J)'$, and $\|\cdot\|_{X(J)}$ its norm.

For each $t$, $r(\cdot,t):\mathcal I_t\to\mathcal I$ is an affine bijection with inverse $r^{-1}(\cdot,t)$, and the map $\Phi(y,t):=(r^{-1}(y,t),t)$ is a $C^1$-diffeomorphism from $\mathcal I\times(\tau,T)$ onto
$Q_{\tau,T}$ whose Jacobian determinant equals $\gamma(t)\in[\gamma_0,\gamma_1]$.

\begin{definition}\label{def:Bochner-t}
Let $\tau<T$ and $q\in[1,\infty]$.
\begin{itemize}
\item Let $u\in L^1_{loc}(Q_{\tau,T})$ and $v=u\circ\Phi$. We say that $u\in L^q(\tau,T;X(\mathcal I_t))$
(respectively, $u\in C([\tau,T];X(\mathcal I_t))$) if $v\in L^q(\tau,T;X(\mathcal I))$ (respectively,
$v\in C([\tau,T];X(\mathcal I))$), in the usual Bochner sense, and we set
\[
\|u\|_{L^q(\tau,T;X(\mathcal I_t))}^q:=\int_\tau^T\|u(t)\|_{X(\mathcal I_t)}^q\,dt
\quad\text{if } q<\infty,
\qquad
\|u\|_{L^\infty(\tau,T;X(\mathcal I_t))}:=\operatorname*{ess\,sup}_{t\in(\tau,T)}\|u(t)\|_{X(\mathcal I_t)} .
\]
We say that a sequence $\{u^m\}$ converges to $u$ in $C([\tau,T];X(\mathcal I_t))$ if
$\{u^m\circ\Phi\}$ converges to $v$ in $C([\tau,T];X(\mathcal I))$.
\item Let $F=\{F(t)\}_{t\in(\tau,T)}$ with $F(t)\in H^1(\mathcal I_t)'$ for a.e. $t$. We say that
$F\in L^q(\tau,T;H^1(\mathcal I_t)')$ if $\widetilde F:=\{\widetilde{F(t)}\}_{t\in(\tau,T)}\in L^q(\tau,T;H^1(\mathcal I)')$,
with $\widetilde{F(t)}$ given by Lemma~\ref{lem:scaling}~(iv), and we define its norm as above.
\end{itemize}
\end{definition}

\begin{proposition}\label{Prop1}
Let $v\in L^{1}_{loc}(\mathcal{I}\times(\tau,T))$ and $u(x,t)=v\big(r(x,t),t\big)$. Then, for
$q\in[1,\infty]$ and $X$ as above, $u\in L^q(\tau,T;X(\mathcal I_t))$ if and only if
$v\in L^q(\tau,T;X(\mathcal I))$. In this case the function $t\mapsto\|u(t)\|_{X(\mathcal I_t)}$ is
measurable, and
\[
c\,\|v\|_{L^q(\tau,T;X(\mathcal I))}\leq\|u\|_{L^q(\tau,T;X(\mathcal I_t))}\leq C\,\|v\|_{L^q(\tau,T;X(\mathcal I))},
\]
where $(c,C)=(\gamma_0^{1/p},\gamma_1^{1/p})$ if $X=L^p$, $p<\infty$, $(c,C)=(1,1)$ if $X=L^\infty$,
$(c,C)=(c_k,C_k)$ if $X\in\{H^k,H^k_N,H^k_{N^2}\}$, and $(c,C)=(C_1^{-1},c_1^{-1})$ if $X=H^1(\cdot)'$. In
particular,
\[
u\in L^{2}(\tau,T;H_N^{k}(\mathcal{I}_{t}))\iff v\in L^{2}(\tau,T;H^{k}_N(\mathcal{I})),
\]
\[
u\in L^{2}(\tau,T;H^{k}_{N^2}(\mathcal{I}_{t}))\iff v\in L^{2}(\tau,T;H_{N^2}^{k}(\mathcal{I})).
\]
The same equivalence and the same bounds (with $q=\infty$) hold for the spaces $C([\tau,T];X(\mathcal I_t))$.
\end{proposition}

\medskip
\subsubsection{\textbf{Time derivative.}}
On the reference interval we use the Gelfand triple
$H^1(\mathcal I)\subset L^2(\mathcal I)\equiv L^2(\mathcal I)'\subset H^1(\mathcal I)'$, with dense and continuous
inclusions. For $v\in L^2(\tau,T;H^1(\mathcal I))$, we write $\partial_tv\in L^2(\tau,T;H^1(\mathcal I)')$ if
there exists $g\in L^2(\tau,T;H^1(\mathcal I)')$ such that
\[
\int_\tau^T(v(t),z)_{L^2(\mathcal I)}\,\theta'(t)\,dt=-\int_\tau^T\langle g(t),z\rangle\,\theta(t)\,dt,
\qquad\forall\,z\in H^1(\mathcal I),\ \theta\in C^\infty_c(\tau,T),
\]
and in this case we set $\partial_tv:=g$. We denote
\[
W(\tau,T;\mathcal I):=\big\{v\in L^2(\tau,T;H^1(\mathcal I))\,:\,\partial_tv\in L^2(\tau,T;H^1(\mathcal I)')\big\}.
\]
It is classical that $W(\tau,T;\mathcal I)\subset C([\tau,T];L^2(\mathcal I))$, and that
for every $v\in W(\tau,T;\mathcal I)$ the function $t\mapsto\|v(t)\|_{L^2(\mathcal I)}^2$ is absolutely continuous
with $\frac{d}{dt}\|v(t)\|_{L^2(\mathcal I)}^2=2\langle\partial_tv(t),v(t)\rangle$ for a.e. $t\in(\tau,T)$.

\begin{definition}\label{def:W-t}
We set $W(\tau,T;\mathcal I_t):=\{u\in L^1_{loc}(Q_{\tau,T})\,:\,v=u\circ\Phi\in W(\tau,T;\mathcal I)\}$. For
$u\in W(\tau,T;\mathcal I_t)$ and a.e. $t\in(\tau,T)$, we define $\frac{\partial u}{\partial t}(t)\in H^1(\mathcal I_t)'$ by
\begin{equation}\label{eq:def-ut}
\Big\langle\frac{\partial u}{\partial t}(t),w\Big\rangle_t
:=\gamma(t)\,\big\langle\partial_tv(t),\widetilde w\big\rangle
-\int_{\mathcal I}\frac{\partial r^{-1}}{\partial t}(y,t)\,\partial_yv(y,t)\,\widetilde w(y)\,dy,
\qquad w\in H^1(\mathcal I_t),
\end{equation}
where $\widetilde w=w(r^{-1}(\cdot,t))$.
\end{definition}

The next result justifies this definition and contains the rule for differentiating
$\|u(t)\|_{L^2(\mathcal I_t)}^2$, which is used throughout this work.

\begin{proposition}\label{Prop7}
Let $u\in L^1_{loc}(Q_{\tau,T})$ and $v=u\circ\Phi$. Then:
\begin{itemize}
\item[(i)] If $u\in W(\tau,T;\mathcal I_t)$, then $u\in L^2(\tau,T;H^1(\mathcal I_t))$,
$\frac{\partial u}{\partial t}\in L^2(\tau,T;H^1(\mathcal I_t)')$, and
\[
\int_\tau^T\Big\|\frac{\partial u}{\partial t}(t)\Big\|_{H^1(\mathcal I_t)'}^2dt
\leq\frac{2}{c_1^2}\int_\tau^T\Big(\gamma_1^2\|\partial_tv(t)\|_{H^1(\mathcal I)'}^2+L^2\|v(t)\|_{H^1(\mathcal I)}^2\Big)dt ,
\]
where $L:=\sup_{t\in\mathbb R}\max\{|\alpha'(t)|,|\beta'(t)|\}$, which is finite by (H3).
\item[(ii)] For $u\in W(\tau,T;\mathcal I_t)$, $\frac{\partial u}{\partial t}$ coincides with the distributional time
derivative of $u$ in $Q_{\tau,T}$, that is,
\[
-\int_\tau^T\int_{\alpha(t)}^{\beta(t)}u(x,t)\,\frac{\partial\varphi}{\partial t}(x,t)\,dx\,dt
=\int_\tau^T\Big\langle\frac{\partial u}{\partial t}(t),\varphi(\cdot,t)\Big\rangle_t\,dt,
\qquad\forall\,\varphi\in C^\infty_c(Q_{\tau,T}).
\]
If, in addition, $v\in H^1(\tau,T;L^2(\mathcal I))$, then $\frac{\partial u}{\partial t}\in L^2(Q_{\tau,T})$ and
$\big\langle\frac{\partial u}{\partial t}(t),w\big\rangle_t=\int_{\alpha(t)}^{\beta(t)}\frac{\partial u}{\partial t}(x,t)\,w(x)\,dx$ for all
$w\in H^1(\mathcal I_t)$.
\item[(iii)] For $u\in W(\tau,T;\mathcal I_t)$, we have $u\in C([\tau,T];L^2(\mathcal I_t))$, the function
$t\mapsto\|u(t)\|_{L^2(\mathcal I_t)}^2$ is absolutely continuous, and
\begin{equation}\label{eq:Reynolds}
\frac{d}{dt}\|u(t)\|_{L^2(\mathcal I_t)}^2=2\Big\langle\frac{\partial u}{\partial t}(t),u(t)\Big\rangle_t
+\beta'(t)\,|u(\beta(t),t)|^2-\alpha'(t)\,|u(\alpha(t),t)|^2
\qquad\text{for a.e. } t\in(\tau,T).
\end{equation}
\end{itemize}
\end{proposition}

\begin{remark}\label{rem:Reynolds-higher}
Proposition~\ref{Prop7}~(iii) applies, in particular, to spatial derivatives of sufficiently regular
solutions. For instance, if $\frac{\partial^2\rho}{\partial x^2}\in W(\tau,T;\mathcal I_t)$, then \eqref{eq:Reynolds} with
$u=\frac{\partial^2\rho}{\partial x^2}$ is the identity used in the proof of Lemma~\ref{AsympCL}.
\end{remark}

\subsection{Theory of non-Autonomous Dynamical Systems }

In this section we will focus and recall some concepts known in the theory of pullback attractors, aiming to present the existence and relationship of the families of pullback attractors associated to closed processes defined on families of phase spaces parameterized in time cf. \cite{ HRHJ-2024, JHOC-2024, KLOEDEN20094702,PPJ-2008}.

Denote by $\mathbb{R}^{2}_{d}=\{(t,\tau)\in\mathbb{R}^{2}:\tau\leq t\}$. Let us consider a family of metric spaces $\{(X_{t},d_{X_{t}})\}_{t\in\mathbb{R}},$ and we remember the notion of process.
\begin{definition}
A bi-parametric family of operators  $\{U(t,\tau):X_{\tau}\rightarrow X_{t}\;|\;(\tau,t)\in\mathbb{R}_{d}^{2}\}$ is a process on $\{X_{t}\}_{t\in\mathbb{R}}$ if \begin{itemize}
\item[(i)] $U(\tau,\tau)=Id_{X_{\tau}}$ (identity on $X_{\tau}$) for any $\tau\in\mathbb R$,

\item[(ii)] $U(t,s)U(s,\tau )=U(t,\tau)$ for any $\tau\le s\le t$.
\end{itemize}
\end{definition}

For convenience we will refer to the process by $U$, instead of $\{U(t,\tau):X_{\tau}\rightarrow X_{t}\;|\;(\tau,t)\in\mathbb{R}_{d}^{2}\}$. Moreover, remember that if the evolution process $U$ satisfies $U(t,\tau):(X_{\tau},d_{X_{\tau}})\rightarrow (X_{t},d_{X_{t}})$ is a continuous map for any $(t,\tau)\in \mathbb{R}^{2}_{d}$, then, we say that $U$ is a continuous process. 

As in the random attractor theory cf. \cite[Ch. 1]{CLR}, let us denote by  $\mathcal{D}_{X}$ a universe on the family of phase spaces $\{X_{t}\}_{t\in\mathbb{R}}$, that is, a nonempty class of families parameterized in time $\widehat{D}=\{D(t): D(t)\subset X_{t},\;D(t)\neq \emptyset,\;t\in\mathbb{R}\}$.

\begin{definition}
  A universe $\mathcal{D}_X$ is inclusion-closed if given $\widehat{D}\in{\mathcal{D}_X}$ and $\widehat{D}'=\{D'(t)\subset X_{t}:t\in\mathbb{R}\}$ with $D'(t)\subset D(t)$ for all $t\in\mathbb{R}$, it fulfils that $\widehat{D}'\in {\mathcal{D}_X}$.
\end{definition}

\begin{remark} We are going to denote the universe of families of bounded sets in $\{X_t\}_{t\in \mathbb{R}}$ by $\mathcal{D}_{F}^{X}$, that is, $\mathcal{D}_{F}^{X} :=\{  \widehat D=\{D(t) \subset B_{X_t}(0,R):t\in\mathbb R\} \text{ for some } R >0 \}$. 
\end{remark}

Now we are going to give two fundamental definitions in the theory of pullback attractors, which are the key to obtaining objects of attraction, which are absorption and asymptotic compactness.

\begin{definition}
A family $\widehat D_0=\{D_0(t)\subset X_{t}:t\in\mathbb R\}$ is pullback $\mathcal D_X$-absorbing for a process $U$ if for any $t\in\mathbb R$ and $\widehat D\in\mathcal {D}_X$ there exists $\tau (\widehat D,t)\le t$ such that $U(t,\tau)D(\tau)\subset D_0(t)$ for all $\tau\le\tau (\widehat D,t).$
\end{definition}

\begin{definition}
Given a family $\widehat D_0=\{D_0(t):t\in\mathbb R\}$, a process $U$ is pullback $\widehat D_0$-asymptotically compact if for any $t\in\mathbb{R}$ and sequences $\{\tau_n\}\subset (-\infty ,t]$ and $\{x_{n}\}_{n\in\mathbb{N}}$ with $\tau_n\to-\infty$ and $x_n\in D_0(\tau_n)$ for all $n,$ it holds that the sequence $\{U(t,\tau_n)x_n\}$ is relatively compact in $X_{t}$.

A process $U$ is said pullback $\mathcal D_X$-asymptotically compact if it is pullback $\widehat D$-asymptotically compact for any $\widehat D\in\mathcal D_X.$
\end{definition}


\begin{definition}\label{Def-PA}
A family $\widehat{\mathcal{A}}=\{\mathcal{A}(t):\mathcal{A}(t)\subset X_{t},\;\mathcal{A}(t)\neq \emptyset,\;t\in\mathbb{R}\}$ is said to be a pullback $\mathcal{D}_X$-attractor for the process $U$, if:
\begin{itemize}
    \item [(i)] $\mathcal{A}(t)$ is compact in $X_{t}$ for all $t\in\mathbb{R}$;
    \item [(ii)] $\widehat{\mathcal{A}}$ is pullback $\mathcal{D}_X$-attracting, i.e., $$\lim_{\tau\rightarrow-\infty}dist_{X_{t}}(U(t,\tau)D(\tau),\mathcal{A}(t))=0\quad\text{ for all } \widehat{D}\in\mathcal{D}_X,\text{ and for all }t\in \mathbb{R},$$
     where for each $t\in\mathbb{R}$, $dist_{X_{t}}$ denotes the Hausdorff semi-distance in $X_{t}$ between two subsets of $X_{t}$;
    \item[(iii)] $\widehat{\mathcal{A}}$ is invariant, i.e., $$U(t,\tau)\mathcal{A}(\tau)=\mathcal{A}(t)\quad\text{ for all }(t,\tau)\in\mathbb{R}^{2}_{d}.$$
\end{itemize}
A pullback $\mathcal{D}_X$-attractor $\widehat{\mathcal{A}}$ is said to be minimal if it satisfies that if there exists another family of closed sets $\widehat{C}=\{C(t):C(t)\subset X_{t},\;C(t)\neq\emptyset,\;t\in\mathbb{R}\}$ such that it is pullback $\mathcal{D}_X$-attracting, then $\mathcal{A}(t)\subset C(t)$ for all $t\in \mathbb{R}$.
\end{definition}
\begin{remark}
Observe that the conditions of the pullback attractor $(i)$-$(iii)$ not guarantee that it be unique, cf. \cite{mrreal09}. The condition of minimality gives uniqueness. On the other hand, when the attractor $\widehat{\mathcal{A}}\in\mathcal D_X,$ then it is also the unique family of closed subsets in $\mathcal D_X$ satisfying $(ii)$-$(iii).$
\end{remark}

As usual, we build our pullback attractors as a union of omega-limit (or $\omega$-limit) sets, since these have good properties of attraction and compactness. Namely, given a family $\widehat D=\{D(t):D(t)\subset X_{t},\;D(t)\neq\emptyset,\;t\in\mathbb{ R}\},$ the omega-limit set of $\widehat D$ by $U$ at time $t$ (when it make sense) is defined by
\[
\Lambda(\widehat D,t)=\bigcap_{s\le t}\overline{\bigcup_{\tau\le s}U(t,\tau)D(\tau)}^{X_{t}}\quad\forall t\in\mathbb{R}.
\]

We are ready to state the theorem of existence of pullback attractors, whose candidate to satisfy the properties (i)-(iii) exposed in the Definition \ref{Def-PA}, is the family (parameterized in time), whose elements are the omega-limit set. The proof can be found in \cite[Theorem 23]{PPJ-2008}.

\begin{theorem}\label{t4}
Consider a continuous process $U$ on the family $\{X_{t}\}_{t\in\mathbb{R}}$, a pullback $\mathcal D$-absorbing family $\widehat D_0=\{D_{0}(t):D_0(t)\subset X_{t},\;D_{0}(t)\neq\emptyset,\; t\in\mathbb{R}\}$ and assume that $U$ is pullback $\widehat D_0$-asymptotically compact. Then the family $\mathcal {A}_{\mathcal D_{X}}=\{\mathcal{A}_{\mathcal D_X }(t)\subset X_{t}:t\in\mathbb{R}\}$ given by
$$\mathcal {A}_{\mathcal D_X}(t)=\overline{\bigcup_{\widehat D\in\mathcal D_X}\Lambda(\widehat D,t)}^{X_{t}}\quad\forall t\in\mathbb{R},$$
is the minimal pullback $\mathcal D_X$-attractor.
Moreover, it holds that $\mathcal{A}_{\mathcal D_X}(t)\subset\overline{D_0(t)}^{X_{t}}$ for any $t\in\mathbb{R}.$ If $\widehat D_0\in\mathcal D_X$ has closed sections and $\mathcal D_X$ is inclusion-closed then $\mathcal{A}_{\mathcal D_X}\in\mathcal D_X$.
\end{theorem}

The next result allows us to compare different attractors, for instance for regularity purposes.

\begin{theorem}\label{compara} (cf. \cite[Theorem 4]{DCDSB2018} or \cite[Theorem 3.12]{HRHJ-2024})
Consider two families of metric spaces $\{(X_t,d_{X_t})\}$ and $\{(Y_t,d_{Y_t})\}$ with continuous embedding $X_t\subset Y_t,$ for all $t\in \mathbb R$, respective universes $\mathcal D_X$ and  $\mathcal D_Y$, and $\mathcal D_X\subset \mathcal D_Y.$ Assume that $\mathcal U$ is a  continuous process in both family of spaces. For each $t\in\mathbb R$ denote
$$\mathcal A_X(t)=\overline{\bigcup_{\widehat D_1\in\mathcal D_{X}}\Lambda_{X_t}(\widehat D_1,t)}^{X_t},\qquad \text{ and }\qquad \mathcal A_Y(t)=\overline{\bigcup_{\widehat D_2\in\mathcal D_{Y}}\Lambda_{Y_t}(\widehat D_2,t)}^{Y_t}.$$
 Then $\mathcal A_{X}(t)\subset\mathcal A_{Y}(t)$ for all $t\in\mathbb R.$
 
Moreover, $\mathcal A_{X}(t)=\mathcal A_{Y}(t)$ for all $t\in\mathbb R$ if the following two conditions hold:
\begin{itemize} 
\item[$(i)$] $\mathcal A_{X}(t)$ is a compact subset of $X_t$ for all $t\in\mathbb R;$ 

\item [$(ii)$] for any $\widehat D_2\in\mathcal D_{Y}$ and $t\in\mathbb R$ there exist a family $\widehat D_1\in\mathcal D_{X}$ and a $t^*_{\widehat D_1}$ such that $\mathcal U$ is pullback $\widehat D_1$-asymptotically compact, and for any $s\le t^*_{\widehat D_1}$ there exists a $\tau_s<s$ such that $\mathcal U(s,\tau)D_2(\tau)\subset D_1(s)$ for all $\tau\le\tau_s.$
\end{itemize}
\end{theorem}

\begin{definition}
    We say that the families of metric spaces $\{X_t\}_{t\in \mathbb R}$ and $\{Y_t\}_{t\in \mathbb R}$ are \textbf{equivalent} if there exists a family of maps $\{\Theta(t)\}_{t\in\mathbb R}$ such that $\Theta(t):X_t \to Y_t$ is a homeomorphism for each $t\in \mathbb R$.
\end{definition}

\begin{definition}
    Let $\mathcal D_X$ and $\mathcal D_Y$ be two universes on the equivalent families of metric spaces $\{X_t\}_{t\in \mathbb R}$ and $\{Y_t\}_{t\in \mathbb R}$, respectively. We say that the universes $\mathcal D_X$ and $\mathcal D_Y$ are equivalent if  for any $\widehat D \in \mathcal D_X$ the family $\{\Theta(t)D(t)\subset Y_t:t\in \mathbb R\} \in \mathcal D_Y$ and for any $\widehat C \in \mathcal D_Y$ the family $\{\Theta^{-1}(t)C(t)\subset X_t : t\in \mathbb R\}\in \mathcal D_X$.
\end{definition}

 \begin{theorem}(cf. \cite[Theorem 3.16]{HRHJ-2024})
    Let $\{X_t\}_{t\in \mathbb R}$ and $\{Y_t\}_{t\in \mathbb R}$ be two equivalent families of metric spaces. Let $U_{X}$ be a continuous process on the family $\{X_t\}_{t\in \mathbb R}$. Then, $U_{Y}$ defined as
    \begin{equation}\label{process-Y}
        U_{Y}(t,\tau)y := \Theta(t) \circ U_{X}(t,\tau) \circ \Theta^{-1}(\tau)y \hspace{0,3cm}\text{ for all }(t,\tau)\in \mathbb R_d^2,\; y\in Y_{\tau},
    \end{equation}
    is a continuous process on the family $\{Y_t\}_{t\in \mathbb R}$.
\end{theorem}

\begin{theorem}\label{t5}
(cf. \cite[Theorem 3.20]{HRHJ-2024}) Let $\{X_t\}_{t\in \mathbb R}$ and $\{Y_t\}_{t\in \mathbb R}$ be two equivalent families of metric spaces. Consider $U_{X}$ a process and a universe $\mathcal D_{X}$ on the family $\{X_t\}_{t\in \mathbb R}$, and let us consider the process $U_{Y}$, given in \eqref{process-Y}, and the universe $\mathcal D_{Y}$. If $\mathcal A_{\mathcal D_{X}}$ is the minimal pullback $\mathcal D_{X}$-attractor for the continuous process $U_{X}$, then the family $\mathcal A_{\mathcal D_{Y}}=\{\mathcal A_{\mathcal D_{Y}}(t): t\in \mathbb R\}$ defined by
      $$  \mathcal A_{\mathcal D_{Y}} (t):=\Theta(t)\mathcal A_{\mathcal D_{X}}(t), \quad t\in \mathbb R, $$ is the minimal pullback $\mathcal D_{Y}$-attractor for the continuous process $U_{Y}$.
\end{theorem}

\subsection{Singular Gronwall inequality with exponential decay} The estimates below rely on a Volterra-type inequality with a kernel that is weakly singular at the origin and exponentially tempered at infinity. The singular part follows the classical iteration technique of Henry \cite[Lemma 7.1.1]{Henry1981}, later refined by Ye, Gao and Ding \cite{Ye2007} and by Webb \cite{Webb2019}; however, since their kernels lack a tempering factor, that approach alone only yields Mittag--Leffler type growth, not decay. We therefore combine it with resolvent and small gain arguments for Volterra equations with $L^1$ kernels in exponentially weighted spaces, which are classical (see, e.g., Gripenberg, Londen and Staffans \cite{GLS1990}). Lemma \ref{lem:main} below should be regarded as a variant of these classical results; we include a short self-contained proof, since it yields an explicit exponential decay rate and explicit constants, which are needed in our estimates.
 
Throughout this subsection, $\alpha$, $\kappa$ and $\gamma$ denote auxiliary parameters, unrelated to the boundary curve $\alpha(t)$, the chemotactic coefficient $\kappa$ and the length $\gamma(t)$. Let $\alpha\in(0,1)$ and $C>0$ be fixed constants, and define the kernel
\begin{equation}\label{eq:kernel}
k(\sigma):=C\,\sigma^{-\alpha}e^{-\sigma},\qquad \sigma>0 .
\end{equation}
This kernel combines a weak (integrable) singularity at $\sigma=0$ with
exponential decay as $\sigma\to\infty$; in particular
\begin{equation}\label{eq:kappa-def}
\kappa:=\int_0^\infty k(\sigma)\,d\sigma
=C\int_0^\infty\sigma^{-\alpha}e^{-\sigma}\,d\sigma
=C\,\Gamma(1-\alpha)\;<\;\infty .
\end{equation}
 
\begin{lemma}[Singular Gronwall inequality with exponential decay]\label{lem:main}
Let $\tau\in\mathbb{R}$, $\lambda>0$, $A,B\geq0$, and let $y:[\tau,\infty)\to[0,\infty)$ be continuous and satisfy
\begin{equation}\label{eq:hyp-ineq}
y(t)\;\leq\; A\,e^{-\lambda(t-\tau)}+B
+\int_\tau^t k(t-s)\,y(s)\,ds,
\qquad t\geq\tau .
\end{equation}
Assume the \emph{small-gain condition}
\begin{equation}\label{eq:small-gain}
\tag{H}
\kappa=C\,\Gamma(1-\alpha)<1 .
\end{equation}
Then there exist a rate $\delta^{\ast }=\delta^{\ast }(C,\alpha,\lambda)\in
\big(0,\min\{1,\lambda\}\big)$ and a constant
$C_0=C_0(C,\alpha,\lambda)\geq1$, both given explicitly in the proof, such
that
\begin{equation}\label{eq:conclusion}
y(t)\;\leq\;C_0\,A\,e^{-\delta^{\ast }(t-\tau)}+C_0\,B,
\qquad \forall\,t\geq\tau .
\end{equation}
\end{lemma}
 \begin{proof}Throughout, all convolutions are taken on $[0,\infty)$:
$(f*g)(\sigma):=\int_0^\sigma f(\sigma-r)g(r)\,dr$, and
$\int_\tau^t k(t-s)f(s)\,ds=(k\ast f)(t-\tau)$ after the shift $s\mapsto
s-\tau$. We identify $y$ with a function of $t-\tau\geq0$ where
convenient.
 
Since $e^{-\lambda(t-\tau)}\leq1$, and, by \eqref{eq:kappa-def},
$\int_\tau^t k(t-s)\,ds\leq\kappa$ for every $t\geq\tau$, inequality
\eqref{eq:hyp-ineq} gives, for every $T\geq\tau$ and every
$t\in[\tau,T]$,
\[
y(t)\;\leq\;A+B+\kappa\sup_{[\tau,T]}y .
\]
Taking the supremum over $t\in[\tau,T]$ (finite because $y$ is
continuous on the compact interval $[\tau,T]$) and using
\eqref{eq:small-gain},
\[
\sup_{[\tau,T]}y\;\leq\;\frac{A+B}{1-\kappa},
\]
a bound independent of $T$. Letting $T\to\infty$,
\begin{equation}\label{eq:Z0}
Z_0:=\sup_{t\geq\tau}y(t)\;\leq\;\frac{A+B}{1-\kappa}\;<\;\infty .
\end{equation}
In particular $y$ is bounded on $[\tau,\infty)$.
 
Now, we are going to remove the constant from the force term. Set
\[
z(t):=y(t)-\frac{B}{1-\kappa},\qquad t\geq\tau .\]
Substituting $y(s)=z(s)+\dfrac{B}{1-\kappa}$ into \eqref{eq:hyp-ineq} and
using $\int_\tau^t k(t-s)\,ds\leq\kappa$,
\begin{align*}
z(t)+\frac{B}{1-\kappa}
&\leq A\,e^{-\lambda(t-\tau)}+B
+\int_\tau^t k(t-s)z(s)\,ds
+\frac{B}{1-\kappa}\int_\tau^t k(t-s)\,ds\\
&\leq A\,e^{-\lambda(t-\tau)}+B+\frac{B\kappa}{1-\kappa}
+\int_\tau^t k(t-s)z(s)\,ds .
\end{align*}
Since $B+\dfrac{B\kappa}{1-\kappa}-\dfrac{B}{1-\kappa}
=B\cdot\dfrac{(1-\kappa)+\kappa-1}{1-\kappa}=0$, the constant terms
cancel exactly, leaving
\begin{equation}\label{eq:z-ineq}
z(t)\;\leq\;\tilde a(t)+\int_\tau^t k(t-s)\,z(s)\,ds,
\qquad \tilde a(t):=A\,e^{-\lambda(t-\tau)} .
\end{equation}
Note that by \eqref{eq:Z0}, $z$ is bounded:
$|z(t)|\leq Z_0+\dfrac{B}{1-\kappa}=:Z_0'<\infty$ for all $t\geq\tau$.
 
On the other hand, since $k\geq0$, inequality \eqref{eq:z-ineq} may be iterated: substituting
the bound \eqref{eq:z-ineq} for $z(s)$ back into the integral (valid
because $k\geq0$ preserves the inequality direction) and applying Fubini theorem, one obtains, for every $n\in\mathbb N$,
\begin{equation}\label{eq:iterate}
z(t)\;\leq\;\sum_{j=0}^{n-1}\big(k^{\ast j}\ast \tilde a\big)(t)
+\big(k^{\ast n}\ast z\big)(t),
\qquad k^{\ast 0}\ast \tilde a:=\tilde a,
\end{equation}
where $k^{\ast j}$ is the $j$-fold Volterra convolution power of $k$, that is,
\[
k^{\ast 1}:=k,
\qquad
k^{\ast n}:=k*k^{\ast (n-1)},\qquad n\geq2,
\]
so that, explicitly,
\[
k^{\ast n}(\sigma)
=\int_0^\sigma\int_0^{\sigma_1}\cdots\int_0^{\sigma_{n-2}}
k(\sigma-\sigma_1)\,k(\sigma_1-\sigma_2)\cdots k(\sigma_{n-2}-\sigma_{n-1})\,
k(\sigma_{n-1})\;d\sigma_{n-1}\cdots d\sigma_1 .
\]
By convention, $k^{*0}*f:=f$ for any function $f$. Thus, since
$\int_0^\infty k^{*j}(\sigma)\,d\sigma=\kappa^{j}$ (a standard property of
convolution integrals), and $|z|\leq Z_0'$ as shown above,
\begin{equation*}
\big|(k^{\ast n}*z)(t)\big|
\leq Z_0'\int_0^{t-\tau}k^{*n}(\sigma)\,d\sigma
\leq Z_0'\,\kappa^{n}
\xrightarrow[n\to\infty]{}0,
\end{equation*}
because $\kappa<1$ by \eqref{eq:small-gain}. Letting $n\to\infty$ in
\eqref{eq:iterate},
\begin{equation}\label{eq:resolvent-bound}
z(t)\;\leq\;\tilde a(t)+(R \ast \tilde a)(t),
\qquad
R:=\sum_{j=1}^\infty k^{*j},
\qquad
\int_0^\infty R(\sigma)\,d\sigma=\sum_{j=1}^\infty\kappa^{j}
=\frac{\kappa}{1-\kappa}<\infty .
\end{equation}
 
Now, for $\theta\in[0,1)$ define 
\begin{equation*}
\kappa(\theta):=\int_0^\infty e^{\theta\sigma}k(\sigma)\,d\sigma
=C\int_0^\infty\sigma^{-\alpha}e^{-(1-\theta)\sigma}\,d\sigma
=C\,\Gamma(1-\alpha)\,(1-\theta)^{-(1-\alpha)} .
\end{equation*}
The map $\theta\mapsto\kappa(\theta)$ is continuous and strictly
increasing on $[0,1)$, with $\kappa(0)=\kappa<1$ by \eqref{eq:small-gain}
and $\kappa(\theta)\to\infty$ as $\theta\to1^{-}$. By the intermediate
value theorem there is a unique $\theta^{\ast}\in(0,1)$ with $\kappa(\theta^{\ast})=1$, and $\kappa(\theta)<1$ for every $\theta\in[0,\theta^{\ast})$. Define
\begin{equation}\label{eq:delta-star}
\delta^{\ast}:=\tfrac12\min\{\theta^{\ast},\,1,\,\lambda\}\;>\;0 .
\end{equation}
Then $\kappa(\delta^{\ast})<1$, and $\delta^{\ast}\leq\lambda$, so that
\begin{equation}\label{eq:atilde-tilt}
\tilde a(t)=A\,e^{-\lambda(t-\tau)}\;\leq\;A\,e^{-\delta^{\ast}(t-\tau)},
\qquad t\geq\tau .
\end{equation}
Using the elementary identity $e^{\gamma\sigma}(f*g)(\sigma)
=(f_\gamma*g_\gamma)(\sigma)$, where $f_\gamma(\sigma):=e^{\gamma\sigma}
f(\sigma)$, applied repeatedly to $k^{*j}$,
\[
\int_0^\infty e^{\delta^{\ast}\sigma}R(\sigma)\,d\sigma
=\sum_{j=1}^\infty\int_0^\infty e^{\delta^{\ast}\sigma}k^{\ast j}(\sigma)\,d\sigma
=\sum_{j=1}^\infty\kappa(\delta^{\ast})^{j}
=\frac{\kappa(\delta^{\ast})}{1-\kappa(\delta^{\ast})}\;<\;\infty .
\]
Hence, using \eqref{eq:atilde-tilt} and the change of variables
$\sigma=t-s$,
\begin{align}\label{Step4}
\begin{array}{l}
\displaystyle (R*\tilde a)(t)
\leq A\int_0^{t-\tau}R(\sigma)\,e^{-\delta^{\ast}(t-\tau-\sigma)}\,d\sigma
=A\,e^{-\delta^{\ast}(t-\tau)}\int_0^{t-\tau}R(\sigma)e^{\delta^{*}\sigma}\,d\sigma\\ \displaystyle \hspace{1.8cm}
\leq A\,\frac{\kappa(\delta^{\ast})}{1-\kappa(\delta^{*})}\,
e^{-\delta^{\ast}(t-\tau)} .
\end{array}
\end{align}
 
Finally, combining \eqref{Step4} with \eqref{eq:resolvent-bound} and
\eqref{eq:atilde-tilt}, we have
\[
z(t)\;\leq\;A\Big[1+\frac{\kappa(\delta^{\ast})}{1-\kappa(\delta^{*})}\Big]
e^{-\delta^{\ast}(t-\tau)}
=\frac{A}{1-\kappa(\delta^{\ast})}\,e^{-\delta^{\ast}(t-\tau)} .
\]
Recalling $y(t)=z(t)+\dfrac{B}{1-\kappa}$ and setting
\begin{equation}\label{eq:C0}
C_0:=\frac{1}{1-\kappa(\delta^{\ast})}\;\geq\;\frac{1}{1-\kappa}\;\geq\;1
\end{equation}
(the middle inequality holds since $\kappa$ is nondecreasing and
$\delta^{*}\geq0$), we obtain
\[
y(t)\;\leq\;C_0\,A\,e^{-\delta^{\ast}(t-\tau)}+\frac{B}{1-\kappa}
\;\leq\;C_0\,A\,e^{-\delta^{\ast}(t-\tau)}+C_0\,B ,
\]
which is \eqref{eq:conclusion}. 
\end{proof}
 
\begin{corollary} If $B=0$ in Lemma \ref{lem:main}, then, under the same hypotheses,
\[
y(t)\;\leq\;C_0\,A\,e^{-\delta^{\ast}(t-\tau)},\qquad t\geq\tau,
\]
with $\delta^{\ast}$ and $C_0$ as in \eqref{eq:delta-star}--\eqref{eq:C0}.
\end{corollary}
 
\begin{remark}[The classical Henry--Gronwall lemma does not yield decay] The classical result of Henry \cite[Lemma~7.1.1]{Henry1981} treats inequalities of the form $u(t)\leq a(t)+b\int_0^t(t-s)^{\beta-1}u(s)\,ds$ and yields $u(t)\leq a(t)E_\beta(b\Gamma(\beta)t^\beta)$, where $E_\beta$ is the Mittag--Leffler function. The kernel $(t-s)^{\beta-1}$ there is not integrable on   $(0,\infty)$, so no analogue of $\kappa$ in \eqref{eq:kappa-def} is finite, and the Neumann series does not converge geometrically: $E_\beta$ grows (sub-exponentially but unboundedly) rather than decaying. The mechanism that produces decay in Lemma \ref{lem:main} is precisely the exponential tempering factor $e^{-\sigma}$ in \eqref{eq:kernel}, which makes $k\in L^1(0,\infty)$ and allows the small-gain condition \eqref{eq:small-gain} together with the tilting argument in \eqref{Step4}.
\end{remark}

\subsection{Semilinear abstract differential equation} Let $\mathcal{V}$ and $\mathcal{H}$ be two separable real Hilbert spaces with dense and compact embedding $\mathcal{V}\hookrightarrow \mathcal{H}$. Identifying $\mathcal{H}$ and its dual $\mathcal{H}^{'}$ and denoting the dual space of $\mathcal{V}$ by $\mathcal{V}^{'}$, we have $\mathcal{V}\hookrightarrow \mathcal{H}\hookrightarrow \mathcal{V}^{'}$. We denote the scalar product of $\mathcal{H}$ by $(\cdot,\cdot)$ and the norm by $|\cdot|$. The duality product between $\mathcal{V}^{'}$ and $\mathcal{V}$ which coincides with the scalar product of $\mathcal{H}$ on $\mathcal{H}\times\mathcal{H}$ is denoted by $\langle \cdot, \cdot\rangle$, and the norms of $\mathcal{V}$ and $\mathcal{V}^{'}$ are denoted by $\|\cdot\|$ and $\|\cdot \|_{\ast}$, respectively.

In this section, we shall first prove the existence and uniqueness of a weak solution for the Cauchy problem of a semilinear abstract differential equation,
\begin{equation}\label{E}\tag{P}
\begin{cases}
    \displaystyle \frac{dX}{dt}+AX=F(t,X), \hspace{0,5cm}t>\tau,\\
    X(\tau) = X_\tau
    \end{cases}
\end{equation}
in the space $\mathcal{V}^{'}$. Here, $A$ is the positive definite self-adjoint operator of $\mathcal{H}$ defined by a symmetric bilinear form $a(X,Y)$ on $\mathcal{V}$, i.e., $\langle AX,Y\rangle= a(X,Y)$, which satisfies 
\begin{equation}\label{ai}
    |a(X,Y)|\leq M \|X\|\|Y\|, \hspace{0,5cm} X,Y \in \mathcal{V},
\end{equation}
\begin{equation}\label{aii}
    a(X,X) \geq \delta \|X\|^2, \hspace{0,5cm} X\in \mathcal{V},
\end{equation}
with some $\delta$ and $M>0$. $A$ is also a bounded operator from $\mathcal{V}$ to $\mathcal{V}'$. $F(\cdot,\cdot)$ is a given continuous function from $\mathbb{R} \times \mathcal{V}$ to $\mathcal{V}'$ satisfying:

\vspace{0,5cm}

\noindent {\bf (F.i)} There is $U \in L^{2}_{loc}(\mathbb{R};\mathcal{V}')$ such that for each $\eta>0$, there exist $c_\eta >0$ and a continuous function $\Psi_{\eta}: \mathbb{R}\times [0,\infty) \to [0,\infty)$ such that for each $t\in \mathbb{R}$, $\Psi_{\eta}(t,\cdot):[0,\infty) \to [0,\infty)$ is increasing and 
\begin{equation*}
    \|F(t,X)\|_{\ast} \leq \eta \|X\| + \Psi_{\eta}(t,|X|) + c_\eta\|U (t)\|_\ast, \hspace{0,5cm} X\in \mathcal{V}.
\end{equation*}

\noindent {\bf (F.ii)} For each $\eta>0$, there exist  continuous function $\Phi_{\eta}: \mathbb{R}\times [0,\infty) \to [0,\infty)$ such that for each $t\in \mathbb{R}$, $\Phi_{\eta}(t,\cdot):[0,\infty) \to [0,\infty)$ is increasing and 
\begin{align*}
    \|F(t,X)- F(t,Y)\|_{\ast} \leq \eta &\|X-Y\|  +(1+\|X\|+\|Y\|)\times\Phi_{\eta}(t,|X| + |Y|)|X-Y|, \hspace{0.5cm} X,Y\in \mathcal{V}.
\end{align*}

\begin{definition}[Weak solution] Let $X_\tau\in\mathcal{H}$ and $\tau<T$. A function $X:[\tau,T]\rightarrow\mathcal{H}$ is called a \emph{weak solution} of problem \eqref{E} on $[\tau,T]$ if
\begin{equation*}
X\in L^{2}(\tau,T;\mathcal{V})\cap C([\tau,T];\mathcal{H}),
\quad \frac{dX}{dt}\in L^{2}(\tau,T;\mathcal{V}'),
\end{equation*}
satisfies the initial condition
\[
X(\tau)=X_\tau \quad \text{in } \mathcal{H},
\]
and
\begin{equation*}
    \left\langle \frac{dX}{dt}(t),Y\right\rangle +a(X(t),Y) = \big\langle F(t,X(t)),Y\big\rangle,
\end{equation*}
for every $Y\in\mathcal{V}$ and for almost every $t\in(\tau,T)$.
\end{definition}

The proof of the following theorem follows along the same lines as that of the corresponding result in \cite[Theorem 2.1]{ryu2001}, taking into account the modifications arising from the non-autonomous term $F(t,X)$ in system \eqref{E}.

\vspace{0,5cm}

\begin{theorem}\label{T-Exis} (cf. \cite[Theorem 2.1]{ryu2001})
    Let \eqref{ai}, \eqref{aii}, {\bf(F.i)} and {\bf(F.ii)} be satisfied. Then, for any $X_0 \in \mathcal{H}$, there exists a unique weak solution 
    \begin{equation}
        X\in H^1 (\tau,T_0 ; \mathcal{V}')\cap  C ([\tau,T_0] ; \mathcal{H})\cap  L^2 (\tau,T_0 ; \mathcal{V})
    \end{equation}
    to \eqref{E}; the number $T_0 = T_0 (|X_0|, \tau, \|U\|_{L^2(\tau,T;\mathcal{V}')})>\tau$, where $U \in L^{2}_{loc}(\mathbb{R};\mathcal{V}')$ was given in {\bf(F.i)}.
\end{theorem}

\subsubsection{ \bf Some bilinear forms} For $m_0,m_1>0$, define
$\mathbf a_{m_0,m_1}:H^1(\mathcal I)\times H^1(\mathcal I)\to\mathbb R$ by
\[
\mathbf a_{m_0,m_1}(u,v):=m_1\int_{\mathcal I}u_xv_x\,dx+m_0\int_{\mathcal I}uv\,dx,
\qquad u,v\in H^1(\mathcal I) .
\]
This form is symmetric, continuous and coercive on $H^1(\mathcal I)$:
there are constants $C,c>0$ with
$|\mathbf a_{m_0,m_1}(u,v)|\leq C\|u\|_{H^1}\|v\|_{H^1}$ and
$\mathbf a_{m_0,m_1}(u,u)\geq c\|u\|_{H^1}^2$. By the Lax--Milgram
theorem, for every $f\in H^1(\mathcal I)'$ there is a unique $u\in
H^1(\mathcal I)$ with $\mathbf a_{m_0,m_1}(u,v)=\langle
f,v\rangle_{H^1(\mathcal I)',H^1(\mathcal I)}$ for all $v\in H^1(\mathcal
I)$: this is the weak formulation of the Neumann problem
$-m_1u_{xx}+m_0u=f$ on $\mathcal I$, $u_x(0)=u_x(1)=0$. The associated
operator
\[
\widetilde A_{m_0,m_1}:H^1(\mathcal I)\to H^1(\mathcal I)',
\qquad
\langle\widetilde A_{m_0,m_1}u,v\rangle_{H^1(\mathcal I)',H^1(\mathcal I)}:=\mathbf a_{m_0,m_1}(u,v),
\]
thus satisfies $\widetilde A_{m_0,m_1}u=-m_1u_{xx}+m_0u$ in the
distributional sense.
 
\paragraph{From the operator back to a realization on $L^2$.}
Identifying $L^2(\mathcal I)$ with its dual, the $L^2$-realization of
$\widetilde A_{m_0,m_1}$ is
\[
A_{m_0,m_1}:D(A_{m_0,m_1})\subset L^2(\mathcal I)\to L^2(\mathcal I),
\qquad
D(A_{m_0,m_1}):=\big\{u\in H^1(\mathcal I):\widetilde A_{m_0,m_1}u\in L^2(\mathcal I)\big\},
\]
with $A_{m_0,m_1}u=\widetilde A_{m_0,m_1}u$. By elliptic regularity,
$D(A_{m_0,m_1})=H^2_N(\mathcal I)$, and $A_{m_0,m_1}u=-m_1u_{xx}+m_0u$
classically. The operator $A_{m_0,m_1}$ is self-adjoint, positive, with
compact resolvent, hence with discrete spectrum
$0<m_0=\lambda_0<\lambda_1\leq\lambda_2\leq\cdots\to\infty$ (the
lowest eigenvalue $\lambda_0=m_0$ corresponds to constant
eigenfunctions).
 
\paragraph{\bf Fractional powers, and back to bilinear forms.} For $\theta\geq0$, the domain of the fractional power
$A_{m_0,m_1}^\theta$ is
\[
D\big(A_{m_0,m_1}^\theta\big)=
\begin{cases}
H^{2\theta}(\mathcal I), & 0\leq\theta<3/4,\\
H^{2\theta}_N(\mathcal I), & 3/4<\theta<7/4,\\
H^{2\theta}_{N^2}(\mathcal I), & 7/4<\theta\leq2 .
\end{cases}
\]
This closes the loop: associated with $A_{m_0,m_1}^\theta$ is the
bilinear form
\[
\mathbf a_{m_0,m_1}^{(\theta)}(u,v):=\big(A_{m_0,m_1}^{\theta/2}u,A_{m_0,m_1}^{\theta/2}v\big)_{L^2(\mathcal I)},
\qquad u,v\in D\big(A_{m_0,m_1}^{\theta/2}\big),
\]
which, by self-adjointness, satisfies
$\mathbf a_{m_0,m_1}^{(\theta)}(u,v)=(A_{m_0,m_1}^\theta u,v)_{L^2(\mathcal
I)}$ for $u\in D(A_{m_0,m_1}^\theta)$. Taking $\theta=1$ recovers the
original form: $\mathbf a_{m_0,m_1}^{(1)}=\mathbf a_{m_0,m_1}$.
 
\bigskip
The next two examples instantiate this construction on the product
spaces used throughout the paper.
 
\paragraph{\bf Example 01 ($\mathcal H=L^2\times H^1$, $\mathcal V=H^1\times H^2_N$).} Given $n_0,n_1,m_0,m_1>0$, let $\widetilde A_{n_0,n_1}:H^1(\mathcal
I)\to H^1(\mathcal I)'$ be as above and let $A_{m_0,m_1}:H^2_N(\mathcal
I)\to L^2(\mathcal I)$ be the $L^2$-realization of $\widetilde
A_{m_0,m_1}$. Since $A_{m_0,m_1}$ identifies $\big(H^2_N(\mathcal
I)\big)'\simeq L^2(\mathcal I)$, we get $\mathcal V'\simeq H^1(\mathcal
I)'\times L^2(\mathcal I)$. The operator matrix
\[
\mathbf A=\begin{pmatrix}\widetilde A_{n_0,n_1}&0\\0&A_{m_0,m_1}\end{pmatrix}:\mathcal V\to\mathcal V'
\]
is associated with the bilinear form
\[
\mathbf a_{\mathbf A}\big((u,\rho),(v,\phi)\big):=\mathbf a_{n_0,n_1}(u,v)+\big(A_{m_0,m_1}\rho,A_{m_0,m_1}\phi\big)_{L^2(\mathcal I)},
\]
i.e.\ $\mathbf a_{\mathbf A}(U,V)=\langle\mathbf AU,V\rangle_{\mathcal
V',\mathcal V}$ for $U,V\in\mathcal V$.
 
\paragraph{\bf Example 02 ($\mathcal H=H^1\times H^2_N$, $\mathcal V=H^2_N\times H^3_N$).}
Let $A_{n_0,n_1}:H^2_N(\mathcal I)\to L^2(\mathcal I)$ be the
$L^2$-realization of $\widetilde A_{n_0,n_1}$, and let $A_{m_0,m_1}:H^3_N(\mathcal I)\to H^1(\mathcal I)$ be the $H^1$-realization of $\widetilde A_{m_0,m_1}$, where $H^1(\mathcal I)$ carries the equivalent inner product
$(\rho,\phi)_{H^1(\mathcal I)}:=\big(A_{m_0,m_1}^{1/2}\rho,A_{m_0,m_1}^{1/2}\phi\big)_{L^2(\mathcal
I)}$ (the case $\theta=1$ of the general form above). As before,
$\big(H^2_N(\mathcal I)\big)'\simeq L^2(\mathcal I)$ and $\big(H^3_N(\mathcal
I)\big)'\simeq H^1(\mathcal I)$, so $\mathcal V'\simeq L^2(\mathcal
I)\times H^1(\mathcal I)$. The operator matrix
\[
\mathbf A=\begin{pmatrix}A_{n_0,n_1}&0\\0&A_{m_0,m_1}\end{pmatrix}:\mathcal V\to\mathcal V'
\]
is associated with the bilinear form
\begin{align*}
\mathbf a_{\mathbf A}\big((u,\rho),(v,\phi)\big)
&:=\big(A_{n_0,n_1}u,A_{n_0,n_1}v\big)_{L^2(\mathcal I)}+\big(A_{m_0,m_1}\rho,A_{m_0,m_1}\phi\big)_{H^1(\mathcal I)}
\\&=\big(A_{n_0,n_1}u,A_{n_0,n_1}v\big)_{L^2(\mathcal I)}+\big(A_{m_0,m_1}^{3/2}\rho,A_{m_0,m_1}^{3/2}\phi\big)_{L^2(\mathcal I)},
\end{align*}
using $(A_{m_0,m_1}\rho,A_{m_0,m_1}\phi)_{H^1}=(A_{m_0,m_1}^{1/2}A_{m_0,m_1}\rho,A_{m_0,m_1}^{1/2}A_{m_0,m_1}\phi)_{L^2}=(A_{m_0,m_1}^{3/2}\rho,A_{m_0,m_1}^{3/2}\phi)_{L^2}$;
as above, $\mathbf a_{\mathbf A}(U,V)=\langle\mathbf AU,V\rangle_{\mathcal
V',\mathcal V}$ for $U,V\in\mathcal V$.

\section{Existence and uniqueness of local weak solutions}\label{sec3}

\subsection{Existence of local weak solution} Now, we will prepare to prove the existence of a unique local weak solution to \eqref{equ4iii}. As a consequence, we will also obtain the existence of a local weak solution to \eqref{equ4} and \eqref{equ5}.

\begin{theorem}\label{teoEx-02}
\textup{(Existence of a local weak solution to \eqref{equ4iii})}
Under the assumptions given in (H1)-(H2), given $(w_0,\phi_0)\in
L^2(\mathcal I)\times H^1(\mathcal I)$ and $\tau\in\mathbb R$, there
exists $T_{w_0,\phi_0}=T_0(\tau,w_0,\phi_0)>\tau$ such that the semilinear
parabolic evolution problem \eqref{equ4iii} has a unique local weak
solution
\[
(w(\cdot),\phi(\cdot))=\big(w(\cdot;\tau,w_0,\phi_0),\phi(\cdot;\tau,w_0,\phi_0)\big):
[\tau,T_{w_0,\phi_0})\to L^2(\mathcal I)\times H^1(\mathcal I)
\]
with
\begin{align*}
w&\in H^1\big(\tau,T_{w_0,\phi_0};H^1(\mathcal I)'\big)\cap\mathcal C([\tau,T_{w_0,\phi_0}];L^2(\mathcal I))\cap L^2(\tau,T_{w_0,\phi_0};H^1(\mathcal I)),\\
\phi&\in H^1(\tau,T_{w_0,\phi_0};L^2(\mathcal I))\cap\mathcal C([\tau,T_{w_0,\phi_0}];H^1(\mathcal I))\cap L^2(\tau,T_{w_0,\phi_0};H^2_N(\mathcal I)),
\end{align*}
satisfying the initial condition $(w(\tau),\phi(\tau))=(w_0,\phi_0)$.
\end{theorem}

\begin{proof}
Let $\mathcal I=(0,1)\subset\mathbb R$ be a bounded interval, and consider
the Hilbert spaces
\[
\mathcal H=L^2(\mathcal I)\times H^1(\mathcal I),
\qquad
\mathcal V=H^1(\mathcal I)\times H^2_N(\mathcal I),
\qquad
\mathcal V'=\big(H^1(\mathcal I)\big)'\times\big(H^2_N(\mathcal I)\big)' .
\]
The embedding $\mathcal V\hookrightarrow\mathcal H$ is dense and
compact: compactness follows from the Rellich--Kondrachov theorem on the
bounded interval $\mathcal I$, and density from the fact that
$H^2_N(\mathcal I)$ contains the Neumann eigenfunctions
$\cos(k\pi y)$, $k\geq0$, whose linear span is dense in
$H^1(\mathcal I)$. The space $\mathcal V$ is endowed with its standard
product norm $\|\cdot\|_{\mathcal V}$, and we write
$\|z\|_{\mathcal H}:=\big(\|w\|_{L^2(\mathcal I)}^2+\|\phi\|_{H^1(\mathcal I)}^2\big)^{1/2}$
for the standard norm of $z=(w,\phi)\in\mathcal H$.
Let $\widetilde A_{a,1}:H^1(\mathcal I)\to(H^1(\mathcal I))'$ be the
operator associated with $-a\Delta+I$, and let
$A_{b,1}:H^2_N(\mathcal I)\to L^2(\mathcal I)$ be the $L^2$-realization
of the operator associated with $-b\Delta+I$. Using $A_{b,1}$, one
identifies $(H^2_N(\mathcal I))'\simeq L^2(\mathcal I)$, so that
$\mathcal V'\simeq(H^1(\mathcal I))'\times L^2(\mathcal I)$, with duality
pairing
\[
\Big\langle\begin{pmatrix}f\\ g\end{pmatrix},\begin{pmatrix}v\\ \xi\end{pmatrix}\Big\rangle_{\mathcal V',\mathcal V}
=\langle f,v\rangle_{H^1(\mathcal I)',H^1(\mathcal I)}+(g,A_{b,1}\xi)_{L^2(\mathcal I)} .
\]
For $(g,\xi)\in H^1(\mathcal I)\times H^2_N(\mathcal I)$ we have
$(g,A_{b,1}\xi)_{L^2(\mathcal I)}=b\,(\partial_yg,\partial_y\xi)_{L^2(\mathcal I)}+(g,\xi)_{L^2(\mathcal I)}$,
since $\partial_y\xi=0$ at $y\in\{0,1\}$. Hence, in order that this
pairing coincide with the scalar product of $\mathcal H$ on
$\mathcal H\times\mathcal H$, as required in Theorem~\ref{T-Exis}, we
endow $\mathcal H$ with the scalar product
\[
\big((w,\phi),(v,\xi)\big)_{\mathcal H}
:=(w,v)_{L^2(\mathcal I)}+b\,(\partial_y\phi,\partial_y\xi)_{L^2(\mathcal I)}+(\phi,\xi)_{L^2(\mathcal I)},
\]
whose norm $|\cdot|$ satisfies
$\min\{1,\sqrt b\}\,\|z\|_{\mathcal H}\leq|z|\leq\max\{1,\sqrt b\}\,\|z\|_{\mathcal H}$.
Define
\[
\mathbf A:=\begin{pmatrix}\widetilde A_{a,1}&0\\0&A_{b,1}\end{pmatrix}:\mathcal V\to\mathcal V',
\]
which is bounded and associated with the symmetric bilinear form
\[
\mathbf a_{\mathbf A}\big((w,\rho),(v,\xi)\big)
:=\langle\widetilde A_{a,1}w,v\rangle_{H^1(\mathcal I)',H^1(\mathcal I)}
+\big(A_{b,1}\rho,A_{b,1}\xi\big)_{L^2(\mathcal I)},
\]
satisfying \eqref{ai} and \eqref{aii}.

By (H1), the function $s(\cdot)$ in \eqref{rescaling} (whose defining
formula makes sense for every $t\in\mathbb R$) is of class $C^1$, with
$s'(t)=\gamma^{-2}(t)\in(\gamma_1^{-2},\gamma_0^{-2})$; hence
$s:\mathbb R\to\mathbb R$ is an increasing $C^1$-diffeomorphism, and its
inverse $t(\cdot)=s^{-1}(\cdot)$ satisfies
\begin{equation*}
t'(s)=\gamma^2(t(s))\in(\gamma_0^2,\gamma_1^2),\qquad s\in\mathbb R .
\end{equation*}
The abstract formulation of \eqref{equ4iii} is then
\begin{equation*}
\begin{cases}
z'+\mathbf A z=F(s,z), & s>\tau,\\
z(\tau)=(w_0,\phi_0)\in\mathcal H,
\end{cases}
\end{equation*}
where $z=(w,\phi)$ and
\begin{equation*}
F(s,z)=\begin{pmatrix}
w-\kappa\,\partial_y\big(w\,\partial_y\phi\big)-\gamma^2(t(s))\,b(y,t(s))\,\partial_yw\\[2pt]
\phi+c\,\gamma^2(t(s))\,w-d\,\gamma^2(t(s))\,\phi-\gamma^2(t(s))\,b(y,t(s))\,\partial_y\phi
\end{pmatrix},
\end{equation*}
with
\begin{equation*}
b(y,t):=-\frac{\alpha'(t)+\gamma'(t)\,y}{\gamma(t)},\qquad (y,t)\in\overline{\mathcal I}\times\mathbb R,
\end{equation*}
where $'$ denotes differentiation with respect to the original time
variable $t$. (Here and below, the diffusion coefficient $b>0$ of the
$\phi$-equation is a constant, while the transport coefficient is always
written with its arguments, $b(y,t)$.)

To apply Theorem~\ref{T-Exis} it remains to verify that
$F:\mathbb R\times\mathcal V\to\mathcal V'$ satisfies conditions
\textbf{(F.i)} and \textbf{(F.ii)}. Throughout, we write
\begin{equation*}
b_\gamma(t):=\frac{|\alpha'(t)|+|\gamma'(t)|}{\gamma(t)},
\qquad
m(s):=\gamma^2(t(s))\,b_\gamma(t(s))=\gamma(t(s))\big(|\alpha'(t(s))|+|\gamma'(t(s))|\big),
\end{equation*}
and we set $B(y,s):=\gamma^2(t(s))\,b(y,t(s))$, the coefficient that
actually appears in $F$. Since $\partial_yb(y,t)=-\gamma'(t)/\gamma(t)$
and $y\in[0,1]$, we have, for all $y\in\overline{\mathcal I}$ and
$s\in\mathbb R$,
\begin{equation}\label{eq:B-bounds}
|B(y,s)|\leq m(s),\qquad
|\partial_yB(y,s)|=\gamma(t(s))\,|\gamma'(t(s))|\leq m(s).
\end{equation}
By (H1), $\alpha',\gamma'$ and $t(\cdot)$ are continuous, so
$s\mapsto m(s)$ and $s\mapsto\gamma^2(t(s))$ are continuous on
$\mathbb R$. We denote by $C_1,C_2$ the constants in the one-dimensional
interpolation inequality
$\|f\|_{L^\infty(\mathcal I)}\leq C_1\|f\|_{L^2(\mathcal I)}^{1/2}\|f\|_{H^1(\mathcal I)}^{1/2}$
and the embedding $H^1(\mathcal I)\hookrightarrow L^\infty(\mathcal I)$
(with $\|v\|_{L^\infty(\mathcal I)}\leq C_2\|v\|_{H^1(\mathcal I)}$),
respectively, and we set $c_b:=\sqrt2\max\{1,b\}$, so that
$\|A_{b,1}\xi\|_{L^2(\mathcal I)}=\|{-b}\,\partial_y^2\xi+\xi\|_{L^2(\mathcal I)}\leq c_b\|\xi\|_{H^2(\mathcal I)}$
for every $\xi\in H^2_N(\mathcal I)$.

\smallskip
\noindent\emph{An auxiliary estimate.} Let $u\in H^1(\mathcal I)$. For
every $\zeta\in H^1(\mathcal I)$ with $\|\zeta\|_{H^1(\mathcal I)}\leq1$,
integrating by parts,
\[
\int_0^1B(y,s)\,\partial_yu\,\zeta\,dy
=\big[B(y,s)\,u\,\zeta\big]_{y=0}^{y=1}
-\int_0^1u\,\big[\partial_yB(y,s)\,\zeta+B(y,s)\,\partial_y\zeta\big]\,dy .
\]
Using \eqref{eq:B-bounds},
$\|\zeta\|_{L^\infty(\mathcal I)}\leq C_2$ and
$\|\zeta\|_{L^2(\mathcal I)},\|\partial_y\zeta\|_{L^2(\mathcal I)}\leq1$,
and then the interpolation inequality,
\begin{equation}\label{eq:aux-B}
\sup_{\|\zeta\|_{H^1(\mathcal I)}\leq1}\Big|\int_0^1B(y,s)\,\partial_yu\,\zeta\,dy\Big|
\leq2C_1C_2\,m(s)\,\|u\|_{L^2(\mathcal I)}^{1/2}\|u\|_{H^1(\mathcal I)}^{1/2}
+2\,m(s)\,\|u\|_{L^2(\mathcal I)} .
\end{equation}
We shall also use Young's inequality in the form
\begin{equation}\label{eq:Young}
A\,\|u\|_{L^2(\mathcal I)}^{1/2}\|u\|_{H^1(\mathcal I)}^{1/2}
\leq\varepsilon\,\|u\|_{H^1(\mathcal I)}+\frac{A^2}{4\varepsilon}\,\|u\|_{L^2(\mathcal I)},
\qquad A\geq0,\ \varepsilon>0 .
\end{equation}

\smallskip
\noindent{\bf Verification of (F.i).} Let $z=(w,\phi)\in\mathcal V$. By
definition of the duality pairing,
\begin{equation*}
\begin{array}{l}
\|F(s,w,\phi)\|_{\mathcal V'}
\leq\underbrace{\sup_{\|\zeta\|_{H^1(\mathcal I)}\leq1}\Big|\int_0^1\Big[w-\kappa\,\partial_y(w\,\partial_y\phi)-B(y,s)\,\partial_yw\Big]\zeta\,dy\Big|}_{=:I_1}
\\ \hspace{2.2cm}+\underbrace{\sup_{\|\xi\|_{H^2_N(\mathcal I)}\leq1}\Big|\int_0^1\Big[\phi+c\,\gamma^2(t(s))\,w-d\,\gamma^2(t(s))\,\phi-B(y,s)\,\partial_y\phi\Big]\,A_{b,1}\xi\,dy\Big|}_{=:I_2} .
\end{array}
\end{equation*}

\smallskip
\noindent{\bf Estimate of $I_1$.} Since $\phi\in H^2_N(\mathcal I)$, we
have $\partial_y\phi=0$ at $y\in\{0,1\}$, so that
$\int_0^1\partial_y(w\,\partial_y\phi)\,\zeta\,dy=-\int_0^1w\,\partial_y\phi\,\partial_y\zeta\,dy$.
Hence
\begin{equation}\label{eq:I1-split}
I_1\;\leq\;\|w\|_{L^2(\mathcal I)}
+\kappa\sup_{\|\zeta\|_{H^1(\mathcal I)}\leq1}\Big|\int_0^1w\,\partial_y\phi\,\partial_y\zeta\,dy\Big|
+\sup_{\|\zeta\|_{H^1(\mathcal I)}\leq1}\Big|\int_0^1B(y,s)\,\partial_yw\,\zeta\,dy\Big| .
\end{equation}
For the second term, by Cauchy--Schwarz, $\|\partial_y\zeta\|_{L^2(\mathcal I)}\leq1$
and the interpolation inequality,
\[
\kappa\sup_{\|\zeta\|_{H^1(\mathcal I)}\leq1}\Big|\int_0^1w\,\partial_y\phi\,\partial_y\zeta\,dy\Big|
\leq\kappa\,\|w\|_{L^\infty(\mathcal I)}\|\partial_y\phi\|_{L^2(\mathcal I)}
\leq\kappa C_1\|\partial_y\phi\|_{L^2(\mathcal I)}\,\|w\|_{L^2(\mathcal I)}^{1/2}\|w\|_{H^1(\mathcal I)}^{1/2} .
\]
The third term is estimated by \eqref{eq:aux-B} with $u=w$. Fix
$\eta>0$ and apply \eqref{eq:Young} with $\varepsilon=\eta/2$ to each of
the two products
$\|w\|_{L^2(\mathcal I)}^{1/2}\|w\|_{H^1(\mathcal I)}^{1/2}$ (with
$A=\kappa C_1\|\partial_y\phi\|_{L^2(\mathcal I)}$ and $A=2C_1C_2\,m(s)$,
respectively), so that their contributions sum to exactly
$\eta\|w\|_{H^1(\mathcal I)}$. We obtain
\begin{equation}\label{eq:I1-final-n}
I_1\;\leq\;\eta\,\|w\|_{H^1(\mathcal I)}
+\Big[1+\frac{\kappa^2C_1^2}{2\eta}\|\partial_y\phi\|_{L^2(\mathcal I)}^2+N_\eta(s)\Big]\|w\|_{L^2(\mathcal I)},
\end{equation}
where
\begin{equation*}
N_\eta(s):=\frac{2C_1^2C_2^2}{\eta}\,m(s)^2+2\,m(s) .
\end{equation*}

\smallskip
\noindent{\bf Estimate of $I_2$.} Since
$\|A_{b,1}\xi\|_{L^2(\mathcal I)}\leq c_b$ whenever
$\|\xi\|_{H^2_N(\mathcal I)}\leq1$, the Cauchy--Schwarz inequality and
\eqref{eq:B-bounds} give directly, without any need to absorb an
$\eta$-weighted term,
\begin{equation}\label{eq:I2-final}
I_2\;\leq\;c_b\Big[\big(1+d\,\gamma^2(t(s))\big)\|\phi\|_{L^2(\mathcal I)}
+c\,\gamma^2(t(s))\,\|w\|_{L^2(\mathcal I)}
+m(s)\,\|\partial_y\phi\|_{L^2(\mathcal I)}\Big] .
\end{equation}

\smallskip
\noindent\emph{Conclusion.} Combining \eqref{eq:I1-final-n} and
\eqref{eq:I2-final}, and using
\[
\|w\|_{H^1(\mathcal I)}\leq\|(w,\phi)\|_{\mathcal V},
\qquad
\|w\|_{L^2(\mathcal I)},\ \|\phi\|_{L^2(\mathcal I)},\ \|\partial_y\phi\|_{L^2(\mathcal I)}\leq\|(w,\phi)\|_{\mathcal H},
\]
we obtain
\begin{equation*}
\|F(s,w,\phi)\|_{\mathcal V'}\;\leq\;\eta\,\|(w,\phi)\|_{\mathcal V}+\Psi_\eta\big(s,\|(w,\phi)\|_{\mathcal H}\big),
\end{equation*}
where, for $s\in\mathbb R$ and $r\geq0$,
\begin{equation*}
\begin{array}{l}
\displaystyle \Psi_\eta(s,r):=\frac{\kappa^2C_1^2}{2\eta}\,r^3
+\big[1+N_\eta(s)+K(s)\big]\,r,\\[10pt] \displaystyle
K(s):=c_b\Big[1+(c+d)\,\gamma^2(t(s))+m(s)\Big] .
\end{array}
\end{equation*}
Since $m(\cdot)$ and $\gamma^2(t(\cdot))$ are continuous, $\Psi_\eta$ is
continuous on $\mathbb R\times[0,\infty)$, and $r\mapsto\Psi_\eta(s,r)$
is increasing for each fixed $s$ (as a sum of nonnegative,
increasing powers of $r$). Finally, since
$\|z\|_{\mathcal H}\leq|z|/\min\{1,\sqrt b\}$, the function
$\widetilde\Psi_\eta(s,r):=\Psi_\eta\big(s,r/\min\{1,\sqrt b\}\big)$
has the same properties and
$\|F(s,z)\|_{\mathcal V'}\leq\eta\|z\|_{\mathcal V}+\widetilde\Psi_\eta(s,|z|)$.
This verifies \textbf{(F.i)} with $U\equiv0$ (so that the term
$c_\eta\|U(s)\|_{\mathcal V'}$ is absent).

\smallskip
\noindent{\bf Verification of (F.ii).}
Let $z_1=(w,\phi)$, $z_2=(v,\psi)\in\mathcal V$. The terms of
$F(s,z_1)-F(s,z_2)$ that are linear in $z$ are estimated exactly as
above: the first component contributes $\|w-v\|_{L^2(\mathcal I)}$ and
the second component contributes, as in \eqref{eq:I2-final},
\[
c_b\Big[\big(1+d\,\gamma^2(t(s))\big)\|\phi-\psi\|_{L^2(\mathcal I)}
+c\,\gamma^2(t(s))\,\|w-v\|_{L^2(\mathcal I)}
+m(s)\,\|\partial_y(\phi-\psi)\|_{L^2(\mathcal I)}\Big]
\leq K(s)\,\|z_1-z_2\|_{\mathcal H} .
\]
It remains to estimate, in $(H^1(\mathcal I))'$, the two terms
\[
\kappa\big[\partial_y(w\,\partial_y\phi)-\partial_y(v\,\partial_y\psi)\big]
\quad\text{and}\quad
B(y,s)\,\partial_y(w-v) .
\]
For the first, writing
$w\,\partial_y\phi-v\,\partial_y\psi=(w-v)\,\partial_y\phi+v\,\partial_y(\phi-\psi)$,
integrating by parts against $\zeta\in H^1(\mathcal I)$ (the boundary
terms vanish since $\phi,\psi\in H^2_N(\mathcal I)$), and applying
Cauchy--Schwarz together with the interpolation inequality (to $w-v$ and
to $\partial_y(\phi-\psi)$),
\begin{align*}
\kappa\sup_{\|\zeta\|_{H^1(\mathcal I)}\leq1}&\Big|\int_0^1\big[(w-v)\,\partial_y\phi+v\,\partial_y(\phi-\psi)\big]\,\partial_y\zeta\,dy\Big|\\
&\leq\kappa\,\|w-v\|_{L^\infty(\mathcal I)}\|\partial_y\phi\|_{L^2(\mathcal I)}
+\kappa\,\|v\|_{L^2(\mathcal I)}\|\partial_y(\phi-\psi)\|_{L^\infty(\mathcal I)}\\
&\leq\kappa C_1\|\partial_y\phi\|_{L^2(\mathcal I)}\,\|w-v\|_{L^2(\mathcal I)}^{1/2}\|w-v\|_{H^1(\mathcal I)}^{1/2}
+\kappa C_1\|v\|_{L^2(\mathcal I)}\,\|\phi-\psi\|_{H^1(\mathcal I)}^{1/2}\|\phi-\psi\|_{H^2(\mathcal I)}^{1/2} .
\end{align*}
For the second, \eqref{eq:aux-B} with $u=w-v$ gives
\[
\sup_{\|\zeta\|_{H^1(\mathcal I)}\leq1}\Big|\int_0^1B(y,s)\,\partial_y(w-v)\,\zeta\,dy\Big|
\leq2C_1C_2\,m(s)\,\|w-v\|_{L^2(\mathcal I)}^{1/2}\|w-v\|_{H^1(\mathcal I)}^{1/2}
+2\,m(s)\,\|w-v\|_{L^2(\mathcal I)} .
\]
We now apply \eqref{eq:Young} with $\varepsilon=\eta/3$ to each of the
three products above: with $u=w-v$ in the first and third ones, and,
in the second one, to the product
$\|\phi-\psi\|_{H^1(\mathcal I)}^{1/2}\|\phi-\psi\|_{H^2(\mathcal I)}^{1/2}$,
for which \eqref{eq:Young} holds verbatim with the pair
$(L^2,H^1)$ replaced by $(H^1,H^2)$. The terms carrying $\mathcal V$-norms add up to
\[
\frac{2\eta}3\,\|w-v\|_{H^1(\mathcal I)}+\frac\eta3\,\|\phi-\psi\|_{H^2(\mathcal I)}
\leq\eta\,\|z_1-z_2\|_{\mathcal V},
\]
while the remaining terms are bounded by
\[
\Big[\frac{3\kappa^2C_1^2}{4\eta}\Big(\|\partial_y\phi\|_{L^2(\mathcal I)}^2+\|v\|_{L^2(\mathcal I)}^2\Big)
+\frac{3C_1^2C_2^2}{\eta}\,m(s)^2+2\,m(s)\Big]\,\|z_1-z_2\|_{\mathcal H} .
\]
Collecting all terms, and using
$\|\partial_y\phi\|_{L^2(\mathcal I)}^2+\|v\|_{L^2(\mathcal I)}^2\leq\big(\|z_1\|_{\mathcal H}+\|z_2\|_{\mathcal H}\big)^2$,
\begin{equation*}
\|F(s,z_1)-F(s,z_2)\|_{\mathcal V'}
\;\leq\;\eta\,\|z_1-z_2\|_{\mathcal V}
+\Phi_\eta\big(s,\|z_1\|_{\mathcal H}+\|z_2\|_{\mathcal H}\big)\,\|z_1-z_2\|_{\mathcal H},
\end{equation*}
where, for $s\in\mathbb R$ and $r\geq0$,
\begin{equation*}
\Phi_\eta(s,r):=\frac{3\kappa^2C_1^2}{4\eta}\,r^2+M_\eta(s),
\qquad
M_\eta(s):=1+\frac{3C_1^2C_2^2}{\eta}\,m(s)^2+2\,m(s)+K(s) .
\end{equation*}
As before, $\Phi_\eta$ is continuous on $\mathbb R\times[0,\infty)$ and
increasing in $r$. Setting
\[
\widetilde\Phi_\eta(s,r):=\frac{1}{\min\{1,\sqrt b\}}\,\Phi_\eta\Big(s,\frac{r}{\min\{1,\sqrt b\}}\Big)
\]
and using $1\leq1+\|z_1\|_{\mathcal V}+\|z_2\|_{\mathcal V}$, we obtain
\[
\|F(s,z_1)-F(s,z_2)\|_{\mathcal V'}
\leq\eta\,\|z_1-z_2\|_{\mathcal V}
+\big(1+\|z_1\|_{\mathcal V}+\|z_2\|_{\mathcal V}\big)\,
\widetilde\Phi_\eta\big(s,|z_1|+|z_2|\big)\,|z_1-z_2| ,
\]
which is \textbf{(F.ii)}.

Therefore, by Theorem~\ref{T-Exis}, for every $\tau\in\mathbb R$ and
$(w_0,\phi_0)\in\mathcal H$ there exist $T_{w_0,\phi_0}=T_0(\tau,w_0,\phi_0)>\tau$
and a unique local weak solution to \eqref{equ4iii},
\[
(w,\phi)\in H^1(\tau,T_{w_0,\phi_0};\mathcal V')\cap \mathcal C([\tau,T_{w_0,\phi_0}];\mathcal H)\cap L^2(\tau,T_{w_0,\phi_0};\mathcal V) ,
\]
which, by the identification $\mathcal V'\simeq(H^1(\mathcal I))'\times L^2(\mathcal I)$,
is precisely the assertion of the theorem.
\end{proof}

\begin{theorem}
{\it (Existence of a regular local weak solution to \eqref{equ4iii})} Under the assumptions given in (H1)-(H2), given $(w_0,\phi_0) \in H^1(\mathcal{I})\times H^{2}_N(\mathcal{I})$ and $\tau \in \mathbb{R}$, then there exists $T_0=T_0(\tau, w_0,\rho_0)> \tau$ such that the semilinear  parabolic evolution problem \eqref{equ4iii} has a unique local weak solution $(w(\cdot),\phi(\cdot)): [\tau,T_{w_0,\phi_0}) \to H^1(\mathcal{I})\times H^{2}_N(\mathcal{I})$ with 
\begin{equation}\label{reg-33}
\begin{array}{l}
  \displaystyle w \in H^1 (\tau , T_{w_0,\phi_0} ; L^{2}(\mathcal{I}))\cap \mathcal{C}([\tau, T_{w_0,\phi_0}]; H^1(\mathcal{I})) \cap L^2 (\tau,T_{w_0,\phi_0} ; H^{2}_N(\mathcal{I})) \\\\
   \displaystyle\phi \in H^1 (\tau , T_{w_0,\phi_0} ; H^{1}(\mathcal{I}))\cap \mathcal{C}([\tau, T_{w_0,\phi_0}]; H^2_N(\mathcal{I})) \cap L^2 (\tau,T_{w_0,\phi_0} ; H^{3}_N(\mathcal{I})),
   \end{array}
\end{equation}
 satisfying the initial condition $(w(\tau ),\phi( \tau ))=(w_0,\phi_0)$.
\end{theorem}
\begin{proof} Consider the Hilbert space
\[
\mathcal{H} =  H^1(\mathcal{I})\times H^{2}_N(\mathcal{I}), \hspace{0.3cm} \mathcal{V}= H^2_N(\mathcal I)\times H^3_N(\mathcal I).
\]
Let $\widetilde{A}_{a,1}:H^1(\mathcal{I}) \to H^{-1}(\mathcal{I})$ be the operator associated with $-a\Delta + I$ (see Example 01).
Thus, consider $A_{1}:H^2_{N}(\mathcal{I}) \to L^2(\mathcal{I})$ the $L^2$-realization of $\widetilde{A}_{1,1}$, and $A_{2}:H^3_{N}(\mathcal{I}) \to H^1(\mathcal{I})$ the $H^1$-realization of $\widetilde{A}_{b,1}$ (see Example 02). Then, the dual space of $\mathcal{V}$ is identify by
\[
\mathcal{V}' = L^{2}(\mathcal{I}) \times H^1(\mathcal{I}) .
\]
Thus, the dual product between $\mathcal{V}'$ and $\mathcal{V}$ is given by
\begin{align*}
    \Big\langle \begin{pmatrix}
        f\\ \phi 
    \end{pmatrix} , \begin{pmatrix}
        u\\ \rho 
    \end{pmatrix} \Big \rangle_{\mathcal{V}', \mathcal{V}} = ( f, A_{1} u )_{L^2(\mathcal I) } + ( A_{2}^{3/2}\phi , A_{2}^{3/2} \rho)_{L^2 (\mathcal I) }. 
\end{align*}
Now, we define the operator
\[
\mathbf{A} :=
\begin{pmatrix}
A_1 & 0 \\
0 & A_2
\end{pmatrix}.
\]
Note that $\mathbf{A}$ is a bounded linear operator from $\mathcal{V}$ to $\mathcal{V}'$ associated with the bilinear form
\begin{equation*}
   \mathbf{a}_{\mathbf{A}}(U,V)= (A_1 u,A_1 v)_{L^2(\mathcal{I})} + (A_2^{3/2}\rho , A_2^{3/2}\phi)_{L^2(\mathcal{I})}, \hspace{0,3cm} U=(u,\rho), V=(v,\phi) \in \mathcal{V} 
\end{equation*}
and satisfies the conditions \eqref{ai}-\eqref{aii}. Then, the abstract formulation of the system \eqref{equ4iii} is given as
\begin{equation*}
  \begin{cases}
 \displaystyle z' + \mathbf{A}z = F(s,z), \ \ \  s> \tau;\\
z(\tau)=\begin{pmatrix}
    w_{0}\\  \phi_0
\end{pmatrix} \in \mathcal{H},
  \end{cases}
\end{equation*}
where  
\begin{align*}
  z=\begin{pmatrix}
    w \\ \phi
\end{pmatrix}, \hspace{0,5cm}  F(s,z)=\begin{pmatrix} w-\kappa  \frac{\partial }{\partial y}\Big(w \frac{\partial \phi}{\partial y}\Big) - \gamma^{2}(t(s))b(y,t(s)) \frac{\partial w}{\partial y} \\ 
       \phi +  c\gamma^{2}(t(s)) w -d\gamma^{2}(t(s))\phi- \gamma^{2}(t(s))b(y,t(s)) \frac{\partial \phi}{\partial y}
    \end{pmatrix}.
\end{align*}

In order to apply Theorem \ref{T-Exis}, we can follow the same way as Theorem \ref{teoEx-02} and \cite[Theorem 3.3]{OY01}. Therefore, it follows from Theorem \ref{T-Exis} that given $\tau \in \mathbb{R}$ and $(w_0,\phi_0) \in \mathcal{H}$, then there exist $T_{w_0,\phi_0}=T_0(\tau, w_0,\phi_0)> \tau$ and a unique local weak solution to \eqref{equ4iii} with regularity given in \eqref{reg-33}.
\end{proof}

\begin{theorem}\label{teoEx-04}
({\it Extra regularity}) Under the assumptions given in (H1)-(H2), given $(w_0,\phi_0) \in H^2_N(\mathcal{I})\times H^{3}_N(\mathcal{I})$ and $\tau \in \mathbb{R}$, then there exists $T_{w_0,\phi_0}=T_0(\tau, w_0,\phi_0)> \tau$ such that the semilinear  parabolic evolution problem \eqref{equ4iii} has a unique local weak solution $(w(\cdot),\phi(\cdot)): [\tau,T_{w_0,\phi_0}) \to H^2_N(\mathcal{I})\times H^{3}_N(\mathcal{I})$ with 
\begin{equation*}
\begin{array}{l}
  \displaystyle w \in H^1 (\tau , T_{w_0,\phi_0} ; H^{1}(\mathcal{I}))\cap \mathcal{C}([\tau, T_{w_0,\phi_0}]; H^2_N(\mathcal{I})) \cap L^2 (\tau,T_{w_0,\phi_0} ; H^{3}_N(\mathcal{I})) \\\\
   \displaystyle\phi \in H^1 (\tau , T_{w_0,\phi_0} ; H^{2}_N(\mathcal{I}))\cap \mathcal{C}([\tau, T_{w_0,\phi_0}]; H^3_N(\mathcal{I})) \cap L^2 (\tau,T_{w_0,\phi_0} ; H^{4}_{N^2}(\mathcal{I})),
   \end{array}
\end{equation*}
 satisfying the initial condition $(w(\tau ),\phi( \tau ))=(w_0,\phi_0)$.
\end{theorem}
\begin{proof}
The proof follows the same lines as the proof of Theorem \ref{teoEx-02} (see also
\cite[Theorem 3.3]{OY01}), now taking
$$H=H^2_N(\mathcal{I})\times H^3_N(\mathcal{I}),\qquad V=H^3_N(\mathcal{I})\times H^4_{N_2}(\mathcal{I}),$$
with $A=\mathrm{diag}(A_1,A_2)$ built from $\tilde A_{a,1}$ and $\tilde
A_{b,1}$, as in Theorem \ref{teoEx-02}, at one order of regularity higher (i.e., $A_1$ is
now the $H^1$-realization of $\tilde A_{a,1}$ on $H^3_N(\mathcal{I})$, and $A_2$ is the
$H^2_N$-realization of $\tilde A_{b,1}$ on $H^4_{N_2}(\mathcal{I})$). The verification of (F.i)-(F.ii) is obtained through the same interpolation and Cauchy--Schwarz
estimates as in Theorem \ref{teoEx-02}, replacing the embedding $H^1(\mathcal{I})\hookrightarrow
L^\infty(\mathcal{I})$ by $H^2_N(\mathcal{ I})\hookrightarrow C^1( \mathcal{\overline I})$ to control the nonlinear term $\partial_y(w\,\partial_y\phi)$. Theorem \ref{T-Exis} then yields, for every $\tau\in\mathbb{R}$ and $(w_0,\phi_0)\in H^2_N(\mathcal{I})\times H^3_N(\mathcal{I})$, a unique
local weak solution with the regularity stated above.
\end{proof}
As a consequence of the local existence theorems for problem \eqref{equ4iii}, and by means of the mapping $r(x,t)$ defined in \eqref{kappa}, the rescaling function $s(\cdot)$ given in \eqref{rescaling}, Proposition \ref{Prop1} and Proposition \ref{Prop7}, we obtain the existence of a weak solution to problem \eqref{equ4}, which is stated below.

\begin{theorem}[Local weak solutions to \eqref{equ4}]\label{teoEx-05-34}
Assume that {\rm(H1)--(H2)} hold. Let $\tau\in\mathbb{R}$ and let
$(u_0,\rho_0)$ be the initial data. Then there exists
\[
T_{u_0,\rho_0}=T_0(\tau,u_0,\rho_0)>\tau
\]
such that the semilinear parabolic evolution problem \eqref{equ4}
admits a unique local weak solution
$(u(\cdot),\rho(\cdot))$ on $[\tau,T_{u_0,\rho_0}]$. Moreover, the following
regularity properties hold.

\begin{itemize}

\item If
$(u_0,\rho_0)\in L^2(\mathcal{I}_\tau)\times H^1(\mathcal{I}_\tau)$,
then
\begin{equation*}
\begin{array}{l}
\displaystyle
u\in H^1(\tau,T_{u_0,\rho_0};H^1(\mathcal{I}_t)')
\cap C([\tau,T_{u_0,\rho_0}];L^2(\mathcal{I}_t))
\cap L^2(\tau,T_{u_0,\rho_0};H^1(\mathcal{I}_t)),\\\\
\displaystyle
\rho\in H^1(\tau,T_{u_0,\rho_0};L^2(\mathcal{I}_t))
\cap C([\tau,T_{u_0,\rho_0}];H^1(\mathcal{I}_t))
\cap L^2(\tau,T_{u_0,\rho_0};H^2_N(\mathcal{I}_t)).
\end{array}
\end{equation*}

\item If
$(u_0,\rho_0)\in H^1(\mathcal{I}_\tau)\times H^2_N(\mathcal{I}_\tau)$,
then
\begin{equation*}
\begin{array}{l}
\displaystyle
u\in H^1(\tau,T_{u_0,\rho_0};L^2(\mathcal{I}_t))
\cap C([\tau,T_{u_0,\rho_0}];H^1(\mathcal{I}_t))
\cap L^2(\tau,T_{u_0,\rho_0};H^2_N(\mathcal{I}_t)),\\\\
\displaystyle
\rho\in H^1(\tau,T_{u_0,\rho_0};H^1(\mathcal{I}_t))
\cap C([\tau,T_{u_0,\rho_0}];H^2_N(\mathcal{I}_t))
\cap L^2(\tau,T_{u_0,\rho_0};H^3_N(\mathcal{I}_t)).
\end{array}
\end{equation*}

\item If
$(u_0,\rho_0)\in H^2_N(\mathcal{I}_\tau)\times H^3_N(\mathcal{I}_\tau)$,
then
\begin{equation*}
\begin{array}{l}
\displaystyle
u\in H^1(\tau,T_{u_0,\rho_0};H^1(\mathcal{I}_t))
\cap C([\tau,T_{u_0,\rho_0}];H^2_N(\mathcal{I}_t))
\cap L^2(\tau,T_{u_0,\rho_0};H^3_N(\mathcal{I}_t)),\\\\
\displaystyle
\rho\in H^1(\tau,T_{u_0,\rho_0};H^2_N(\mathcal{I}_t))
\cap C([\tau,T_{u_0,\rho_0}];H^3_N(\mathcal{I}_t))
\cap L^2(\tau,T_{u_0,\rho_0};H^4_{N^2}(\mathcal{I}_t)).
\end{array}
\end{equation*}

\end{itemize}
Furthermore, satisfying the initial condition
\[
(u(\tau),\rho(\tau))=(u_0,\rho_0).
\]
\end{theorem}
\subsection{Existence of a global weak solution}

Consider the following interpolation inequality
\begin{equation}\label{LpLq}
    \|\cdot\|_{L^q(\mathcal{I})} \leq C_{p,q}\|\cdot\|_{H^1 (\mathcal{I})}^{\sigma} \|\cdot\|_{L^p(\mathcal{I})}^{1-\sigma}, \hspace{0,5cm} 1\leq p< q\leq \infty,
\end{equation}
where $\sigma=(1/p -1/q)/(1/p +1/2)$.

For $0\leq s_0 <s <s_1 <\infty$, $H^s(\mathcal{I})=[H^{s_0}(\mathcal{I}), H^{s_1}(\mathcal{I})]_{\theta}$ with $s=(1-\theta)s_0 + \theta s_1$, and
\begin{equation}\label{Hs}
    \|\cdot\|_{H^s(\mathcal{I})} \leq C_{s_0 s_1} \|\cdot\|_{H^{s_0}(\mathcal{I})}^{1-\theta} \|\cdot\|_{H^{s_1}(\mathcal{I})}^{\theta}.
\end{equation}

For $s>1/2$, $H^s(\mathcal{I}) \subset \mathcal{C}(\overline{\mathcal{I}})$ with
\begin{equation}\label{CH1}
    \|\cdot\|_{\mathcal{C}(\overline{\mathcal{I}})} \leq C_s \|\cdot\|_{H^s(\mathcal{I})}.
\end{equation}

The next result proves that solutions to problem \eqref{equ4iii} remain non negative. This property is essential for the global well-posedness analysis, since it allows us to exclude the possibility of finite time blow up.

\begin{lemma}\label{lem:positivity}
Assume \textup{(H1)--(H2)}, and let $L^{2}(\mathcal{I})\ni w_0 \geq0$ and
$H^{1}(\mathcal{I})\ni \phi_0 \geq0$ \textup{(}a.e.\ in $\mathcal I$\textup{)}.
Then the unique weak solution $(w,\phi)$ to \eqref{equ4iii} on the
interval $[\tau,T_{w_0,\phi_0})$, given in Theorem~\ref{teoEx-02},
satisfies
\begin{equation}\label{equ38}
w(\cdot,s)\geq0 \quad\text{and}\quad \phi(\cdot,s)\geq0
\quad\text{a.e. in }\mathcal I,\ \text{ for every } s\in[\tau,T_{w_0,\phi_0}).
\end{equation}
\end{lemma}

\begin{proof}
Throughout the proof we use the notation of the proof of
Theorem~\ref{teoEx-02}; in particular
$B(y,s)=\gamma^2(t(s))\,b(y,t(s))=-\gamma(t(s))\big(\alpha'(t(s))+\gamma'(t(s))\,y\big)$,
so that, since $\alpha+\gamma=\beta$,
\begin{equation}\label{eq:B-boundary}
B(0,s)=-\gamma(t(s))\,\alpha'(t(s)),\qquad
B(1,s)=-\gamma(t(s))\,\beta'(t(s)),\qquad
\partial_yB(y,s)=-\gamma(t(s))\,\gamma'(t(s)) .
\end{equation}
For $T>\tau$ we set
\[
\ell_T:=\max_{s\in[\tau,T]}\gamma(t(s))\,|\gamma'(t(s))|
\;\leq\;\gamma_1\max_{r\in[\tau,t(T)]}|\gamma'(r)|<\infty ,
\]
which is finite by (H1), since $\gamma'$ is continuous; under (H3) one
may simply take $\ell_T=\gamma_1\|\gamma'\|_{L^\infty(\mathbb R)}$.

Let $H:\mathbb{R} \to [0,\infty)$ be a nonincreasing $\mathcal{C}^3$-function
such that $H(r)>0$ if $r<0$ and $H(r)=0$ if $r\geq 0$, and such that,
for some constant $c_H>0$,
\begin{equation}\label{eq:H-props}
\begin{cases}
0\leq H(r) \leq c_H\,r^2, & r\in \mathbb{R},\\
0\leq H'(r)\,r \leq c_H\,H(r), & r\in \mathbb{R},\\
0\leq H''(r)\,r^2 \leq c_H\,H(r), & r\in \mathbb{R},\\
\Big| \dfrac{d(H'')^{\frac{1}{2}}}{dr}(r)\,r\Big|\leq c_H\,(H''(r))^{\frac{1}{2}}, & r\in \mathbb{R},
\end{cases}
\end{equation}
see \cite[Chapter 12]{yagi2009} or \cite[Theorem 3.5]{OY01}. Note that
$H''\geq0$, and that $0\leq H''(r)\leq c_H^2$ for all $r$ (first and
third properties, plus continuity at $r=0$); since $H'(0)=0$, the
function $H'$ is globally Lipschitz with $|H'(r)|\leq c_H^2|r|$.

\medskip
\noindent{\bf Step 1: nonnegativity of $w$ for regular data.}
Assume first that $0\leq w_0\in H_N^2(\mathcal{I})$ and
$0\leq\phi_0\in H_N^3(\mathcal{I})$. By Theorem~\ref{teoEx-04}, the
corresponding solution is a strong solution; in particular, for every
$T\in(\tau,T_{w_0,\phi_0})$,
\begin{align*}
w&\in H^1(\tau,T;L^2(\mathcal{I}))\cap L^2(\tau,T;H_N^2(\mathcal{I})),\\
\phi&\in H^1(\tau,T;L^2(\mathcal{I}))\cap \mathcal C([\tau,T];H_N^2(\mathcal{I}))\cap L^2(\tau,T;H_N^3(\mathcal{I})),
\end{align*}
and the equations in \eqref{equ4iii} hold in $L^2(\mathcal I)$ for a.e.\
$s\in(\tau,T)$. Fix such a $T$ and consider
\begin{equation*}
\Xi(s):=\int_{\mathcal{I}}H(w(y,s))\,dy,\qquad s\in[\tau,T].
\end{equation*}
Since $H'$ is globally Lipschitz and $w\in H^1(\tau,T;L^2(\mathcal I))$,
$\Xi$ is absolutely continuous on $[\tau,T]$ and, for a.e.\ $s$,
$\Xi'(s)=\int_{\mathcal I}H'(w)\,\partial_sw\,dy$. Using the first equation
of \eqref{equ4iii},
\[
\partial_sw=a\,\partial_y^2w-\kappa\,\partial_y\big(w\,\partial_y\phi\big)-B(y,s)\,\partial_yw,
\]
and integrating by parts (the boundary terms in the first two integrals
vanish because $\partial_yw=\partial_y\phi=0$ at $y\in\{0,1\}$), we get
\begin{equation}\label{eq:Xi-prime}
\Xi'(s)=-a\int_{\mathcal{I}}H''(w)\,|\partial_yw|^2\,dy
+\kappa\int_{\mathcal{I}}H''(w)\,w\,\partial_yw\,\partial_y\phi\,dy
-\int_{\mathcal{I}}B(y,s)\,H'(w)\,\partial_yw\,dy .
\end{equation}

\emph{The chemotactic term.} Since $H''\geq0$, Young's inequality gives,
for every $\varepsilon>0$,
\[
\kappa\Big|\int_{\mathcal I}H''(w)\,w\,\partial_yw\,\partial_y\phi\,dy\Big|
\leq\kappa\varepsilon\int_{\mathcal I}H''(w)\,|\partial_yw|^2\,dy
+\frac{\kappa}{4\varepsilon}\int_{\mathcal I}H''(w)\,w^2\,|\partial_y\phi|^2\,dy .
\]
Choosing $\varepsilon=a/\kappa$, the first term is absorbed by the
first term in \eqref{eq:Xi-prime}, and by \eqref{eq:H-props},
\begin{align*}
-a\int_{\mathcal{I}}H''(w)\,|\partial_yw|^2\,dy
+\kappa\int_{\mathcal{I}}H''(w)\,w\,\partial_yw\,\partial_y\phi\,dy
&\leq\frac{\kappa^2}{4a}\int_{\mathcal I}H''(w)\,w^2\,|\partial_y\phi|^2\,dy\\
&\leq\frac{\kappa^2c_H}{4a}\,\|\partial_y\phi(s)\|_{L^\infty(\mathcal I)}^2\,\Xi(s).
\end{align*}

\emph{The transport term.} Since $H'(w)\,\partial_yw=\partial_y\big(H(w)\big)$,
integrating by parts and using \eqref{eq:B-boundary},
\begin{align*}
-\int_{\mathcal I}B(y,s)\,H'(w)\,\partial_yw\,dy
&=-\big[B(y,s)\,H(w(y,s))\big]_{y=0}^{y=1}
+\int_{\mathcal I}\partial_yB(y,s)\,H(w)\,dy\\
&=\gamma(t(s))\,\beta'(t(s))\,H(w(1,s))-\gamma(t(s))\,\alpha'(t(s))\,H(w(0,s))\\
&\qquad-\gamma(t(s))\,\gamma'(t(s))\,\Xi(s) .
\end{align*}
By (H2), $\alpha'(t(s))\geq0$ and $\beta'(t(s))\leq0$; since $H\geq0$,
the boundary contribution is nonpositive. Hence
\[
-\int_{\mathcal I}B(y,s)\,H'(w)\,\partial_yw\,dy\leq\gamma(t(s))\,|\gamma'(t(s))|\,\Xi(s)\leq\ell_T\,\Xi(s).
\]

\emph{Conclusion of Step 1.} Since $\phi\in\mathcal C([\tau,T];H^2_N(\mathcal I))$
and $H^2(\mathcal I)\hookrightarrow C^1(\overline{\mathcal I})$, the
quantity
\[
C_T:=\frac{\kappa^2c_H}{4a}\,\sup_{s\in[\tau,T]}\|\partial_y\phi(s)\|_{L^\infty(\mathcal I)}^2+\ell_T
\]
is finite, and we have shown that
\[
\Xi'(s)\leq C_T\,\Xi(s)\qquad\text{for a.e. } s\in(\tau,T).
\]
Since $w_0\geq0$, $\Xi(\tau)=0$, and Gronwall's inequality gives
$\Xi\equiv0$ on $[\tau,T]$, that is, $H(w(\cdot,s))=0$ and hence
$w(\cdot,s)\geq0$ a.e.\ in $\mathcal I$, for every $s\in[\tau,T]$. As
$T\in(\tau,T_{w_0,\phi_0})$ is arbitrary, $w\geq0$ on
$[\tau,T_{w_0,\phi_0})$.

\medskip
\noindent{\bf Step 2: nonnegativity of $\phi$ for regular data.}
With the same data and $T$ as in Step 1, let
$\Xi_2(s):=\int_{\mathcal I}H(\phi(y,s))\,dy$. Using the second equation
of \eqref{equ4iii} and $\partial_y\phi=0$ at $y\in\{0,1\}$,
\begin{align*}
\Xi_2'(s)={}&-b\int_{\mathcal I}H''(\phi)\,|\partial_y\phi|^2\,dy
+c\,\gamma^2(t(s))\int_{\mathcal I}H'(\phi)\,w\,dy\\
&-d\,\gamma^2(t(s))\int_{\mathcal I}H'(\phi)\,\phi\,dy
-\int_{\mathcal I}B(y,s)\,H'(\phi)\,\partial_y\phi\,dy .
\end{align*}
The first term is nonpositive since $H''\geq0$. The second one is
nonpositive because $H'\leq0$ and, by Step 1, $w\geq0$ (this is why $w$
is treated first). The third one is nonpositive since $H'(r)\,r\geq0$
by \eqref{eq:H-props} and $d>0$. The last one is estimated exactly as the
transport term in Step 1, with $w$ replaced by $\phi$. Therefore
$\Xi_2'(s)\leq\ell_T\,\Xi_2(s)$ for a.e.\ $s\in(\tau,T)$, and, since
$\Xi_2(\tau)=0$, Gronwall's inequality yields $\phi\geq0$ on
$[\tau,T]$, for every $T\in(\tau,T_{w_0,\phi_0})$.

\medskip
\noindent{\bf Step 3: general data.}
Now, given initial data conditions $0 \leq w_0 \in L^2(\mathcal{I})$ and
$0 \leq \phi_0 \in H^1(\mathcal{I})$, we take sequences $\{w_{0,n}\}$ and $\{\phi_{0,n}\}$ such that
\[
0 \leq w_{0,n} \in H_N^2(\mathcal{I})
\quad\text{and}\quad
0 \leq \phi_{0,n} \in H^3_N(\mathcal{I}),
\]
which converge to $w_0$ in $L^2(\mathcal{I})$ and to $\phi_0$ in $H^1(\mathcal{I})$,
respectively.

We already know that, for each initial data $(w_{0,n},\phi_{0,n})$, there exists a solution
$(w_n,\phi_n)$ to \eqref{equ4iii} satisfying \eqref{equ38}
on some interval $[\tau,T_{w_{0,n},\phi_{0,n}})$.  Since
\[
w_{0,n}\to w_0 \quad\text{in }L^2(\mathcal I),
\qquad
\phi_{0,n}\to\phi_0 \quad\text{in }H^1(\mathcal I),
\]
there exists a constant $M>0$ such that
\[
\|w_{0,n}\|_{L^2(\mathcal I)}
+\|\phi_{0,n}\|_{H^1(\mathcal I)}
\le M,
\qquad\forall n.
\]
By the local existence Theorem \ref{T-Exis}, the existence time depends only on this bound. Hence there exists a positive constant $T(M)>0$ such that
\[
T_{w_{0,n},\phi_{0,n}}\ge T(M),
\qquad\forall n.
\]
Therefore,
\[
T_0:=\inf_nT_{w_{0,n},\phi_{0,n}}
\ge T(M)>\tau.
\]

On the other hand, the solution depends continuously on the
initial data. Therefore, as $n\to\infty$,
\[
w_n \longrightarrow w
\quad\text{in}\quad
\mathcal{C}([\tau,T_0);L^2(\mathcal{I})),
\]
and
\[
\phi_n \longrightarrow \phi
\quad\text{in}\quad
\mathcal{C}([\tau,T_0);H^1(\mathcal{I})).
\]
Hence, $(w,\phi)$ satisfies estimate \eqref{equ38} for
$\tau\leq s\leq T_0$. Finally, by the uniqueness of the local
solution, estimate \eqref{equ38} remains valid for every
$s\in[\tau,T_{w_0,\phi_0})$.
\end{proof}

\begin{corollary}
Given $\tau\in \mathbb{R}$, $L^{2}(\mathcal{I}_\tau)\ni u_0 \geq0$ and $H^{1}(\mathcal{I}_\tau)\ni \rho_0 \geq0$. Then, the unique weak solutions $u$ and $\rho$ to \eqref{equ4} on the interval $[0,T_{u_0,\rho_0}]$, given in Theorem \ref{teoEx-05-34}, are such that
\begin{equation*}
u(x,t)\geq0, \hspace{0,3cm} \text{ and } \hspace{0,3cm} \rho(x,t)\geq0,
\quad \text{ for all }
x \in \mathcal{I}_t \ \text{ with } t\in [\tau,T_{u_0,\rho_0}).
\end{equation*}
\end{corollary}
\begin{proof}
   Combining the change of variables $r(\cdot,t)$ defined in \eqref{kappa}, the time rescaling \eqref{rescaling}, and Lemma \ref{lem:positivity}, it follows that
\[
u(x,t)\geq0,\qquad \rho(x,t)\geq0,
\]
for all $x\in\mathcal{I}_t$ and $t\in[\tau,T_{u_0,\rho_0})$.
\end{proof}

\begin{theorem}\label{teoEx-05} {\it (Existence of global weak solution to \eqref{equ4iii})} Under the assumptions given in (H1)-(H2), given $(w_0,\varphi_0) \in L^2(\mathcal{I})\times H^{1}(\mathcal{I})$, with $w_0\geq 0$, $\varphi_0 \geq 0$ and $T, \tau \in \mathbb{R}$ with $T>\tau$, then the semilinear  parabolic evolution problem \eqref{equ4iii} has a unique global weak solution $(w(\cdot),\varphi(\cdot))=(w(\cdot; \tau ,w_0,\varphi_0),\varphi(\cdot; \tau ,w_0,\varphi_0)): [\tau,T] \to L^2(\mathcal{I})\times H^{1}(\mathcal{I})$ with
\begin{align*}
   &0\leq w \in H^1 \big(\tau , T ; H^{1}(\mathcal{I})'\big)\cap \mathcal{C}([\tau, T]; L^2(\mathcal{I})) \cap L^2 (\tau,T ; H^{1}(\mathcal{I})) \\&
  0\leq  \varphi \in H^1 (\tau , T ; L^{2}(\mathcal{I}))\cap \mathcal{C}([\tau, T]; H^1(\mathcal{I})) \cap L^2 (\tau,T ; H^{2}_N(\mathcal{I})),
\end{align*}
 satisfying the initial condition $(w(\tau ),\varphi( \tau ))=(w_0,\varphi_0)$.
\end{theorem}
\begin{proof}
We will make formal computations that can be justified by using the approximations of the solutions $w$ and $\varphi$ of \eqref{equ4iii} associated with regular initial data as in Theorem \ref{teoEx-04}.

Fix $T>\tau$. Throughout the proof, $s\in[\tau,T]\cap[\tau,T_{w_0,\varphi_0})$, where $[\tau,T_{w_0,\varphi_0})$ is the maximal interval of existence, and we write
\[
B(y,s):=\gamma^2(t(s))\,b(y,t(s))=-\gamma(t(s))\big(\alpha'(t(s))+\gamma'(t(s))\,y\big),
\]
\[
m_T:=\max_{s\in[\tau,T]}\gamma(t(s))\big(|\alpha'(t(s))|+|\gamma'(t(s))|\big),
\]
so that $|B(y,s)|\leq m_T$ for $(y,s)\in\overline{\mathcal I}\times[\tau,T]$; $m_T<\infty$ by (H1), since $\alpha'$ and $\gamma'$ are continuous. (Under (H3) one may take $m_T=\gamma_1(\|\alpha'\|_{\infty}+\|\gamma'\|_{\infty})$.)

Now, integrating on $\mathcal{I}$ the first equation of \eqref{equ4iii}, we have
\begin{align*}
   \frac{d}{ds} \int_{\mathcal{I}}w\,dy &= a\int_{\mathcal{I}} \frac{\partial^{2}w}{\partial y^{2}}dy - \kappa\int_{\mathcal{I}}\frac{\partial }{\partial y}\Big(w \frac{\partial \varphi}{\partial y}\Big)dy - \gamma^{2}(t(s))\int_{\mathcal{I}}b(y,t(s)) \frac{\partial w}{\partial y}dy\\
   & = a\frac{\partial w}{\partial y}\Big|_{y=0}^{y=1} - \kappa\, w \frac{\partial \varphi}{\partial y}\Big|_{y=0}^{y=1} - \gamma^{2}(t(s))\int_{\mathcal{I}}b(y,t(s)) \frac{\partial w}{\partial y}dy\\
   &= -\gamma^{2}(t(s))\int_{\mathcal{I}}b(y,t(s)) \frac{\partial w}{\partial y}dy\\
   & = \gamma(t(s))\int_{\mathcal{I}}\big(\alpha'(t(s))+y\,\gamma'(t(s))\big)\frac{\partial w}{\partial y}dy\\
   &= \gamma(t(s))\alpha'(t(s))[w(1,s)-w(0,s)]+ \gamma(t(s))\gamma'(t(s))\int_{\mathcal{I}}y\frac{\partial w}{\partial y}dy \\
   &= \gamma(t(s))\alpha'(t(s))[w(1,s)-w(0,s)]+ \gamma(t(s))\gamma'(t(s))w(1,s) - \gamma(t(s))\gamma'(t(s))\int_{\mathcal{I}}w\, dy.
\end{align*}
Since $t'(s)=\gamma^2(t(s))$, we have $\frac{d}{ds}\gamma(t(s))=\gamma^2(t(s))\gamma'(t(s))$. Then, using $\alpha'+\gamma'=\beta'$,
\begin{align*}
   \frac{d}{ds}\Big(\gamma(t(s))\int_{\mathcal{I}}w(y,s)\,dy\Big) &= \gamma(t(s))\Big[\frac{d}{ds} \int_{\mathcal{I}}w\,dy+\gamma(t(s))\gamma'(t(s))\int_{\mathcal{I}}w\,dy\Big]\\
   &=\gamma(t(s))\Big[\gamma(t(s))\alpha'(t(s))[w(1,s)-w(0,s)]+ \gamma(t(s))\gamma'(t(s))w(1,s)\Big] \\
   & = \gamma^2(t(s))\Big[ \beta'(t(s))w(1,s) -\alpha'(t(s))w(0,s) \Big].
\end{align*}
It follows from (H2) and Lemma \ref{lem:positivity} that
\begin{align}\label{w-gamma}
   \frac{d}{ds}\Big(\gamma(t(s))\int_{\mathcal{I}}w(y,s)\,dy\Big) \leq 0.
\end{align}
Therefore, since $t(\tau)=\tau$ and $w\geq0$, we get
\begin{equation}\label{w-L1-00}
    \|w(s)\|_{L^1(\mathcal{I})}\leq \frac{\gamma(t(\tau))}{\gamma(t(s))}\|w_0\|_{L^1(\mathcal{I})}\leq C_{\gamma}\|w_0\|_{L^1(\mathcal{I})},
\end{equation}
for all $ s\in [\tau,T_{w_0,\varphi_0})$, where $C_{\gamma}:=\gamma_1/\gamma_0$.

Now, multiplying the second equation of \eqref{equ4iii} by $\varphi$ and integrating on $\mathcal{I}$, we obtain
\begin{align*}
\begin{array}{l}
  \displaystyle  \frac{1}{2}\frac{d}{ds}\|\varphi\|_{L^2(\mathcal{I})}^{2} + b\Big\|\frac{ \partial \varphi}{\partial y}\Big\|_{L^2(\mathcal{I})}^{2} = c\gamma^2(t(s))\int_{\mathcal{I}} w\varphi\, dy - d\gamma^{2}(t(s)) \|\varphi\|^2_{L^2(\mathcal{I})}  \\ \displaystyle \hspace{5cm}- \int_{\mathcal{I}} B(y,s)\frac{\partial \varphi}{\partial y} \varphi\, dy     \\ \displaystyle \hspace{4.6cm}\leq \frac{b}{4}\Big\|\frac{ \partial \varphi}{\partial y}\Big\|_{L^2(\mathcal{I})}^{2} + \frac{m_T^2}{b}  \|\varphi\|_{L^2(\mathcal{I})}^2+   \varepsilon\|\varphi\|_{L^\infty(\mathcal{I})}^2 \\ \displaystyle \hspace{5cm} + \frac{c^2\gamma_1^4}{4\varepsilon} \|w\|_{L^1(\mathcal{I})}^{2} - d\gamma^{2}(t(s)) \|\varphi\|^2_{L^2(\mathcal{I})},
    \end{array}
\end{align*}
where we used $|B|\leq m_T$, $\gamma\leq\gamma_1$ and Young's inequality. Then, by \eqref{CH1} with $\varepsilon=b/(4c_1^2)$, and since $d\gamma^2(t(s))\geq d\gamma_0^2$, we get
\begin{align*}
\begin{array}{l}
  \displaystyle  \frac{1}{2}\frac{d}{ds}\|\varphi\|_{L^2(\mathcal{I})}^{2} + d\gamma^{2}_0 \|\varphi\|^2_{L^2(\mathcal{I})} + \frac{3b}{4}\Big\|\frac{ \partial \varphi}{\partial y}\Big\|_{L^2(\mathcal{I})}^{2} \leq  \frac{m_T^2}{b}  \|\varphi\|_{L^2(\mathcal{I})}^2+   \frac{b}{4}\|\varphi\|_{H^1(\mathcal{I})}^2 + \frac{C_b}{2} \|w\|_{L^1(\mathcal{I})}^{2},
    \end{array}
\end{align*}
where $C_{b}:=2c^2\gamma_1^4c_1^2/b$ and $c_1>0$ is the embedding constant of $H^1(\mathcal{I}) \subset L^{\infty}(\mathcal{I})$. Since $\|\varphi\|_{H^1(\mathcal I)}^2=\|\varphi\|_{L^2(\mathcal I)}^2+\|\partial_y\varphi\|_{L^2(\mathcal I)}^2$, we have
\begin{align}\label{fw-00}
\begin{array}{l}
  \displaystyle  \frac{d}{ds}\|\varphi\|_{L^2(\mathcal{I})}^{2} +2d\gamma_0^2\|\varphi\|_{L^2(\mathcal{I})}^{2} + b  \Big\|\frac{ \partial \varphi}{\partial y}\Big\|_{L^2(\mathcal{I})}^{2} \leq  C_1  \|\varphi\|_{L^2(\mathcal{I})}^2 + C_b \|w\|_{L^1(\mathcal{I})}^{2},
    \end{array}
\end{align}
where $C_1:=\frac{2m_T^2}{b}+\frac b2$. From \eqref{w-L1-00} and Gronwall's inequality, we obtain
\begin{align}\label{fw-01}
\begin{array}{l}
  \displaystyle  \|\varphi(s)\|_{L^2(\mathcal{I})}^{2}  \leq   F_0[T,\tau, \|\varphi_0\|_{L^2(\mathcal{I})}, \|w_0\|_{L^1(\mathcal{I})}] , \hspace{0,3cm} \text{for all} \ s\in [\tau,T]\cap[\tau,T_{w_0,\varphi_0}),
    \end{array}
\end{align}
where
\begin{equation*}
    F_0[T,\tau, \|\varphi_0\|_{L^2(\mathcal{I})}, \|w_0\|_{L^1(\mathcal{I})}]= e^{C_1(T-\tau)}\Big( \|\varphi_0\|_{L^2(\mathcal{I})}^2 + \frac{C_b C_{\gamma}^2\|w_0\|_{L^1(\mathcal{I})}^{2}}{C_1} \Big).
\end{equation*}
Moreover, integrating \eqref{fw-00} on $[\tau,s]$ and using \eqref{fw-01}, we obtain
\begin{align*}
\begin{array}{l}
  \displaystyle   \int_{\tau}^{s}\|\varphi(r)\|_{H^1(\mathcal{I})}^{2}dr \leq F_1[T,\tau, \|\varphi_0\|_{L^2(\mathcal{I})}, \|w_0\|_{L^1(\mathcal{I})}], \hspace{0,3cm} \text{for all} \ s\in [\tau,T]\cap[\tau,T_{w_0,\varphi_0}),
    \end{array}
\end{align*}
where
\begin{align*}
    F_1[T,\tau, \|\varphi_0\|_{L^2(\mathcal{I})}, \|w_0\|_{L^1(\mathcal{I})}] =& C_{b,\gamma}^{-1}\Big(\|\varphi_0\|_{L^2(\mathcal{I})}^2 +  C_1(T-\tau) F_0[T,\tau, \|\varphi_0\|_{L^2(\mathcal{I})}, \|w_0\|_{L^1(\mathcal{I})}] \\&+ C_bC_\gamma^2 (T-\tau)\|w_0\|_{L^1(\mathcal{I})}^{2}\Big),
\end{align*}
and $C_{b,\gamma}:=\min\{2d\gamma_0^2,b\}$.

On the other hand, consider the operator
\[
\widetilde A_{b,1}:H^1(\mathcal I)\to H^1(\mathcal I)'
\]
which is defined by
\[
\left\langle \widetilde A_{b,1}u,v\right\rangle
=
b\int_{\mathcal I}u_yv_y\,dy
+
\int_{\mathcal I}uv\,dy,
\qquad u,v\in H^1(\mathcal I),
\]
and corresponds to the operator $-b\Delta+I$ endowed with homogeneous Neumann boundary conditions. Moreover, $-\widetilde A_{b,1}$ generates an analytic semigroup on $H^1(\mathcal I)'$ satisfying
\begin{equation}\label{Ab1}
\Big\|e^{-\sigma\widetilde A_{b,1}}\Big\|_{\mathscr{L}(H^1(\mathcal I)')}
\leq C e^{-\sigma},\qquad
\Big\|\widetilde A_{b,1}^{15/16}e^{-\sigma\widetilde A_{b,1}}\Big\|_{\mathscr{L}(H^1(\mathcal I)')}
\leq C\sigma^{-\frac{15}{16}},\qquad 0<\sigma\leq T-\tau,
\end{equation}
where $C$ may depend on $T-\tau$. By the equivalence of weak and mild solutions for the analytic semigroup $\{e^{-s\widetilde A_{b,1}}\}_{s\geq0}$ (cf.\ \cite{Pazy1983}), we have
\begin{align*}
    \varphi(s)= e^{-\widetilde{A}_{b,1}(s-\tau)}\varphi_0 + \int_\tau^s e^{-(s-r)\widetilde{A}_{b,1}} \big[ c\gamma^2(t(r))w(r) -(d\gamma^2(t(r))-1)\varphi(r) -B(y,r)\, \partial_y \varphi(r)\big] dr.
\end{align*}
Now, using the triangle inequality and the norm equivalence $\|\widetilde A_{b,1}\psi\|_{H^1(\mathcal I)'}\sim\|\psi\|_{H^1(\mathcal I)}$,
\begin{align*}
\|\varphi(s)\|_{H^1(\mathcal I)}
&\leq C\big\|\widetilde A_{b,1}e^{-\widetilde A_{b,1}(s-\tau)}\varphi_0\big\|_{H^1(\mathcal I)'}\\
&\hspace{0.5cm}+C\int_\tau^s\Big\|\widetilde A_{b,1}e^{-(s-r)\widetilde A_{b,1}}
\Big[c\gamma^2(t(r))w-(d\gamma^2(t(r))-1)\varphi-B(y,r)\,\partial_y\varphi\Big]\Big\|_{H^1(\mathcal I)'}\,dr .
\end{align*}

\medskip
\noindent By the semigroup bound \eqref{Ab1},
the first term is bounded by $C\|\varphi_0\|_{H^1(\mathcal I)}$.
For the integral, split the bracket into its three terms and factor
$\widetilde A_{b,1}=\widetilde A_{b,1}^{15/16}\widetilde A_{b,1}^{1/16}$:
\begin{align*}
\|\varphi(s)\|_{H^1(\mathcal I)}
&\leq C\|\varphi_0\|_{H^1(\mathcal I)} +C\int_\tau^s\Big\|\widetilde A_{b,1}^{15/16}e^{-(s-r)\widetilde A_{b,1}}
\widetilde A_{b,1}^{1/16}\big[c\gamma^2(t(r))w\big]\Big\|_{H^1(\mathcal I)'}\,dr\\
&\hspace{0.5cm}+C\int_\tau^s\Big\|\widetilde A_{b,1}^{15/16}e^{-(s-r)\widetilde A_{b,1}}
\widetilde A_{b,1}^{1/16}\big[-(d\gamma^2(t(r))-1)\varphi\big]\Big\|_{H^1(\mathcal I)'}\,dr\\
&\hspace{0.5cm}+C\int_\tau^s\Big\|\widetilde A_{b,1}^{15/16}e^{-(s-r)\widetilde A_{b,1}}
\widetilde A_{b,1}^{1/16}\big[-B(y,r)\,\partial_y\varphi\big]\Big\|_{H^1(\mathcal I)'}\,dr .
\end{align*}

\medskip
\noindent Applying the smoothing estimate in \eqref{Ab1}, and the embeddings
$\|\widetilde A_{b,1}^{1/16}w\|_{H^1(\mathcal I)'}\leq C\|w\|_{L^1(\mathcal I)}$ (for $w$; recall that $\widetilde A_{b,1}^{1/16}$ maps $(H^{7/8}(\mathcal I))'$ into $H^1(\mathcal I)'$ and $L^1(\mathcal I)\hookrightarrow (H^{7/8}(\mathcal I))'$, since $H^{7/8}(\mathcal I)\hookrightarrow L^\infty(\mathcal I)$)
 and  $\|\widetilde A_{b,1}^{1/16}\psi\|_{H^1(\mathcal I)'}\leq C\|\psi\|_{L^2(\mathcal I)}$ (for $\varphi$ and
$B\,\partial_y\varphi$), together with $\gamma\leq\gamma_1$ and $|B|\leq m_T$:
\begin{align*}
\|\varphi(s)\|_{H^1(\mathcal I)}
&\leq C\|\varphi_0\|_{H^1(\mathcal I)}
+Cc\gamma_1^2\int_\tau^s(s-r)^{-\frac{15}{16}}\|w(r)\|_{L^1(\mathcal I)}\,dr\\
&\hspace{0.5cm}+C(d\gamma_1^2+1)\int_\tau^s(s-r)^{-\frac{15}{16}}\|\varphi(r)\|_{L^2(\mathcal I)}\,dr\\
&\hspace{0.5cm}+C\,m_T\int_\tau^s(s-r)^{-\frac{15}{16}}\|\partial_y\varphi(r)\|_{L^2(\mathcal I)}\,dr .
\end{align*}
Substituting the a priori bounds \eqref{w-L1-00} and \eqref{fw-01}, and using
$\int_\tau^s(s-r)^{-15/16}\,dr\leq C_\ast:=16(T-\tau)^{\frac1{16}}$ for the two
resulting terms, together with $\|\partial_y\varphi\|_{L^2(\mathcal I)}\leq\|\varphi\|_{H^1(\mathcal I)}$:
\begin{align*}
\begin{array}{l}
\displaystyle \|\varphi(s)\|_{H^1(\mathcal I)}
\leq C\|\varphi_0\|_{H^1(\mathcal I)} +Cc\gamma_1^2C_{\gamma}C_\ast\|w_0\|_{L^1(\mathcal I)}\\ \displaystyle \hspace{3cm}
+C(d\gamma_1^2+1)C_\ast\,F_0^{1/2}\big[T,\tau,\|\varphi_0\|_{L^2(\mathcal I)},\|w_0\|_{L^1(\mathcal I)}\big]\\ \displaystyle \hspace{3cm}
+C\,m_T\int_\tau^s(s-r)^{-\frac{15}{16}}\|\varphi(r)\|_{H^1(\mathcal I)}\,dr .
\end{array}
\end{align*}

To estimate $\|\varphi(s)\|_{H^1(\mathcal I)}$, we introduce
\[
Y(s):=\|\varphi(s)\|_{H^1(\mathcal I)},\qquad s\geq \tau.
\]
Moreover, let
\begin{align*}
A_{T,\tau}:={}&
C\|\varphi_0\|_{H^1(\mathcal I)}
+Cc\gamma_1^{2}C_{\gamma}C_\ast
\|w_0\|_{L^1(\mathcal I)}+
C(d\gamma_1^{2}+1)C_\ast
F_{0}^{1/2}\big[
T,\tau,\|\varphi_0\|_{L^2(\mathcal I)},
\|w_0\|_{L^1(\mathcal I)}
\big],
\end{align*}
and
\[
K_T:=C\,m_T .
\]
Then the above estimate can be written as
\[
Y(s)
\leq
A_{T,\tau}
+
K_T\int_\tau^s
(s-r)^{-\frac{15}{16}}Y(r)\,dr.
\]
Notice that $(s-r)^{-\frac{15}{16}}
=(s-r)^{\frac{1}{16}-1}$. Thus, the kernel $k(\sigma):=\sigma^{\frac{1}{16}-1},$ $\sigma>0$,
is locally integrable, since
\[
\int_0^{T-\tau} k(\sigma)\,d\sigma =
\int_0^{T-\tau} \sigma^{-\frac{15}{16}}\,d\sigma
=
16(T-\tau)^{\frac{1}{16}}
<\infty.
\]
Therefore, the above inequality is a Volterra integral inequality with a weakly singular kernel. By the singular Gronwall inequality (or Henry-Gronwall inequality, cf. \cite[Lemma 7.1.1]{Henry1981}), we obtain
\[
Y(s)
\leq
A_{T,\tau}
E_{\frac{1}{16}}
\left(
K_T\Gamma\left(\frac{1}{16}\right)
(s-\tau)^{\frac{1}{16}}
\right),
\qquad s\in[\tau,T]\cap[\tau,T_{w_0,\varphi_0}),
\]
where $E_\alpha$ denotes the Mittag--Leffler function
\[
E_\alpha(z)
:=
\sum_{k=0}^{\infty}
\frac{z^k}{\Gamma(\alpha k+1)},
\qquad \alpha>0.
\]
Since $E_{1/16}$ is increasing on $[0,\infty)$ and $s-\tau\leq T-\tau$, we consequently have
\begin{align}\label{fw-03}
\|\varphi(s)\|_{H^1(\mathcal I)}
\leq F_{2}[T,\tau,\|\varphi_0\|_{H^1(\mathcal{I})},\|w_0\|_{L^1(\mathcal{I})}],
\qquad s\in[\tau,T]\cap[\tau,T_{w_0,\varphi_0}),
\end{align}
where
 \begin{align*}
F_{2}[T,\tau,\|\varphi_0\|_{H^1(\mathcal{I})},\|w_0\|_{L^1(\mathcal{I})}]:=
A_{T,\tau}
E_{\frac{1}{16}}
\left(
K_T\Gamma\left(\frac{1}{16}\right)
(T-\tau)^{\frac{1}{16}}
\right).
 \end{align*}

Now, multiplying the first equation of \eqref{equ4iii} by $w$ and integrating on $\mathcal{I}$, we obtain
\begin{align*}
\begin{array}{l}
  \displaystyle
    \frac{1}{2}\frac{d}{ds}\|w\|_{L^2(\mathcal{I})}^{2} + a\Big\|\frac{ \partial w}{\partial y}\Big\|_{L^2(\mathcal{I})}^{2}  =  \kappa\int_{\mathcal{I}} w\frac{\partial w}{\partial y}\frac{\partial \varphi}{\partial y}dy - \int_{\mathcal{I}} B(y,s)\frac{\partial w}{\partial y} w\, dy\\ \hspace{5cm}\displaystyle \leq \kappa\int_{\mathcal{I}} w\frac{\partial w}{\partial y}\frac{\partial \varphi}{\partial y}dy + C_2\|w\|_{L^2(\mathcal{I})}\Big\|\frac{\partial w}{\partial y}\Big\|_{L^2(\mathcal{I})},
    \end{array}
\end{align*}
where $C_2:= m_T$. Then, by Young's inequality, we get
\begin{align}\label{w-000}
\begin{array}{l}
  \displaystyle
    \frac{d}{ds}\|w\|_{L^2(\mathcal{I})}^{2} + a\Big\|\frac{ \partial w}{\partial y}\Big\|_{L^2(\mathcal{I})}^{2} \leq 2\kappa\int_{\mathcal{I}} w\frac{\partial w}{\partial y}\frac{\partial \varphi}{\partial y}dy + \frac{C_2^2}{a}\|w\|_{L^2(\mathcal{I})}^2.
    \end{array}
\end{align}

By \eqref{LpLq} (namely, $\|w\|_{L^4(\mathcal I)}^2\leq C\|w\|_{H^1(\mathcal I)}\|w\|_{L^1(\mathcal I)}$ and $\|\partial_y\varphi\|_{L^4(\mathcal I)}^2\leq C\|\varphi\|_{H^2(\mathcal I)}^{1/2}\|\partial_y\varphi\|_{L^2(\mathcal I)}^{3/2}$) we have, for any $\varepsilon>0$,
\begin{align}\label{w-001}
\begin{array}{l}
  \displaystyle
   2\kappa\int_{\mathcal{I}} w  \frac{\partial \varphi}{\partial y} \frac{\partial w} {\partial y} dy \leq \frac{2\kappa^2}{a}\|w\|_{L^4(\mathcal{I})}^2\Big\| \frac{\partial \varphi}{\partial y} \Big\|^2_{L^4(\mathcal{I})} + \frac{a}{2}\Big\|\frac{\partial w}{\partial y}\Big\|^2_{L^2(\mathcal{I})}\\ \hspace{3cm}\displaystyle \leq \frac{C\kappa^2}{a} \|w\|_{H^1(\mathcal{I})} \|w\|_{L^1(\mathcal{I})}\|\varphi\|^{1/2}_{H^2(\mathcal{I})}\Big\|\frac{\partial \varphi}{\partial y} \Big\|^{3/2}_{L^2(\mathcal{I})} +  \frac{a}{2}\Big\|\frac{\partial w}{\partial y}\Big\|^2_{L^2(\mathcal{I})}\\  \hspace{3cm}\displaystyle\leq \varepsilon( \|w\|_{H^1(\mathcal{I})}^2 + \|\varphi\|^{2}_{H^2(\mathcal{I})} ) + C_{\varepsilon,\kappa} \|w\|_{L^1(\mathcal{I})}^4 \|\varphi\|_{H^1(\mathcal{I})}^{6}+ \frac{a}{2}\Big\|\frac{\partial w}{\partial y}\Big\|^2_{L^2(\mathcal{I})},
   \end{array}
\end{align}
where in the last step we used Young's inequality with exponents $2$, $4$ and $4$, and $\|\partial_y\varphi\|_{L^2(\mathcal I)}\leq\|\varphi\|_{H^1(\mathcal I)}$.

Now, multiplying the second equation of \eqref{equ4iii} by $- \frac{\partial^2\varphi}{\partial y^2}$ and integrating on $\mathcal{I}$, we obtain
\begin{align*}
\begin{array}{l}
  \displaystyle
    \frac{1}{2}\frac{d}{ds}\Big\|\frac{\partial \varphi}{\partial y}\Big \|_{L^2(\mathcal{I})}^{2} + b\Big\|\frac{ \partial^2 \varphi}{\partial y^2}\Big\|_{L^2(\mathcal{I})}^{2} = -c\gamma^2(t(s))\int_{\mathcal{I}} w\frac{\partial^2 \varphi}{\partial y^2} dy - d\gamma^{2}(t(s)) \Big\|\frac{\partial \varphi}{\partial y}\Big \|_{L^2(\mathcal{I})}^{2}\\ \hspace{6cm} \displaystyle  + \int_{\mathcal{I}} B(y,s)\frac{\partial \varphi}{\partial y} \frac{\partial^2 \varphi}{\partial y^2}dy     \\ \hspace{4.5cm} \displaystyle \leq  \frac{b}{2}\Big\|\frac{ \partial^2 \varphi}{\partial y^2 }\Big\|_{L^2(\mathcal{I})}^{2} + \frac{C_3}{2}  \Big\|\frac{\partial \varphi}{\partial y}\Big \|_{L^2(\mathcal{I})}^{2} + \frac{C_4}{2} \|w \|_{L^2(\mathcal{I})}^{2} - d\gamma^{2}(t(s)) \Big\|\frac{\partial \varphi}{\partial y}\Big \|_{L^2(\mathcal{I})}^{2},
    \end{array}
\end{align*}
where $C_3:=\frac{2 C_2^2}{b} $, $C_4:=\frac{2c^2 \gamma_1^4}{b}$. Hence, since $d\gamma^2(t(s))\geq d\gamma_0^2$,
\begin{align}\label{w-002}
\begin{array}{l}
  \displaystyle
    \frac{d}{ds}\Big\|\frac{\partial \varphi}{\partial y}\Big \|_{L^2(\mathcal{I})}^{2} + 2d\gamma_0^2\Big\|\frac{\partial \varphi}{\partial y}\Big \|_{L^2(\mathcal{I})}^{2}  + b\Big\|\frac{ \partial^2 \varphi}{\partial y^2}\Big\|_{L^2(\mathcal{I})}^{2}  \leq C_3  \Big\|\frac{\partial \varphi}{\partial y}\Big \|_{L^2(\mathcal{I})}^{2} + C_4\|w \|_{L^2(\mathcal{I})}^{2}.
    \end{array}
\end{align}

Thus, adding \eqref{w-000} and \eqref{w-002}, using \eqref{w-001} with $\varepsilon:=\frac{\min\{a,b\}}{4}$, and
$\|w\|_{H^1(\mathcal I)}^2=\|w\|_{L^2(\mathcal I)}^2+\|\partial_yw\|_{L^2(\mathcal I)}^2$,
$\|\varphi\|_{H^2(\mathcal I)}^2=\|\varphi\|_{L^2(\mathcal I)}^2+\|\partial_y\varphi\|_{L^2(\mathcal I)}^2+\|\partial_y^2\varphi\|_{L^2(\mathcal I)}^2$, we have
\begin{align}\label{dphi-0-n}
\begin{array}{l}
  \displaystyle \frac{d}{ds}\Big[ \|w\|_{L^2(\mathcal{I})}^{2} + \Big\|\frac{\partial \varphi}{\partial y}\Big \|_{L^2(\mathcal{I})}^{2}\Big] + \frac a4\Big\|\frac{\partial w}{\partial y}\Big\|^2_{L^2(\mathcal{I})}  + \frac{3b}{4}\Big\|\frac{ \partial^2 \varphi}{\partial y^2}\Big\|_{L^2(\mathcal{I})}^{2}\\ \displaystyle\hspace{2cm}\leq \varepsilon\|\varphi\|_{L^2(\mathcal{I})}^2 + C_5 \Big[ \|w\|_{L^2(\mathcal{I})}^{2} + \Big\|\frac{\partial \varphi}{\partial y}\Big \|_{L^2(\mathcal{I})}^{2}\Big]  +C_{\varepsilon,\kappa} \|w\|_{L^1(\mathcal{I})}^4 \|\varphi\|_{H^1(\mathcal{I})}^{6},
  \end{array}
\end{align}
where $C_5:=\frac{C_2^2}{a}+C_3+C_4+\varepsilon$ (we discarded the nonnegative term $2d\gamma_0^2\|\partial_y\varphi\|^2_{L^2(\mathcal I)}$). Now, define
\[
Y(s):= \|w(s)\|_{L^2(\mathcal I)}^2 +\left\|\frac{\partial\varphi}{\partial y}(s)\right\|_{L^2(\mathcal I)}^2.
\]
Then, from \eqref{dphi-0-n}, discarding the nonnegative terms on the left-hand side, we obtain
\[
Y'(s) \leq C_5Y(s) + \varepsilon\|\varphi(s)\|_{L^2(\mathcal I)}^2 + C_{\varepsilon,\kappa}\|w(s)\|_{L^1(\mathcal I)}^4
\|\varphi(s)\|_{H^1(\mathcal I)}^6.
\]
By \eqref{w-L1-00}, \eqref{fw-01} and \eqref{fw-03},
\[
\|\varphi(s)\|_{L^2(\mathcal I)}^2
\leq F_0,
\qquad
\|\varphi(s)\|_{H^1(\mathcal I)}
\leq F_2,
\qquad
\|w(s)\|_{L^1(\mathcal I)}
\leq
C_{\gamma}\|w_0\|_{L^1(\mathcal I)},
\]
for all $s\in[\tau,T]\cap[\tau,T_{w_0,\varphi_0})$. Consequently,
\[
Y'(s)
\leq
C_5Y(s)+G_{T,\tau},
\qquad
G_{T,\tau}
:=
\varepsilon F_0
+
C_{\varepsilon,\kappa}C_{\gamma}^4
\|w_0\|_{L^1(\mathcal I)}^4F_2^6.
\]
Therefore, by Gronwall's inequality,
\[
Y(s)
\leq
e^{C_5(s-\tau)}Y(\tau)
+
G_{T,\tau}
\int_\tau^s e^{C_5(s-r)}\,dr
\leq
e^{C_5(T-\tau)}Y(\tau)
+
\frac{G_{T,\tau}}{C_5} e^{C_5(T-\tau)}.
\]
Since $Y(\tau)=\|w_0\|_{L^2(\mathcal I)}^2 + \left\|\frac{\partial\varphi_0}{\partial y}\right\|_{L^2(\mathcal I)}^2
\leq \|w_0\|_{L^2(\mathcal I)}^2 +\|\varphi_0\|_{H^1(\mathcal I)}^2$ and $\|w(s)\|_{L^2(\mathcal I)}^2\leq Y(s)$, we obtain
\[
\|w(s)\|_{L^2(\mathcal I)} \leq F_3\big[
T,\tau,\|\varphi_0\|_{H^1(\mathcal I)}, \|w_0\|_{L^2(\mathcal I)} \big],
\qquad s\in[\tau,T]\cap[\tau,T_{w_0,\varphi_0}),
\]
where
\[
F_3 ^2\big[ T,\tau,\|\varphi_0\|_{H^1(\mathcal I)},
\|w_0\|_{L^2(\mathcal I)} \big] :=  e^{C_5(T-\tau)}
\left[ \|w_0\|_{L^2(\mathcal I)}^2
+ \|\varphi_0\|_{H^1(\mathcal I)}^2
+ \frac{G_{T,\tau}}{C_5}\right].
\]
Since the above a priori estimates \eqref{fw-03} and the bound for $\|w(s)\|_{L^2(\mathcal I)}$ are uniform on every finite time interval, the solution cannot blow up in finite time (by Theorem~\ref{T-Exis}, the local existence time depends only on an upper bound for the norm of the initial datum). Consequently, the local solution can be extended to any arbitrary time $s>\tau$, and hence it exists globally in time.
\end{proof}

We are now ready to establish the main result of this section, namely the existence of global weak solutions to problem \eqref{equ4}.

\begin{theorem}\label{teoEx-05-It} {\it (Existence of global weak solution to \eqref{equ4})} Under the assumptions given in (H1)-(H2), given $(u_0,\rho_0) \in L^2(\mathcal{I}_\tau)\times H^{1}(\mathcal{I}_\tau)$, with $u_0\geq 0$, $\rho_0\geq 0$ and $T, \tau \in \mathbb{R}$ with $T>\tau$, then the  problem \eqref{equ4} has a unique global weak solution $(u(\cdot),\rho(\cdot))=(u(\cdot; \tau ,u_0,\varphi_0),\rho(\cdot; \tau ,u_0,\rho_0))$ with 
\begin{align*}
   &u \in H^1 (\tau , T ; H^{1}(\mathcal{I}_t)')\cap \mathcal{C}([\tau, T]; L^2(\mathcal{I}_t)) \cap L^2 (\tau,T ; H^{1}(\mathcal{I}_t)) \\&
   \rho \in H^1 (\tau , T ; L^{2}(\mathcal{I}_t))\cap \mathcal{C}([\tau, T]; H^1(\mathcal{I}_t)) \cap L^2 (\tau,T ; H^{2}_N(\mathcal{I}_t)),
\end{align*}
 satisfying the initial condition $(u(\tau ),\rho( \tau ))=(u_0,\rho_0)$.
\end{theorem}

\section{Existence of Pullback Attractors}\label{sec4}

In this section we are interested in studying the asymptotic behavior, in the pullback sense, of the weak solutions of the problem \eqref{equ4}, on a subfamily of the family of Banach spaces $\mathcal{H}_t=L^2(\mathcal{I}_t)\times H^1(\mathcal{I}_t)$. For this, we will consider the hypotheses (H1)-(H4). Then, Theorem \ref{teoEx-05-It} guarantees the existence and uniqueness of weak global solution of the problem \eqref{equ4} whose initial condition belongs to $\mathcal{H}_\tau$ with $\tau \in \mathbb{R}$. 

Now, integrating the first equation of \eqref{equ4}, and taking into account the Neumman condition, we have
\begin{equation*}
\int_{\alpha(t)}^{\beta(t)}\frac{\partial u}{\partial t}(x,t)dx =0.
\end{equation*}
Then, we obtain
\begin{equation}\label{eq:balance}
\int_{\alpha(t)}^{\beta(t)}\frac{\partial u}{\partial t}(x,t)dx= \frac{d}{dt}\int_{\alpha(t)}^{\beta(t)}u(x,t)dx - \beta'(t)u(\beta(t),t) + \alpha'(t)u(\alpha(t),t)=0.
\end{equation}
From (H2) we have $\beta'(t)u(\beta(t),t) - \alpha'(t)u(\alpha(t),t)\leq 0$ then, we obtain
\begin{equation}\label{pos-0}
   \int_{\alpha(t)}^{\beta(t)}u(x,t)dx \leq \int_{\alpha(\tau)}^{\beta(\tau)}u_0(x)dx.
\end{equation}

Given $M>0$ and $t\in \mathbb{R}$, denote by
\begin{equation*}
    \mathbb{H}^{M}_{1,t} :=\bigg\{ u\in L^2(\mathcal{I}_t): u\geq 0  \ \text{ and } \ \int_{\mathcal{I}_t} u(x)dx \leq M \bigg\}\times \Big\{\rho \in H^1(\mathcal{I}_t): \rho\geq 0\Big\},
\end{equation*}
with the topology induced from $\mathcal{H}_t = L^2 (\mathcal{I}_t) \times H^1(\mathcal{I}_t)$. 

Now, from Theorem \ref{teoEx-05-It} and taking into account \eqref{pos-0}, let us define the bi-parametric families of mappings $\{U_1(t,\tau): \mathbb{H}^M_{1,\tau} \to \mathbb{H}^M_{1,t}:(t,\tau)\in \mathbb{R}_{d}^{2}\}$ 
\begin{equation}\label{PU} \displaystyle 
 U_1(t,\tau)(u_0,\rho_0)=\big(u(t;\tau,u_0,\rho_0),\rho(\cdot;\tau,u_0,\rho_0)\big),
\end{equation}
where $(u,\rho)$ is the unique weak solution of \eqref{equ4} associated to the initial condition $(u_0,\rho_0)\in \mathbb{H}^M_{1,\tau}$ for any $(t,\tau)\in \mathbb{R}_{d}^{2}$.

\begin{theorem}
  Under the hypotheses (H1)-(H2), the bi-parametric family $U_{1}$, defined in \eqref{PU} is a continuous process.
\end{theorem}
\begin{proof} It follows the same line regarding uniqueness (continuous dependence on initial data) guaranteed by Theorem \ref{teoEx-05-It}.
\end{proof}

\begin{lemma}\label{lemmaR11}
Under the hypotheses (H1)-(H4) the evolution process $U_1$ is dissipative.
\end{lemma}
\begin{proof} Consider $t\in \mathbb{R}$ and $\widehat{D}\in \mathcal{D}_{F}(\mathbb{H}_1^M)$. Now, by \eqref{pos-0}, we have
\begin{equation}\label{M0}
   \int_{\mathcal{I}_t}u(x,t)dx  \leq \int_{\mathcal{I}_\tau}u_\tau(x)dx \leq M.
\end{equation}
From the second equation of \eqref{equ4}, using Proposition~\ref{Prop7}~(iii) and the homogeneous Neumann boundary condition, we deduce
\begin{align*}
    \frac{1}{2}\frac{d}{ dt}\|\rho\|_{L^2(\mathcal{I}_t)}^2 +b\Big\|\frac{\partial \rho}{\partial x}\Big\|^2_{L^2(\mathcal{I}_t)}=  c\int_{\mathcal{I}_t}u\rho\,dx -d\int_{\mathcal{I}_t}\rho^2\,dx
    \;+\;\frac12\Big(\beta'(t)\,|\rho(\beta(t),t)|^2-\alpha'(t)\,|\rho(\alpha(t),t)|^2\Big) .
\end{align*}
By (H2), the last term is nonpositive. Thus, by the H\"older inequality, \eqref{M0} and $H^1(\mathcal{I}_t) \hookrightarrow L^\infty(\mathcal{I}_t)$ (with embedding constant $c_\infty>0$ independent of $t$, by Lemma~\ref{lem:scaling} and (H1)), we have
\begin{align*}
    \frac{1}{2}\frac{d}{ dt}\|\rho\|_{L^2(\mathcal{I}_t)}^2 + d\|\rho\|_{L^2(\mathcal{I}_t)}^2+b\Big\|\frac{\partial \rho}{\partial x}\Big\|^2_{L^2(\mathcal{I}_t)}&\leq c\|\rho\|_{L^\infty(\mathcal{I}_t)} \|u\|_{L^1(\mathcal{I}_t)}\\& \leq \frac{\min\{b,d\}}{2}\|\rho\|_{H^1(\mathcal{I}_t)}^2 + \frac{C_{b}}{2}M^2,
\end{align*}
with $C_b:=c^2c_\infty^2/\min\{b,d\}$. Since $\|\rho\|_{H^1(\mathcal I_t)}^2=\|\rho\|_{L^2(\mathcal I_t)}^2+\|\partial_x\rho\|_{L^2(\mathcal I_t)}^2$ and $\min\{b,d\}\leq b,d$, we get
\begin{align*}
    \frac{d}{ dt}\|\rho\|_{L^2(\mathcal{I}_t)}^2 + d\,\|\rho\|_{L^2(\mathcal{I}_t)}^2+b\,\Big\|\frac{\partial \rho}{\partial x}\Big\|^2_{L^2(\mathcal{I}_t)}& \leq  C_{b}M^2 .
\end{align*}
Setting $\nu_d:=d$, by the Gr\"onwall inequality
\begin{align}\label{ab-001}
\begin{array}{l}
  \displaystyle  \|\rho(t)\|_{L^2(\mathcal{I}_t)}^{2}  \leq   \|\rho_0\|_{L^2(\mathcal{I}_\tau)}^2e^{-\nu_{d}(t-\tau)} + \frac{C_b}{\nu_{d}} M^{2}.
    \end{array}
\end{align}

\medskip
\noindent\textbf{Estimate of $\|\rho(t)\|_{H^1(\mathcal I_t)}$.} This estimate is carried out in the rescaled time $s$, where the semigroup acts, and only at the end we return to the original time $t$. Consider
\[
\rho(x,t)=\varphi\big(r(x,t),s(t)\big)\qquad\text{and}\qquad u(x,t)=w\big(r(x,t),s(t)\big),
\]
where $(w,\varphi)$ is the solution of \eqref{equ4iii} with initial datum $(w_0,\varphi_0)=(u_0\circ r^{-1}(\cdot,\tau),\rho_0\circ r^{-1}(\cdot,\tau))$ (recall that $s(\tau)=\tau$). Denote by $c_1,C_1$ the constants in \eqref{eq:equiv-Hk} with $k=1$. By Lemma~\ref{lem:scaling}, for every $\sigma\geq\tau$,
\begin{align}
\|w(\sigma)\|_{L^1(\mathcal I)}&=\frac{\|u(t(\sigma))\|_{L^1(\mathcal I_{t(\sigma)})}}{\gamma(t(\sigma))}\leq\frac{M}{\gamma_0},\label{eq:w-L1-s}\\
\|\varphi(\sigma)\|_{L^2(\mathcal I)}&=\frac{\|\rho(t(\sigma))\|_{L^2(\mathcal I_{t(\sigma)})}}{\gamma(t(\sigma))^{1/2}}\leq\frac{\|\rho(t(\sigma))\|_{L^2(\mathcal I_{t(\sigma)})}}{\gamma_0^{1/2}},
\qquad
\|\rho_0\|_{L^2(\mathcal I_\tau)}\leq\gamma_1^{1/2}\|\varphi_0\|_{L^2(\mathcal I)}.\label{eq:phi-L2-s}
\end{align}

\smallskip
\noindent\emph{Semigroup estimates, see \cite{Pazy1983}.} The operator $\widetilde A_{b,1}$ is self-adjoint and positive, and its spectrum is contained in $[1,\infty)$, the lowest eigenvalue $\lambda_0=1$ corresponding to the constant functions. Hence, by spectral calculus, for $\theta\in(0,1]$ and $\sigma>0$,
\[
\sup_{\lambda\geq1}\lambda^\theta e^{-\sigma\lambda}
\leq e^{-\sigma/2}\sup_{\lambda\geq1}\lambda^\theta e^{-\sigma\lambda/2}
\leq\Big(\frac{2\theta}{e\,\sigma}\Big)^{\theta}e^{-\sigma/2}.
\]
Note that a bound of the form $C\sigma^{-\theta}e^{-\sigma}$ cannot hold for all $\sigma>0$, since $\lambda_0=1$ gives $\|\widetilde A_{b,1}^\theta e^{-\sigma\widetilde A_{b,1}}\|\geq e^{-\sigma}$; this is why the exponent $\sigma/2$ appears. Since $2\theta/e<1$ for $\theta=15/16$, there is $C_S\geq1$ (which only accounts for the equivalence between the norm of $H^1(\mathcal I)'$ and the one induced by $\widetilde A_{b,1}$) such that
\begin{equation}\label{eq:smoothing-new}
\big\|e^{-\sigma\widetilde A_{b,1}}\big\|_{\mathscr L(H^1(\mathcal I)')}\leq C_Se^{-\sigma},
\qquad
\big\|\widetilde A_{b,1}^{15/16}e^{-\sigma\widetilde A_{b,1}}\big\|_{\mathscr L(H^1(\mathcal I)')}\leq C_S\,\sigma^{-\frac{15}{16}}e^{-\frac{\sigma}{2}},
\qquad\sigma>0 .
\end{equation}
Moreover, as in the proof of Theorem~\ref{teoEx-05}, there is $C_E>0$ with $\|\widetilde A_{b,1}^{1/16}w\|_{H^1(\mathcal I)'}\leq C_E\|w\|_{L^1(\mathcal I)}$ and $\|\widetilde A_{b,1}^{1/16}\psi\|_{H^1(\mathcal I)'}\leq C_E\|\psi\|_{L^2(\mathcal I)}$, and, with $c_A:=\min\{1,b\}$ and $C_A:=\max\{1,b\}$,
$c_A\|\psi\|_{H^1(\mathcal I)}\leq\|\widetilde A_{b,1}\psi\|_{H^1(\mathcal I)'}\leq C_A\|\psi\|_{H^1(\mathcal I)}$.
Finally, by (H3), for all $(y,\sigma)$,
\[
|B(y,\sigma)|\leq m_\infty:=\gamma_1\big(\|\alpha'\|_{\infty}+\|\gamma'\|_{\infty}\big),
\qquad B(y,\sigma)=\gamma^2(t(\sigma))\,b(y,t(\sigma)) .
\]

\smallskip
\noindent\emph{Variation of constants in the time $s$.} The second equation of \eqref{equ4iii} reads $\varphi'+\widetilde A_{b,1}\varphi=c\gamma^2(t(\sigma))w-(d\gamma^2(t(\sigma))-1)\varphi-B(y,\sigma)\,\partial_y\varphi$, so that
\[
\varphi(s)=e^{-(s-\tau)\widetilde A_{b,1}}\varphi_0+\int_\tau^se^{-(s-\sigma)\widetilde A_{b,1}}\Big[c\gamma^2(t(\sigma))w(\sigma)-(d\gamma^2(t(\sigma))-1)\varphi(\sigma)-B(y,\sigma)\,\partial_y\varphi(\sigma)\Big]d\sigma .
\]
Set $Y(\sigma):=\|\varphi(\sigma)\|_{H^1(\mathcal I)}$, which is continuous on $[\tau,\infty)$. Applying $\widetilde A_{b,1}=\widetilde A_{b,1}^{15/16}\widetilde A_{b,1}^{1/16}$, the estimates above and $\gamma\leq\gamma_1$, we obtain, with $K:=C_SC_E/c_A$,
\begin{align}\label{eq:main-estimate}
\begin{array}{l}
\displaystyle Y(s)\leq\frac{C_SC_A}{c_A}e^{-(s-\tau)}Y(\tau)
+K\,c\gamma_1^2\int_\tau^s(s-\sigma)^{-\frac{15}{16}}e^{-\frac{s-\sigma}{2}}\|w(\sigma)\|_{L^1(\mathcal I)}\,d\sigma\\[10pt]
\displaystyle\hspace{1.4cm}+K(d\gamma_1^2+1)\int_\tau^s(s-\sigma)^{-\frac{15}{16}}e^{-\frac{s-\sigma}{2}}\|\varphi(\sigma)\|_{L^2(\mathcal I)}\,d\sigma
+K\,m_\infty\int_\tau^s(s-\sigma)^{-\frac{15}{16}}e^{-\frac{s-\sigma}{2}}\,Y(\sigma)\,d\sigma .
\end{array}
\end{align}
With the change of variable $\sigma'=s-\sigma$,
\begin{equation*}
\int_\tau^{s}(s-\sigma)^{-\frac{15}{16}}e^{-\frac{s-\sigma}{2}}\,d\sigma
=\int_0^{s-\tau}\sigma'^{-\frac{15}{16}}e^{-\frac{\sigma'}{2}}\,d\sigma'
\;\leq\;
\int_0^{\infty}\sigma'^{-\frac{15}{16}}e^{-\frac{\sigma'}{2}}\,d\sigma'
=2^{\frac1{16}}\,\Gamma\Big(\frac{1}{16}\Big)<\infty .
\end{equation*}
Hence, by \eqref{eq:w-L1-s}, the $w$-term in \eqref{eq:main-estimate} is bounded by $2^{1/16}\Gamma(\frac1{16})\,K\,c\gamma_1^2\,M/\gamma_0$.

\smallskip
\noindent\emph{The $L^2$-term.} Since $a,b\geq 0$ implies $\sqrt{a+b}\leq \sqrt a+\sqrt b$, estimate \eqref{ab-001}, \eqref{eq:phi-L2-s} and $t(\sigma)-\tau\geq\gamma_0^2(\sigma-\tau)$ (because $t'(\sigma)=\gamma^2(t(\sigma))\geq\gamma_0^2$ and $t(\tau)=\tau$) yield
\begin{equation}\label{eq:sqrt-split}
\|\varphi(\sigma)\|_{L^2(\mathcal{I})}
\;\leq\;
\Big(\frac{\gamma_1}{\gamma_0}\Big)^{1/2}\|\varphi_0\|_{L^2(\mathcal{I})}\,
e^{-\nu(\sigma-\tau)}
+
\Big(\frac{C_b}{\nu_d\gamma_0}\Big)^{1/2}M,
\qquad \sigma\geq \tau,\qquad \nu:=\frac{\nu_d\gamma_0^2}{2}.
\end{equation}
Substituting \eqref{eq:sqrt-split} into the $L^2$-integral in \eqref{eq:main-estimate} and using linearity,
\begin{equation*}
\int_\tau^{s}(s-\sigma)^{-\frac{15}{16}}e^{-\frac{s-\sigma}{2}}\|\varphi(\sigma)\|_{L^2(\mathcal I)}\,d\sigma
\;\leq\;
I_1(s)+I_2(s),
\end{equation*}
where
\begin{align*}
I_1(s)&:=\Big(\frac{\gamma_1}{\gamma_0}\Big)^{1/2}\|\varphi_0\|_{L^2(\mathcal{I})}
\int_\tau^{s}(s-\sigma)^{-\frac{15}{16}}e^{-\frac{s-\sigma}{2}}
e^{-\nu(\sigma-\tau)}\,d\sigma,
\\
I_2(s)&:=\Big(\frac{C_b}{\nu_d\gamma_0}\Big)^{1/2}M
\int_\tau^{s}(s-\sigma)^{-\frac{15}{16}}e^{-\frac{s-\sigma}{2}}\,d\sigma .
\end{align*}
As above,
\begin{equation}\label{eq:I2-bound}
I_2(s)\;\leq\;2^{\frac1{16}}\Gamma\Big(\frac{1}{16}\Big)\Big(\frac{C_b}{\nu_d\gamma_0}\Big)^{1/2}M.
\end{equation}
Set $S:=s-\tau$ and perform the change of variable $\sigma'=s-\sigma\in[0,S]$, so that $\sigma-\tau=S-\sigma'$. Then
\[
I_1(s)=\Big(\frac{\gamma_1}{\gamma_0}\Big)^{1/2}\|\varphi_0\|_{L^2(\mathcal{I})}
\int_0^{S}\sigma'^{-\frac{15}{16}}
\exp\!\Big[-\Big(\frac{\sigma'}{2}+\nu(S-\sigma')\Big)\Big]\,d\sigma'.
\]
The exponent
\[
p(\sigma'):=\frac{\sigma'}{2}+\nu(S-\sigma')
=\Big(\frac12-\nu\Big)\sigma'+\nu S
\]
is affine in $\sigma'$, hence its minimum over $[0,S]$ is attained at one of the endpoints:
$p(0)=\nu S$ and $p(S)=S/2$, so that, with $\delta_1:=\min\{\frac12,\nu\}$,
\begin{equation}\label{eq:phi-lower}
p(\sigma')\;\geq\;\delta_1\,S,
\qquad \forall\,\sigma'\in[0,S].
\end{equation}
Combining \eqref{eq:phi-lower} with $\int_0^S\sigma'^{-15/16}\,d\sigma'=16\,S^{1/16}$ gives
\begin{equation}\label{eq:I1-bound-raw}
I_1(s)\;\leq\;16\,\Big(\frac{\gamma_1}{\gamma_0}\Big)^{1/2}\|\varphi_0\|_{L^2(\mathcal{I})}\,
S^{\frac{1}{16}}\,e^{-\delta_1 S}.
\end{equation}
Fix $\eta\in(0,\delta_1)$. The function $S\mapsto S^{1/16}e^{-\eta S}$ attains its maximum on $[0,\infty)$ at $S^{\ast}=\dfrac{1}{16\eta}$, with maximal value $(16 e\eta)^{-1/16}$; hence
\begin{equation}\label{eq:absorb}
S^{\frac{1}{16}}e^{-\delta_1 S}
=\big(S^{\frac{1}{16}}e^{-\eta S}\big)\,e^{-(\delta_1-\eta)S}
\;\leq\;
C_\eta\, e^{-(\delta_1-\eta)S},
\qquad
C_\eta=\Big(\frac{1}{16\,e\,\eta}\Big)^{1/16}.
\end{equation}
Substituting \eqref{eq:absorb} into \eqref{eq:I1-bound-raw} and recalling $S=s-\tau$,
\begin{equation}\label{eq:I1-final}
I_1(s)\;\leq\;16\,C_\eta\Big(\frac{\gamma_1}{\gamma_0}\Big)^{1/2}\|\varphi_0\|_{L^2(\mathcal{I})}\,
e^{-(\delta_1-\eta)(s-\tau)}.
\end{equation}

\smallskip
\noindent\emph{The Volterra inequality in the time $s$.} Since $\|\varphi_0\|_{L^2(\mathcal I)}\leq Y(\tau)$ and, for $\lambda:=\min\{1,\delta_1-\eta\}$, $e^{-(s-\tau)}\leq e^{-\lambda(s-\tau)}$ and $e^{-(\delta_1-\eta)(s-\tau)}\leq e^{-\lambda(s-\tau)}$, inequality \eqref{eq:main-estimate}, together with \eqref{eq:I2-bound} and \eqref{eq:I1-final}, gives
\begin{equation}\label{eq:main-rho}
Y(s)\;\leq\;A\,e^{-\lambda(s-\tau)} + B
+K\,m_\infty\int_\tau^s(s-\sigma)^{-\frac{15}{16}}e^{-\frac{s-\sigma}{2}}\,Y(\sigma)\,d\sigma,
\qquad s\geq\tau,
\end{equation}
where
\begin{align*}
A&:=A_0\,Y(\tau),\qquad A_0:=\frac{C_SC_A}{c_A}+16\,K\,C_\eta\,(d\gamma_1^2+1)\Big(\frac{\gamma_1}{\gamma_0}\Big)^{1/2},
\\
B&:=B_0M,\qquad B_0:=2^{\frac1{16}}\Gamma\Big(\frac1{16}\Big)K\Big[\frac{c\gamma_1^2}{\gamma_0}+(d\gamma_1^2+1)\Big(\frac{C_b}{\nu_d\gamma_0}\Big)^{1/2}\Big].
\end{align*}
The kernel in \eqref{eq:main-rho} is tempered by $e^{-(s-\sigma)/2}$, while Lemma~\ref{lem:main} is stated for $e^{-(s-\sigma)}$. We reduce to it by the change of time $s=\tau+2p$, $\sigma=\tau+2q$: setting $\widehat Y(p):=Y(\tau+2p)$, $p\geq0$, and using $d\sigma=2\,dq$,
\[
\int_\tau^s(s-\sigma)^{-\frac{15}{16}}e^{-\frac{s-\sigma}{2}}Y(\sigma)\,d\sigma
=2^{\frac1{16}}\int_0^p(p-q)^{-\frac{15}{16}}e^{-(p-q)}\,\widehat Y(q)\,dq,
\]
so that \eqref{eq:main-rho} becomes
\[
\widehat Y(p)\;\leq\;A\,e^{-2\lambda p}+B+2^{\frac1{16}}K\,m_\infty\int_0^p(p-q)^{-\frac{15}{16}}e^{-(p-q)}\,\widehat Y(q)\,dq,
\qquad p\geq0 .
\]

Now, by (H4) we know that
\begin{equation}\label{eq:H5}
    2^{\frac1{16}}K\,m_\infty\,\Gamma(1-\alpha)<1, \text{ where } \ \alpha=\frac{15}{16},
\end{equation}
that is, (H4) holds with $\theta=15/16$ and $C=2^{1/16}K=2^{1/16}C_SC_E/c_A$.
By Lemma \ref{lem:main} (applied to $\widehat Y$, with initial time $0$, $C=2^{1/16}K\,m_\infty$, $\alpha=15/16$ and rate $2\lambda$), hypothesis \eqref{eq:H5} guarantees the existence of $\widehat\delta\in(0,\min\{1,2\lambda\})$ and $C_0\geq1$ (both independent of $\rho_0$, $u_0$ and $M$) such that $\widehat Y(p)\leq C_0A\,e^{-\widehat\delta p}+C_0B$ for $p\geq0$, that is,
\[
\|\varphi(s)\|_{H^1(\mathcal I)}\leq C_0A_0\,\|\varphi_0\|_{H^1(\mathcal I)}\,e^{-\frac{\widehat\delta}{2}(s-\tau)}+C_0B_0M,
\qquad s\geq\tau .
\]

\smallskip
\noindent\emph{Back to the original time $t$.} By Lemma~\ref{lem:scaling}, $\|\rho(t)\|_{H^1(\mathcal I_t)}\leq C_1\|\varphi(s(t))\|_{H^1(\mathcal I)}$ and $\|\varphi_0\|_{H^1(\mathcal I)}\leq c_1^{-1}\|\rho_0\|_{H^1(\mathcal I_\tau)}$; moreover, $s(t)-\tau\geq(t-\tau)/\gamma_1^2$, because $s'(t)=\gamma^{-2}(t)\geq\gamma_1^{-2}$. Therefore, setting $\delta^\ast:=\widehat\delta/(2\gamma_1^2)$,
\begin{equation}\label{abs-01}
\|\rho(t)\|_{H^1(\mathcal I_t)}\;\leq\;
\widehat C_{0,\gamma}\,\|\rho_0\|_{H^1(\mathcal I_\tau)}\,e^{-\delta^{\ast}(t-\tau)}
+\widehat C_{1,\gamma}\,M,
\qquad\forall\,t\geq\tau,
\end{equation}
where
\begin{equation*}
\widehat C_{0,\gamma}:=\frac{C_1C_0A_0}{c_1},
\qquad
\widehat C_{1,\gamma}:=C_1C_0B_0
\end{equation*}
depend only on $a,b,c,d$, $\gamma_0$, $\gamma_1$, $m_\infty$ and $\eta$; in particular $\widehat C_{0,\gamma}$ and $\widehat C_{1,\gamma}$ are independent of $\rho_0$, $u_0$ and $M$.

Now, multiplying the first equation of \eqref{equ4} by $u$ and integrating on $\mathcal{I}_t$, by Proposition~\ref{Prop7}~(iii) and (H2) we obtain
\begin{align}\label{ab-w-002}
    \frac{1}{2}\frac{d}{dt}\|u\|_{L^2(\mathcal{I}_t)}^{2} + a\Big\|\frac{ \partial u}{\partial x}\Big\|_{L^2(\mathcal{I}_t)}^{2} &=\kappa\int_{\mathcal{I}_t} u\frac{\partial u}{\partial x}\frac{\partial \rho}{\partial x}dx+\frac12\Big(\beta'(t)|u(\beta(t),t)|^2-\alpha'(t)|u(\alpha(t),t)|^2\Big)\notag\\
    &\leq \kappa\int_{\mathcal{I}_t} u\frac{\partial u}{\partial x}\frac{\partial \rho}{\partial x}dx.
\end{align}
By \eqref{LpLq} and the same way as \eqref{w-001}, we have
\begin{align}\label{ab-w-001}
\begin{array}{l}
  \displaystyle
   \kappa\int_{\mathcal{I}_t} u  \frac{\partial \rho}{\partial x} \frac{\partial u} {\partial x} dx \leq  \varepsilon( \|u\|_{H^1(\mathcal{I}_t)}^2 + \|\rho\|^{2}_{H^2(\mathcal{I}_t)} ) + C_{\varepsilon} \|u\|_{L^1(\mathcal{I}_t)}^4 \|\rho\|_{H^1(\mathcal{I}_t)}^{6},
   \end{array}
\end{align}
for any $\varepsilon>0$. Now, multiplying the second equation of
\eqref{equ4} by $-\partial_x^2\rho$ and integrating on $\mathcal I_t$ (no boundary term appears, since $\partial_x\rho=0$ at $x\in\{\alpha(t),\beta(t)\}$),
\begin{align*}
\begin{array}{l}
  \displaystyle
    \frac{1}{2}\frac{d}{dt}\Big\|\frac{\partial \rho}{\partial x}\Big \|_{L^2(\mathcal{I}_t)}^{2} + b\Big\|\frac{ \partial^2 \rho}{\partial x^2}\Big\|_{L^2(\mathcal{I}_t)}^{2} = -c\int_{\mathcal{I}_t} u\frac{\partial^2 \rho}{\partial x^2} dx - d \Big\|\frac{\partial \rho}{\partial x}\Big \|_{L^2(\mathcal{I}_t)}^{2}   \\ \hspace{5.5cm} \displaystyle \leq  \frac{b}{2}\Big\|\frac{ \partial^2 \rho}{\partial x^2 }\Big\|_{L^2(\mathcal{I}_t)}^{2}  + \frac{c^2}{2b} \|u \|_{L^2(\mathcal{I}_t)}^{2} - d \Big\|\frac{\partial \rho}{\partial x}\Big \|_{L^2(\mathcal{I}_t)}^{2}.
    \end{array}
\end{align*}
Again, by \eqref{LpLq},
\begin{equation}\label{u2}
    \frac{c^2}{2b} \|u \|_{L^2(\mathcal{I}_t)}^{2} \leq \varepsilon\|u\|_{H^1(\mathcal{I}_t)}^2 +C_{\varepsilon}\|u\|_{L^1(\mathcal{I}_t)}^2
    \;\leq\;\varepsilon\|u\|_{H^1(\mathcal{I}_t)}^2 +C_{\varepsilon}M^2,
\end{equation}
where, here and below, $C_\varepsilon$ denotes the largest of the constants in \eqref{ab-w-001} and \eqref{u2}.

\smallskip
\noindent\emph{Zero-order term of $u$.} The first equation of \eqref{equ4} has no zero-order damping, so the $L^2$-norm of $u$ must be controlled through the mass. By \eqref{LpLq} with $p=1$, $q=2$ (so that $\sigma=1/3$), Young's inequality and \eqref{M0}, there is $C_g>0$ such that $\|u\|_{L^2(\mathcal I_t)}^2\leq\frac12\|u\|_{H^1(\mathcal I_t)}^2+C_gM^2$. Since $\|u\|_{H^1(\mathcal I_t)}^2=\|u\|_{L^2(\mathcal I_t)}^2+\|\partial_xu\|_{L^2(\mathcal I_t)}^2$, this gives
\begin{equation}\label{eq:u-mass}
\|u\|_{L^2(\mathcal I_t)}^2\leq\Big\|\frac{\partial u}{\partial x}\Big\|_{L^2(\mathcal I_t)}^2+2C_gM^2,
\qquad
\|u\|_{H^1(\mathcal I_t)}^2\leq2\Big\|\frac{\partial u}{\partial x}\Big\|_{L^2(\mathcal I_t)}^2+2C_gM^2 .
\end{equation}

\smallskip
\noindent\emph{Combination.} Multiplying \eqref{ab-w-002} and the estimate for $\partial_x\rho$ by $2$, adding them, and using \eqref{ab-w-001}, \eqref{u2} and $\|u\|_{L^1(\mathcal I_t)}\leq M$,
\begin{align*}
&\frac{d}{dt}\Big[\|u\|_{L^2(\mathcal I_t)}^2+\|\partial_x\rho\|_{L^2(\mathcal I_t)}^2\Big]
+2a\|\partial_xu\|_{L^2(\mathcal I_t)}^2+b\|\partial_x^2\rho\|_{L^2(\mathcal I_t)}^2+2d\|\partial_x\rho\|_{L^2(\mathcal I_t)}^2\\
&\hspace{2cm}\leq4\varepsilon\|u\|_{H^1(\mathcal I_t)}^2+2\varepsilon\|\rho\|_{H^2(\mathcal I_t)}^2+2C_\varepsilon M^4\|\rho\|_{H^1(\mathcal I_t)}^6+2C_\varepsilon M^2 .
\end{align*}
Choose $\varepsilon:=\min\{a,b,d\}/16$. By \eqref{eq:u-mass}, $4\varepsilon\|u\|_{H^1}^2\leq8\varepsilon\|\partial_xu\|^2+8\varepsilon C_gM^2\leq\frac a2\|\partial_xu\|^2+8\varepsilon C_gM^2$, and, writing $\|\rho\|_{H^2}^2=\|\rho\|^2+\|\partial_x\rho\|^2+\|\partial_x^2\rho\|^2$, $2\varepsilon\|\rho\|_{H^2}^2\leq2\varepsilon\|\rho\|^2+\frac d8\|\partial_x\rho\|^2+\frac b8\|\partial_x^2\rho\|^2$ (all norms on $\mathcal I_t$). Hence
\[
\frac{d}{dt}\Big[\|u\|^2+\|\partial_x\rho\|^2\Big]
+\frac{3a}{2}\|\partial_xu\|^2+\frac{7b}{8}\|\partial_x^2\rho\|^2+\frac{15d}{8}\|\partial_x\rho\|^2
\leq2\varepsilon\|\rho\|^2+2C_\varepsilon M^4\|\rho\|_{H^1}^6+(2C_\varepsilon+8\varepsilon C_g)M^2 .
\]
Set $\nu_{abd}:=\min\{a,d\}/2$ and $\nu'_{abd}:=\min\{a,b,d\}/4$. By \eqref{eq:u-mass},
$\frac a2\|\partial_xu\|^2\geq\nu_{abd}\|u\|^2-aC_gM^2$ and $a\|\partial_xu\|^2\geq\nu'_{abd}\|u\|_{H^1}^2-aC_gM^2$; moreover
$\frac{15d}{8}\|\partial_x\rho\|^2\geq(\nu_{abd}+\nu'_{abd})\|\partial_x\rho\|^2$, $\frac{7b}{8}\|\partial_x^2\rho\|^2\geq\nu'_{abd}\|\partial_x^2\rho\|^2$, and
$\nu'_{abd}\big(\|\partial_x\rho\|^2+\|\partial_x^2\rho\|^2\big)=\nu'_{abd}\big(\|\rho\|_{H^2}^2-\|\rho\|^2\big)$. Therefore,
\begin{align}\label{dphi-0}
\begin{array}{l}
  \displaystyle \frac{d}{dt}\Big[ \|u\|_{L^2(\mathcal{I}_t)}^{2} + \Big\|\frac{\partial \rho}{\partial x}\Big \|_{L^2(\mathcal{I}_t)}^{2}\Big] + \nu_{abd}\Big[ \|u\|_{L^2(\mathcal{I}_t)}^{2} + \Big\|\frac{\partial \rho}{\partial x}\Big \|_{L^2(\mathcal{I}_t)}^{2}\Big]  + \nu'_{abd}\Big[\|u\|_{H^1(\mathcal{I}_t)}^{2} + \|\rho\|_{H^2(\mathcal{I}_t)}^{2} \Big]\\ \displaystyle\leq  C_7\|\rho \|_{L^2(\mathcal{I}_t)}^{2}   +2C_{\varepsilon} M^4 \|\rho\|_{H^1(\mathcal{I}_t)}^{6}\;+\;C_8M^2,
  \end{array}
\end{align}
where $C_7:=2\varepsilon+\nu'_{abd}$ and $C_8:=2C_\varepsilon+8\varepsilon C_g+2aC_g$.

From \eqref{ab-001}, \eqref{abs-01} and \eqref{dphi-0}, the Gr\"onwall inequality gives
\begin{align*}
\|u(t)\|_{L^2(\mathcal I_t)}^2+\Big\|\frac{\partial\rho}{\partial x}(t)\Big\|_{L^2(\mathcal I_t)}^2
&\leq e^{-\nu_{abd}(t-\tau)}\Big(\|u_0\|_{L^2(\mathcal I_\tau)}^2
+\|\rho_0\|_{H^1(\mathcal I_\tau)}^2\Big)\\
&\hspace{0.5cm}
+C_7\int_\tau^t e^{-\nu_{abd}(t-s)}\|\rho(s)\|_{L^2(\mathcal I_s)}^2\,ds
+2C_\varepsilon M^4\int_\tau^t e^{-\nu_{abd}(t-s)}\|\rho(s)\|_{H^1(\mathcal I_s)}^6\,ds\\
&\hspace{0.5cm}+\frac{C_8}{\nu_{abd}}M^2 .
\end{align*}

Thus, by \eqref{ab-001}, we get
\begin{equation*}
C_7\int_\tau^t e^{-\nu_{abd}(t-s)}\|\rho(s)\|_{L^2(\mathcal I_s)}^2\,ds
\leq
C_7\|\rho_0\|_{L^2(\mathcal I_\tau)}^2\int_\tau^t e^{-\nu_{abd}(t-s)}e^{-\nu_d(s-\tau)}\,ds
+\frac{C_bC_7}{\nu_d}M^2\int_\tau^t e^{-\nu_{abd}(t-s)}\,ds .
\end{equation*}
With $\sigma=t-s$, $T=t-\tau$, since
$\nu_{abd}\sigma+\nu_d(T-\sigma)\geq\nu_1T$ for $\sigma\in[0,T]$, where
$\nu_1:=\min\{\nu_{abd},\nu_d\}$,
\begin{equation*}
\int_\tau^t e^{-\nu_{abd}(t-s)}e^{-\nu_d(s-\tau)}\,ds
=\int_0^T e^{-\nu_{abd}\sigma}e^{-\nu_d(T-\sigma)}\,d\sigma
\;\leq\;T\,e^{-\nu_1T} .
\end{equation*}
Fixing $\mu_1\in(0,\nu_1)$ and using
$Te^{-(\nu_1-\mu_1)T}\leq\frac{1}{e(\nu_1-\mu_1)}=:C_{\mu_1}$,
\begin{equation}\label{eq:cross-1}
C_7\int_\tau^t e^{-\nu_{abd}(t-s)}\|\rho(s)\|_{L^2(\mathcal I_s)}^2\,ds
\;\leq\;C_7C_{\mu_1}\|\rho_0\|_{L^2(\mathcal I_\tau)}^2e^{-\mu_1(t-\tau)}
+\frac{C_bC_7}{\nu_d\nu_{abd}}M^2 .
\end{equation}

 Now, by \eqref{abs-01}, we have
\begin{equation*}
\|\rho(s)\|_{H^1(\mathcal I_s)}^6
\;\leq\;32\Big[\widehat C_{0,\gamma}^6\|\rho_0\|_{H^1(\mathcal I_\tau)}^6
e^{-6\delta^\ast(s-\tau)}+\widehat C_{1,\gamma}^6M^6\Big] .
\end{equation*}
As above, with $\nu_2:=\min\{\nu_{abd},6\delta^\ast\}$,
$\mu_2\in(0,\nu_2)$, $C_{\mu_2}:=\frac{1}{e(\nu_2-\mu_2)}$,
\begin{equation}\label{eq:cross-2}
2C_\varepsilon M^4\int_\tau^t e^{-\nu_{abd}(t-s)}\|\rho(s)\|_{H^1(\mathcal I_s)}^6\,ds
\;\leq\;64\,C_\varepsilon M^4\widehat C_{0,\gamma}^6\|\rho_0\|_{H^1(\mathcal I_\tau)}^6\,C_{\mu_2}\,e^{-\mu_2(t-\tau)}
+\frac{64\,C_\varepsilon\widehat C_{1,\gamma}^6}{\nu_{abd}}M^{10} .
\end{equation}

Finally, combining \eqref{eq:cross-1} and
\eqref{eq:cross-2}, and setting $\delta_u:=\min\{\nu_{abd},\mu_1,\mu_2\}>0$,
\begin{equation}\label{fw-03-new}
\begin{aligned}
\|u(t)\|_{L^2(\mathcal I_t)}^2+\Big\|\frac{\partial\rho}{\partial x}(t)\Big\|_{L^2(\mathcal I_t)}^2
\;\leq\;{}&
C_{1,M}\Big(\|u_0\|_{L^2(\mathcal I_\tau)}^2+\|\rho_0\|_{H^1(\mathcal I_\tau)}^2
+\|\rho_0\|_{L^2(\mathcal I_\tau)}^2\Big)e^{-\delta_u(t-\tau)}\\
&+C_{2,M}\,\|\rho_0\|_{H^1(\mathcal I_\tau)}^6\,e^{-\delta_u(t-\tau)}
+\frac{C_bC_7}{\nu_d\nu_{abd}}M^2
+\frac{64\,C_\varepsilon\widehat C_{1,\gamma}^6}{\nu_{abd}}M^{10}
\;+\;\frac{C_8}{\nu_{abd}}M^2,
\end{aligned}
\end{equation}
for all $\tau<t$ and $(u_0,\rho_0)\in\widehat D(\tau)$, where
\begin{equation*}
C_{1,M}:=1+C_7C_{\mu_1},
\qquad
C_{2,M}:=64\,C_\varepsilon M^4\widehat C_{0,\gamma}^6\,C_{\mu_2} .
\end{equation*}
In \eqref{fw-03-new} all dependence on the initial data $(u_0,\rho_0)$ is confined to the two exponentially decaying terms, while the last three terms depend only on $M$. Since $\widehat D\in\mathcal D_F(\mathbb
H_1^M)$ is a bounded family, there is $R>0$ with $\|u_0\|_{L^2(\mathcal I_\tau)},\|\rho_0\|_{H^1(\mathcal I_\tau)}\leq R$ for all $(u_0,\rho_0)\in\widehat D(\tau)$ and all $\tau$; hence there exists $\tau_1(t,\widehat D)<t$ such that the decaying terms in \eqref{fw-03-new} are $\leq1$ for all $\tau\leq\tau_1(t,\widehat D)$, giving
\begin{equation}\label{fw-04}
\|u(t)\|_{L^2(\mathcal I_t)}\;\leq\;R_2(M):=\sqrt{1+\frac{C_bC_7}{\nu_d\nu_{abd}}M^2
+\frac{64\,C_\varepsilon\widehat C_{1,\gamma}^6}{\nu_{abd}}M^{10}+\frac{C_8}{\nu_{abd}}M^2},
\end{equation}
for all $\tau\leq\tau_1(t,\widehat D)$ and $(u_0,\rho_0)\in\widehat D(\tau)$.

Similarly, by \eqref{abs-01}, since $\widehat C_{1,\gamma}M$ depends
only on $M$, there exists $\tau_2(t,\widehat D)<t$ such that
$\widehat C_{0,\gamma}\|\rho_0\|_{H^1(\mathcal I_\tau)}e^{-\delta^\ast(t-\tau)}\leq1$
for all $\tau\leq\tau_2(t,\widehat D)$, whence
\begin{equation}\label{fw-04b}
\|\rho(t)\|_{H^1(\mathcal I_t)}\;\leq\;R_1(M):=1+\widehat C_{1,\gamma}M,
\qquad\forall\,\tau\leq\tau_2(t,\widehat D),\ (u_0,\rho_0)\in\widehat D(\tau).
\end{equation}

Combining \eqref{fw-04} and \eqref{fw-04b}, and setting
$\tau_\ast(t,\widehat D):=\min\{\tau_1(t,\widehat D),\tau_2(t,\widehat D)\}$,
\begin{equation*}
\big\|(u(t),\rho(t))\big\|_{L^2(\mathcal I_t)\times H^1(\mathcal{I}_t)}
\;\leq\; R_M:=\max\{R_1(M),R_2(M)\},
\end{equation*}
for all $\tau\leq\tau_\ast(t,\widehat D)$ and $(u_0,\rho_0)\in\widehat
D(\tau)$. Since $R_M$ depends only on $M$, this proves that
$U_{1}$ is dissipative.
\end{proof}

\begin{remark}\label{decae-0} For any integer $q\geq1$ and any bound of the form $x(s)\leq Xe^{-r(s-\tau)}+Y$ ($X,Y,r>0$), the elementary inequality
$(x+y)^q\leq2^{q-1}(x^q+y^q)$ gives
$x(s)^q\leq2^{q-1}\big(X^qe^{-qr(s-\tau)}+Y^q\big)$, and consequently,
exactly as in the proof of Lemma \ref{lemmaR11}, we have
\[
\int_\tau^te^{-\varrho(t-s)}x(s)^q\,ds
\;\leq\;2^{q-1}X^qC_{\mu}\,e^{-\mu(t-\tau)}+\frac{2^{q-1}Y^q}{\varrho},
\]
for any $\mu\in(0,\min\{\varrho,qr\})$, with $C_\mu=\big[e(\min\{\varrho,qr\}-\mu)\big]^{-1}$. 
\end{remark}

\begin{lemma}\label{AsympCL}
Suppose that conditions (H1)-(H4) hold. Then the evolution process
$U_{1}$, given in \eqref{PU}, is pullback
$\mathcal D_F(\mathbb H_1^M)$-asymptotically compact.
\end{lemma}
\begin{proof} We divide the proof into three steps. In Step~1 we obtain an estimate in $H^1(\mathcal I_t)\times H^2(\mathcal I_t)$ for regular initial data;
in Step~2 we extend it, by a smoothing argument, to initial data in $L^2(\mathcal I_\tau)\times H^1(\mathcal I_\tau)$; and in Step~3 we conclude the
asymptotic compactness. 

\vspace{0,5cm}

\noindent $\bullet$ {\bf Step 1}: We first assume that the initial data satisfy $(u_0,\rho_0) \in H^2_N(\mathcal{I}_\tau)\times H^3_N(\mathcal{I}_\tau)$.
By the regularity established in Theorem~\ref{teoEx-05-34}, we have
$\frac{\partial\rho}{\partial t}\in L^2(\tau,T;H^2_N(\mathcal I_t))$ and $\rho\in L^2(\tau,T;H^4_{N^2}(\mathcal I_t))$,
so that we may differentiate the second equation of \eqref{equ4} with respect to $x$, obtaining
\begin{equation}\label{rho-x}
\frac{\partial}{\partial t}\Big(\frac{\partial\rho}{\partial x}\Big)
=b\,\frac{\partial^3\rho}{\partial x^3}+c\,\frac{\partial u}{\partial x}-d\,\frac{\partial\rho}{\partial x}
\quad\text{in } \mathcal I_t ,
\end{equation}
and the traces of $\frac{\partial^2\rho}{\partial t\,\partial x}$ on $x=\alpha(t)$ and $x=\beta(t)$ are well defined.
Since $\frac{\partial\rho}{\partial x}(\alpha(t),t)=\frac{\partial\rho}{\partial x}(\beta(t),t)=0$ for all $t\geq\tau$,
differentiating these identities with respect to $t$ gives
\begin{equation}\label{rho-xt-bd}
\frac{\partial^2\rho}{\partial t\,\partial x}(\alpha(t),t)=-\alpha'(t)\frac{\partial^2\rho}{\partial x^2}(\alpha(t),t),
\qquad
\frac{\partial^2\rho}{\partial t\,\partial x}(\beta(t),t)=-\beta'(t)\frac{\partial^2\rho}{\partial x^2}(\beta(t),t).
\end{equation}
Multiplying \eqref{rho-x} by $-\partial_x^3\rho$, integrating over $\mathcal I_t$, and integrating by parts
on the left-hand side with the help of \eqref{rho-xt-bd}, we get
\[
-\int_{\alpha(t)}^{\beta(t)}\frac{\partial^2\rho}{\partial t\,\partial x}\frac{\partial^3\rho}{\partial x^3}\,dx
=\int_{\alpha(t)}^{\beta(t)}\frac{\partial^3\rho}{\partial t\,\partial x^2}\frac{\partial^2\rho}{\partial x^2}\,dx
+\beta'(t)\Big|\frac{\partial^2\rho}{\partial x^2}(\beta(t),t)\Big|^2-\alpha'(t)\Big|\frac{\partial^2\rho}{\partial x^2}(\alpha(t),t)\Big|^2 .
\]
On the other hand, by the Reynolds transport theorem,
\[
\frac12\frac{d}{dt}\Big\|\frac{\partial^2\rho}{\partial x^2}\Big\|_{L^2(\mathcal I_t)}^2
=\int_{\alpha(t)}^{\beta(t)}\frac{\partial^3\rho}{\partial t\,\partial x^2}\frac{\partial^2\rho}{\partial x^2}\,dx
+\frac12\Big(\beta'(t)\Big|\frac{\partial^2\rho}{\partial x^2}(\beta(t),t)\Big|^2-\alpha'(t)\Big|\frac{\partial^2\rho}{\partial x^2}(\alpha(t),t)\Big|^2\Big).
\]
Moreover, by the Neumann boundary condition,
$\int_{\alpha(t)}^{\beta(t)}\partial_x\rho\,\partial_x^3\rho\,dx=-\|\partial_x^2\rho\|_{L^2(\mathcal I_t)}^2$.
Therefore, the following energy equality holds:
\begin{align*}
\frac12\frac{d}{dt}\Big\|\frac{\partial^2\rho}{\partial x^2}\Big\|_{L^2(\mathcal I_t)}^2
+b\Big\|\frac{\partial^3\rho}{\partial x^3}\Big\|_{L^2(\mathcal I_t)}^2
+d\Big\|\frac{\partial^2\rho}{\partial x^2}\Big\|_{L^2(\mathcal I_t)}^2
=&-c\int_{\alpha(t)}^{\beta(t)}\frac{\partial u}{\partial x}\frac{\partial^3\rho}{\partial x^3}\,dx\\
&-\frac12\Big(\beta'(t)\Big|\frac{\partial^2\rho}{\partial x^2}(\beta(t),t)\Big|^2-\alpha'(t)\Big|\frac{\partial^2\rho}{\partial x^2}(\alpha(t),t)\Big|^2\Big).
\end{align*}
For the first term on the right-hand side, Young's inequality gives
\begin{align*}
-c\int_{\alpha(t)}^{\beta(t)}\frac{\partial u}{\partial x}\frac{\partial^3\rho}{\partial x^3}\,dx
&\leq\frac b2\Big\|\frac{\partial^3\rho}{\partial x^3}\Big\|_{L^2(\mathcal I_t)}^2
+\frac{c^2}{2b}\Big\|\frac{\partial u}{\partial x}\Big\|_{L^2(\mathcal I_t)}^2\\
&\leq\frac b2\Big\|\frac{\partial^3\rho}{\partial x^3}\Big\|_{L^2(\mathcal I_t)}^2
+\varepsilon\|u\|_{H^2(\mathcal I_t)}^2+C_\varepsilon\|u\|_{L^2(\mathcal I_t)}^2,
\end{align*}
where in the last line we used the interpolation inequality
$\|\partial_xu\|_{L^2(\mathcal I_t)}^2\leq\varepsilon\|u\|_{H^2(\mathcal I_t)}^2
+C_\varepsilon\|u\|_{L^2(\mathcal I_t)}^2$ (with $\varepsilon$ replaced by
$2b\varepsilon/c^2$). 

To estimate the boundary term, we note that, by (H3), $L:=\sup_{t\in\mathbb R}\max\{|\alpha'(t)|,|\beta'(t)|\}<\infty.$ By the Agmon-type inequality
$\|g\|_{L^\infty(\mathcal I_t)}^2\leq C\|g\|_{L^2(\mathcal I_t)}\|g\|_{H^1(\mathcal I_t)}$
with $g=\partial_x^2\rho$, the interpolation inequality
$\|\rho\|_{H^2(\mathcal I_t)}\leq C\|\rho\|_{H^3(\mathcal I_t)}^{1/2}\|\rho\|_{H^1(\mathcal I_t)}^{1/2}$
and Young's inequality (with exponents $4/3$ and $4$), we get
\begin{align*}
\frac12\Big|\beta'(t)\Big|\frac{\partial^2\rho}{\partial x^2}(\beta(t),t)\Big|^2-\alpha'(t)\Big|\frac{\partial^2\rho}{\partial x^2}(\alpha(t),t)\Big|^2\Big|
&\leq L\Big\|\frac{\partial^2\rho}{\partial x^2}\Big\|_{L^\infty(\mathcal I_t)}^2
\leq CL\|\rho\|_{H^2(\mathcal I_t)}\|\rho\|_{H^3(\mathcal I_t)}\\
&\leq CL\|\rho\|_{H^3(\mathcal I_t)}^{3/2}\|\rho\|_{H^1(\mathcal I_t)}^{1/2}
\leq\varepsilon\|\rho\|_{H^3(\mathcal I_t)}^2+C_\varepsilon\|\rho\|_{H^1(\mathcal I_t)}^2 .
\end{align*}
Then, we obtain
\begin{equation}\label{rho2x2}
\begin{aligned}
\frac{d}{dt}\Big\|\frac{\partial^2\rho}{\partial x^2}\Big\|_{L^2(\mathcal I_t)}^2
+b\Big\|\frac{\partial^3\rho}{\partial x^3}\Big\|_{L^2(\mathcal I_t)}^2
&+2d\Big\|\frac{\partial^2\rho}{\partial x^2}\Big\|_{L^2(\mathcal I_t)}^2\\
&\leq2\varepsilon\big(\|u\|_{H^2(\mathcal I_t)}^2+\|\rho\|_{H^3(\mathcal I_t)}^2\big)
+2C_\varepsilon\big(\|u\|_{L^2(\mathcal I_t)}^2+\|\rho\|_{H^1(\mathcal I_t)}^2\big) .
\end{aligned}
\end{equation}
Similarly, multiplying the first equation of \eqref{equ4} by $-\partial_x^2 u$
and integrating over $\mathcal{I}_t$, we obtain the energy equality
\begin{equation}\label{u2x2}
\begin{array}{l}
\displaystyle \frac12\frac{d}{dt}\Big\|\frac{\partial u}{\partial x}\Big\|_{L^2(\mathcal I_t)}^2
+a\Big\|\frac{\partial^2u}{\partial x^2}\Big\|_{L^2(\mathcal I_t)}^2
=\int_{\alpha(t)}^{\beta(t)}\frac{\partial}{\partial x}
\Big(u\frac{\partial\rho}{\partial x}\Big)\frac{\partial^2u}{\partial x^2}\,dx
\\ \hspace{5.5cm} \displaystyle\leq\;\frac a2\Big\|\frac{\partial^2u}{\partial x^2}\Big\|_{L^2(\mathcal I_t)}^2
+\frac{1}{2a}\int_{\alpha(t)}^{\beta(t)}
\Big|\frac{\partial}{\partial x}\Big(u\frac{\partial\rho}{\partial x}\Big)\Big|^2\,dx .
\end{array}
\end{equation}

By H\"older's inequality and $H^1(\mathcal{I}_t) \hookrightarrow L^\infty(\mathcal{I}_t)$, we have
\begin{align*}
\int_{\alpha(t)}^{\beta(t)}\Big|\frac{\partial}{\partial x}
\Big(u\frac{\partial\rho}{\partial x}\Big)\Big|^2dx
&\leq C\Big(\|u\|_{H^1(\mathcal I_t)}^2\Big\|\frac{\partial\rho}{\partial x}\Big\|_{L^\infty(\mathcal I_t)}^2
+\|u\|_{L^2(\mathcal I_t)}^2\Big\|\frac{\partial^2\rho}{\partial x^2}\Big\|_{L^\infty(\mathcal I_t)}^2\Big) .
\end{align*}
By the one-dimensional interpolation (Agmon-type) inequality
$\|g\|_{L^\infty(\mathcal I_t)}^2\leq C\|g\|_{L^2(\mathcal I_t)}\|g\|_{H^1(\mathcal I_t)}$
applied to $g=\partial_x\rho$ and $g=\partial_x^2\rho$,
\begin{equation*}
\Big\|\frac{\partial\rho}{\partial x}\Big\|_{L^\infty(\mathcal I_t)}^2\leq C\|\rho\|_{H^1(\mathcal I_t)}\|\rho\|_{H^2(\mathcal I_t)},
\qquad
\Big\|\frac{\partial^2\rho}{\partial x^2}\Big\|_{L^\infty(\mathcal I_t)}^2\leq C\|\rho\|_{H^2(\mathcal I_t)}\|\rho\|_{H^3(\mathcal I_t)},
\end{equation*}
so that
\begin{equation}\label{u2x2-pre}
\int_{\alpha(t)}^{\beta(t)}\Big|\frac{\partial}{\partial x}
\Big(u\frac{\partial\rho}{\partial x}\Big)\Big|^2dx
\;\leq\;C\Big(\|u\|_{H^1(\mathcal I_t)}^2\|\rho\|_{H^2(\mathcal I_t)}\|\rho\|_{H^1(\mathcal I_t)}
+\|u\|_{L^2(\mathcal I_t)}^2\|\rho\|_{H^3(\mathcal I_t)}\|\rho\|_{H^2(\mathcal I_t)}\Big) .
\end{equation}
By the interpolation inequalities
$\|u\|_{H^1(\mathcal I_t)}^2\leq C\|u\|_{H^2(\mathcal I_t)}\|u\|_{L^2(\mathcal I_t)}$
and $\|\rho\|_{H^2(\mathcal I_t)}\leq C\|\rho\|_{H^3(\mathcal I_t)}^{1/2}\|\rho\|_{H^1(\mathcal I_t)}^{1/2}$,
the two terms on the right-hand side of \eqref{u2x2-pre} are bounded by
\begin{equation*}
C\Big(\underbrace{\|u\|_{H^2(\mathcal I_t)}\|u\|_{L^2(\mathcal I_t)}\|\rho\|_{H^3(\mathcal I_t)}^{1/2}\|\rho\|_{H^1(\mathcal I_t)}^{3/2}}_{=:A}
+\underbrace{\|u\|_{L^2(\mathcal I_t)}^2\|\rho\|_{H^3(\mathcal I_t)}^{3/2}\|\rho\|_{H^1(\mathcal I_t)}^{1/2}}_{=:B}\Big) .
\end{equation*}

\noindent $\circ$ \emph{Estimate of $A$.} Applying Young's inequality twice in
succession, first to isolate $\|u\|_{H^2(\mathcal I_t)}$ and then to
isolate the remaining factor of $\|\rho\|_{H^3(\mathcal I_t)}$, we get
\[
\|u\|_{H^2(\mathcal I_t)}\Big(\|u\|_{L^2(\mathcal I_t)}\|\rho\|_{H^3(\mathcal I_t)}^{1/2}\|\rho\|_{H^1(\mathcal I_t)}^{3/2}\Big)
\leq\varepsilon\|u\|_{H^2(\mathcal I_t)}^2
+C_\varepsilon\|u\|_{L^2(\mathcal I_t)}^2\|\rho\|_{H^3(\mathcal I_t)}\|\rho\|_{H^1(\mathcal I_t)}^3,
\]
\[
C_\varepsilon\|\rho\|_{H^3(\mathcal I_t)}\Big(\|u\|_{L^2(\mathcal I_t)}^2\|\rho\|_{H^1(\mathcal I_t)}^3\Big)
\leq\varepsilon\|\rho\|_{H^3(\mathcal I_t)}^2+C_\varepsilon'\|u\|_{L^2(\mathcal I_t)}^4\|\rho\|_{H^1(\mathcal I_t)}^6,
\]
so that
\begin{equation*}
A\;\leq\;\varepsilon\|u\|_{H^2(\mathcal I_t)}^2+\varepsilon\|\rho\|_{H^3(\mathcal I_t)}^2
+C_\varepsilon\|u\|_{L^2(\mathcal I_t)}^4\|\rho\|_{H^1(\mathcal I_t)}^6 .
\end{equation*}

\noindent $\circ$ \emph{Estimate of $B$.} Writing
$\|\rho\|_{H^3(\mathcal I_t)}^{3/2}=\|\rho\|_{H^3(\mathcal I_t)}\cdot\|\rho\|_{H^3(\mathcal I_t)}^{1/2}$
and applying Young's inequality twice, as for $A$ (first
isolating one power of $\|\rho\|_{H^3(\mathcal I_t)}$, then the
remaining one), we get
\[
\|\rho\|_{H^3(\mathcal I_t)}\Big(\|u\|_{L^2(\mathcal I_t)}^2\|\rho\|_{H^3(\mathcal I_t)}^{1/2}\|\rho\|_{H^1(\mathcal I_t)}^{1/2}\Big)
\leq\varepsilon\|\rho\|_{H^3(\mathcal I_t)}^2+C_\varepsilon\|u\|_{L^2(\mathcal I_t)}^4\|\rho\|_{H^3(\mathcal I_t)}\|\rho\|_{H^1(\mathcal I_t)},
\]
\[
C_\varepsilon\|\rho\|_{H^3(\mathcal I_t)}\Big(\|u\|_{L^2(\mathcal I_t)}^4\|\rho\|_{H^1(\mathcal I_t)}\Big)
\leq\varepsilon\|\rho\|_{H^3(\mathcal I_t)}^2+C_\varepsilon'\|u\|_{L^2(\mathcal I_t)}^8\|\rho\|_{H^1(\mathcal I_t)}^2,
\]
so that
\begin{equation}\label{eq:B-bound}
B\;\leq\;2\varepsilon\|\rho\|_{H^3(\mathcal I_t)}^2+C_\varepsilon\|u\|_{L^2(\mathcal I_t)}^8\|\rho\|_{H^1(\mathcal I_t)}^{2} .
\end{equation}

Combining \eqref{u2x2-pre} and \eqref{eq:B-bound}, we obtain
\begin{equation}\label{u2x2-01}
\int_{\alpha(t)}^{\beta(t)}\Big|\frac{\partial}{\partial x}
\Big(u\frac{\partial\rho}{\partial x}\Big)\Big|^2dx
\;\leq\;\varepsilon\big(\|u\|_{H^2(\mathcal I_t)}^2+\|\rho\|_{H^3(\mathcal I_t)}^2\big)
+C_\varepsilon\Big(\|u\|_{L^2(\mathcal I_t)}^4\|\rho\|_{H^1(\mathcal I_t)}^6
+\|u\|_{L^2(\mathcal I_t)}^8\|\rho\|_{H^1(\mathcal I_t)}^{2}\Big) .
\end{equation}

Now we add \eqref{rho2x2} and twice \eqref{u2x2}, and use \eqref{u2x2-01}. Since
\[
\|u\|_{H^2(\mathcal I_t)}^2\leq C\Big(\Big\|\frac{\partial^2u}{\partial x^2}\Big\|_{L^2(\mathcal I_t)}^2+\|u\|_{L^2(\mathcal I_t)}^2\Big)
\quad\text{and}\quad
\|\rho\|_{H^3(\mathcal I_t)}^2\leq C\Big(\Big\|\frac{\partial^3\rho}{\partial x^3}\Big\|_{L^2(\mathcal I_t)}^2+\|\rho\|_{L^2(\mathcal I_t)}^2\Big),
\]
we can choose $\varepsilon>0$ small enough that the terms
$\varepsilon\|u\|_{H^2(\mathcal I_t)}^2$ and $\varepsilon\|\rho\|_{H^3(\mathcal I_t)}^2$
are absorbed by the coercive terms $a\|\partial_x^2u\|_{L^2(\mathcal I_t)}^2$ and
$b\|\partial_x^3\rho\|_{L^2(\mathcal I_t)}^2$, the lower-order terms being moved to
the right-hand side. Thus, there exist constants $C_{ab},\widehat C>0$ such that
\begin{equation*}
\begin{aligned}
\frac{d}{dt}\Big(\Big\|\frac{\partial u}{\partial x}\Big\|_{L^2(\mathcal I_t)}^2
+\Big\|\frac{\partial^2\rho}{\partial x^2}\Big\|_{L^2(\mathcal I_t)}^2\Big)
&+C_{ab}\big(\|u\|_{H^2(\mathcal I_t)}^2+\|\rho\|_{H^3(\mathcal I_t)}^2\big)\\
&\leq\widehat C\Big(\|u\|_{L^2(\mathcal I_t)}^2+\|\rho\|_{H^1(\mathcal I_t)}^2
+\|u\|_{L^2(\mathcal I_t)}^4\|\rho\|_{H^1(\mathcal I_t)}^6
+\|u\|_{L^2(\mathcal I_t)}^8\|\rho\|_{H^1(\mathcal I_t)}^{2}\Big) .
\end{aligned}
\end{equation*}
Applying Young's inequality once more to each product term on the right-hand
side gives the fully separated form
\begin{equation}\label{u2x2-rho2x2-b}
\frac{d}{dt}\Big(\Big\|\frac{\partial u}{\partial x}\Big\|_{L^2(\mathcal I_t)}^2
+\Big\|\frac{\partial^2\rho}{\partial x^2}\Big\|_{L^2(\mathcal I_t)}^2\Big)
+C_{ab}\big(\|u\|_{H^2(\mathcal I_t)}^2+\|\rho\|_{H^3(\mathcal I_t)}^2\big)
\;\leq\;\widehat C\,\Xi(t),
\end{equation}
where
\begin{equation*}
\Xi(t):=\|u\|_{L^2(\mathcal I_t)}^2+\|u\|_{L^2(\mathcal I_t)}^8+\|u\|_{L^2(\mathcal I_t)}^{16}
+\|\rho\|_{H^1(\mathcal I_t)}^{2}+\|\rho\|_{H^1(\mathcal I_t)}^{4}+\|\rho\|_{H^1(\mathcal I_t)}^{12} .
\end{equation*}

Since $\|u\|_{H^2(\mathcal I_t)}^2+\|\rho\|_{H^3(\mathcal I_t)}^2\geq
\|\partial_xu\|_{L^2(\mathcal I_t)}^2+\|\partial_x^2\rho\|_{L^2(\mathcal I_t)}^2$,
the Gr\"onwall inequality applied to \eqref{u2x2-rho2x2-b} on the interval $[s,t]$ gives
\begin{equation}\label{eq:duhamel}
\begin{aligned}
\Big\|\frac{\partial u}{\partial x}(t)\Big\|_{L^2(\mathcal I_t)}^2
+\Big\|\frac{\partial^2\rho}{\partial x^2}(t)\Big\|_{L^2(\mathcal I_t)}^2
\;\leq\;&e^{-C_{ab}(t-s)}\Big(\Big\|\frac{\partial u}{\partial x}(s)\Big\|_{L^2(\mathcal I_s)}^2
+\Big\|\frac{\partial^2\rho}{\partial x^2}(s)\Big\|_{L^2(\mathcal I_s)}^2\Big)\\
&+\widehat C\int_\tau^t e^{-C_{ab}(t-r)}\,\Xi(r)\,dr ,
\end{aligned}
\end{equation}
for all $\tau\leq s\leq t$, where we also used that $\Xi\geq0$ to replace the
integral over $(s,t)$ by the integral over $(\tau,t)$.

To estimate the last integral we use the bounds already established for
$\|u(r)\|_{L^2(\mathcal I_r)}$ and $\|\rho(r)\|_{H^1(\mathcal I_r)}$.
Recall (from the proof of Lemma~\ref{lemmaR11}, inequalities \eqref{abs-01} and \eqref{fw-03-new}) that
\begin{equation}\label{eq:recall-bounds}
\begin{array}{l}
\displaystyle\|u(r)\|_{L^2(\mathcal I_r)}^2\leq\mathcal{Y}\big[\|u_0\|_{L^2(\mathcal{I}_\tau)},\|\rho_0\|_{H^1(\mathcal{I}_\tau)}\big] e^{-\delta_u(r-\tau)}+K(M),\\\\ \displaystyle
\|\rho(r)\|_{H^1(\mathcal I_r)}\leq\widehat C_{0,\gamma}\|\rho_0\|_{H^1(\mathcal I_\tau)}e^{-\delta^\ast(r-\tau)}+\widehat C_{1,\gamma}M,
\end{array}
\end{equation}
for all $r\geq\tau$, where $\mathcal{Y}\big[\|u_0\|_{L^2(\mathcal{I}_\tau)},\|\rho_0\|_{H^1(\mathcal{I}_\tau)}\big]$
collects the (decaying) initial-data dependence and $K(M)$ depends only on $M$
(see \eqref{fw-03-new}).

Fix $\delta_2\in\big(0,\min\{C_{ab},\delta_u,2\delta^\ast\}\big)$. We apply
Remark~\ref{decae-0}, with $\varrho=C_{ab}$ and $\mu=\delta_2$, to each of the six
terms in $\Xi(r)$: to $x(r)=\|u(r)\|_{L^2(\mathcal I_r)}^2$ with $q=1,4,8$, using the
first bound in \eqref{eq:recall-bounds} (rate $\delta_u$), and to
$x(r)=\|\rho(r)\|_{H^1(\mathcal I_r)}$ with $q=2,4,12$, using the second bound in
\eqref{eq:recall-bounds} (rate $\delta^\ast$). In all cases
$\min\{C_{ab},q\delta_u\}$ and $\min\{C_{ab},q\delta^\ast\}$ are greater than or
equal to $\min\{C_{ab},\delta_u,2\delta^\ast\}>\delta_2$, so the constant given by
Remark~\ref{decae-0} is bounded by
$C_\ast:=\big[e\big(\min\{C_{ab},\delta_u,2\delta^\ast\}-\delta_2\big)\big]^{-1}$.
Since $\|u\|_{L^2}^8=(\|u\|_{L^2}^2)^4$ and $\|u\|_{L^2}^{16}=(\|u\|_{L^2}^2)^8$, we obtain
\begin{equation}\label{eq:integral-final}
\widehat C\int_\tau^te^{-C_{ab}(t-r)}\Xi(r)\,dr
\;\leq\;\mathcal{J}\big[\|u_0\|_{L^2(\mathcal{I}_\tau)},\|\rho_0\|_{H^1(\mathcal{I}_\tau)}\big]\,e^{-\delta_2(t-\tau)}+K_2(M),
\end{equation}
where
\begin{equation*}
\begin{aligned}
\mathcal{J}\big[\|u_0\|_{L^2(\mathcal{I}_\tau)},\|\rho_0\|_{H^1(\mathcal{I}_\tau)}\big]
:=\widehat C\,C_\ast\Big(&\mathcal{Y}\big[\|u_0\|_{L^2(\mathcal{I}_\tau)},\|\rho_0\|_{H^1(\mathcal{I}_\tau)}\big]
+8\,\mathcal{Y}\big[\|u_0\|_{L^2(\mathcal{I}_\tau)},\|\rho_0\|_{H^1(\mathcal{I}_\tau)}\big]^4\\
&+128\,\mathcal{Y}\big[\|u_0\|_{L^2(\mathcal{I}_\tau)},\|\rho_0\|_{H^1(\mathcal{I}_\tau)}\big]^8
+2\widehat C_{0,\gamma}^2\|\rho_0\|_{H^1(\mathcal I_\tau)}^2\\
&+8\widehat C_{0,\gamma}^4\|\rho_0\|_{H^1(\mathcal I_\tau)}^4
+2^{11}\widehat C_{0,\gamma}^{12}\|\rho_0\|_{H^1(\mathcal I_\tau)}^{12}\Big)
\end{aligned}
\end{equation*}
and
\begin{equation*}
K_2(M):=\frac{\widehat C}{C_{ab}}\Big(K(M)+8K(M)^4+128K(M)^8
+2\big(\widehat C_{1,\gamma}M\big)^{2}+8\big(\widehat C_{1,\gamma}M\big)^{4}+2^{11}\big(\widehat C_{1,\gamma}M\big)^{12}\Big)
\end{equation*}
depends \emph{only} on $M$. Since $\mathcal Y$ is continuous and bounded on
bounded sets, so is $\mathcal J$.

Finally, by \eqref{eq:recall-bounds} with $r=t$,
\[
\|u(t)\|_{L^2(\mathcal I_t)}^2+\|\rho(t)\|_{H^1(\mathcal I_t)}^2
\leq\Big(\mathcal{Y}\big[\|u_0\|_{L^2(\mathcal{I}_\tau)},\|\rho_0\|_{H^1(\mathcal{I}_\tau)}\big]
+2\widehat C_{0,\gamma}^2\|\rho_0\|_{H^1(\mathcal I_\tau)}^2\Big)e^{-\delta_2(t-\tau)}
+K(M)+2\widehat C_{1,\gamma}^2M^2 ,
\]
since $\delta_2\leq\min\{\delta_u,2\delta^\ast\}$. Adding this inequality to
\eqref{eq:duhamel}, using \eqref{eq:integral-final}, and still denoting by
$\mathcal J$ and $K_2(M)$ the functions obtained by adding
$\mathcal{Y}\big[\|u_0\|_{L^2(\mathcal{I}_\tau)},\|\rho_0\|_{H^1(\mathcal{I}_\tau)}\big]
+2\widehat C_{0,\gamma}^2\|\rho_0\|_{H^1(\mathcal I_\tau)}^2$ to $\mathcal J$ and
$K(M)+2\widehat C_{1,\gamma}^2M^2$ to $K_2(M)$, we conclude that
\begin{equation}\label{ab-uH1-rhoH2}
\begin{aligned}
\|u(t)\|_{H^1(\mathcal I_t)}^2+\|\rho(t)\|_{H^2(\mathcal I_t)}^2
\;\leq\;&\big(\|u(s)\|_{H^1(\mathcal I_s)}^2+\|\rho(s)\|_{H^2(\mathcal I_s)}^2\big)e^{-C_{ab}(t-s)}\\
&+\mathcal{J}\big[\|u_0\|_{L^2(\mathcal{I}_\tau)},\|\rho_0\|_{H^1(\mathcal{I}_\tau)}\big]e^{-\delta_2(t-\tau)}
+K_2(M),
\end{aligned}
\end{equation}
for all $\tau\leq s\leq t$.

\vspace{0,5cm}

\noindent $\bullet$ {\bf Step 2}: Consider now $(u_0,\rho_0) \in L^2(\mathcal{I}_\tau)\times H^1(\mathcal{I}_\tau)$.
By density, there exists a sequence $(u_{0,n},\rho_{0,n}) \in H^2_N(\mathcal{I}_\tau)\times H^3_N(\mathcal{I}_\tau)$ such that
\[  (u_{0,n},\rho_{0,n}) \longrightarrow (u_{0},\rho_{0}) \quad \text{ as } \quad n\to \infty, \quad \text{ in } L^2(\mathcal{I}_\tau)\times H^1(\mathcal{I}_\tau) .\]
Let $(u_n(t),\rho_n(t))= U_1(t,\tau)(u_{0,n},\rho_{0,n})$. By Step~1, we have
\begin{equation}\label{eq:step2-n}
\begin{aligned}
\|u_n(t)\|_{H^1(\mathcal I_t)}^2+\|\rho_n(t)\|_{H^2(\mathcal I_t)}^2
\;\leq\;&\big(\|u_n(s)\|_{H^1(\mathcal I_s)}^2+\|\rho_n(s)\|_{H^2(\mathcal I_s)}^2\big)e^{-C_{ab}(t-s)}\\
&+\mathcal{J}\big[\|u_{0,n}\|_{L^2(\mathcal{I}_\tau)},\|\rho_{0,n}\|_{H^1(\mathcal{I}_\tau)}\big]e^{-\delta_2(t-\tau)}
+K_2(M),
\end{aligned}
\end{equation}
for all $\tau\leq s\leq t$. The first term on the right-hand side is not
controlled by the norm of the initial data in
$L^2(\mathcal I_\tau)\times H^1(\mathcal I_\tau)$. To remove it, we integrate
\eqref{eq:step2-n} with respect to $s$, and for this we first estimate
the integral of $\|u_n(s)\|_{H^1(\mathcal I_s)}^2+\|\rho_n(s)\|_{H^2(\mathcal I_s)}^2$.

\noindent $\circ$ \emph{Integral estimate.} Set $\tilde\nu:=\min\{\nu_{abd},C_{ab}\}$.
Since $\nu_{abd}\geq\tilde\nu$, inequality \eqref{dphi-0} for $(u_n,\rho_n)$ yields
\begin{equation}\label{eq:dphi-n}
\begin{array}{l}
\displaystyle \frac{d}{ds}\Big[ \|u_n\|_{L^2(\mathcal{I}_s)}^{2} + \Big\|\frac{\partial \rho_n}{\partial x}\Big \|_{L^2(\mathcal{I}_s)}^{2}\Big]
+ \tilde\nu\Big[ \|u_n\|_{L^2(\mathcal{I}_s)}^{2} + \Big\|\frac{\partial \rho_n}{\partial x}\Big \|_{L^2(\mathcal{I}_s)}^{2}\Big]
\\ \displaystyle \hspace{5cm}+ \nu'_{abd}\Big[\|u_n\|_{H^1(\mathcal{I}_s)}^{2} + \|\rho_n\|_{H^2(\mathcal{I}_s)}^{2} \Big]\\ \displaystyle \hspace{5cm}
\leq  C_7\|\rho_n \|_{L^2(\mathcal{I}_s)}^{2}   +2C_{\varepsilon} M^4 \|\rho_n\|_{H^1(\mathcal{I}_s)}^{6}.
\end{array}
\end{equation}
By the second inequality in \eqref{eq:recall-bounds}, together with
$\|\rho_n\|_{L^2(\mathcal I_s)}\leq\|\rho_n\|_{H^1(\mathcal I_s)}$ and
$(x+y)^p\leq 2^{p-1}(x^p+y^p)$, we have
\[
\|\rho_n(s)\|_{H^1(\mathcal I_s)}^2\leq 2\widehat C_{0,\gamma}^2\|\rho_{0,n}\|_{H^1(\mathcal I_\tau)}^2e^{-2\delta^\ast(s-\tau)}+2\widehat C_{1,\gamma}^2M^2,
\]
\[
\|\rho_n(s)\|_{H^1(\mathcal I_s)}^6\leq 32\widehat C_{0,\gamma}^6\|\rho_{0,n}\|_{H^1(\mathcal I_\tau)}^6e^{-6\delta^\ast(s-\tau)}+32\widehat C_{1,\gamma}^6M^6 .
\]
Hence, since $e^{-6\delta^\ast(s-\tau)}\leq e^{-2\delta^\ast(s-\tau)}$, the right-hand side of
\eqref{eq:dphi-n} is bounded by
\[
\Big(2C_7\widehat C_{0,\gamma}^2\|\rho_{0,n}\|_{H^1(\mathcal I_\tau)}^2
+64C_\varepsilon\widehat C_{0,\gamma}^6M^4\|\rho_{0,n}\|_{H^1(\mathcal I_\tau)}^6\Big)e^{-2\delta^\ast(s-\tau)}
+K_1(M),
\]
where $K_1(M):=2C_7\widehat C_{1,\gamma}^2M^2+64C_\varepsilon\widehat C_{1,\gamma}^6M^{10}$.
Multiplying \eqref{eq:dphi-n} by $e^{\tilde\nu(s-\tau)}$, integrating over $(\tau,t)$,
multiplying the result by $e^{-\tilde\nu(t-\tau)}$, discarding the
(nonnegative) term at time $t$, and using
$\|u_{0,n}\|_{L^2(\mathcal I_\tau)}^2+\|\partial_x\rho_{0,n}\|_{L^2(\mathcal I_\tau)}^2
\leq\|u_{0,n}\|_{L^2(\mathcal I_\tau)}^2+\|\rho_{0,n}\|_{H^1(\mathcal I_\tau)}^2$, we obtain
\begin{align*}
\nu'_{abd}\int_\tau^t e^{-\tilde\nu(t-s)}&\Big[\|u_n(s)\|_{H^1(\mathcal I_s)}^2+\|\rho_n(s)\|_{H^2(\mathcal I_s)}^2\Big]ds
\leq\Big(\|u_{0,n}\|_{L^2(\mathcal I_\tau)}^2+\|\rho_{0,n}\|_{H^1(\mathcal I_\tau)}^2\Big)e^{-\tilde\nu(t-\tau)}\\
&+\int_\tau^t e^{-\tilde\nu(t-s)}\Big(C_7\|\rho_n(s)\|_{L^2(\mathcal{I}_s)}^{2}+2C_{\varepsilon}M^4\|\rho_n(s)\|_{H^1(\mathcal{I}_s)}^{6}\Big)ds.
\end{align*}
Fix $\delta_3\in\big(0,\min\{\tilde\nu,2\delta^\ast\}\big)$ and set
$C_{\delta_3}:=\big[e\big(\min\{\tilde\nu,2\delta^\ast\}-\delta_3\big)\big]^{-1}$. By the bound obtained
above for the integrand, Remark~\ref{decae-0} (with $q=1$, $\varrho=\tilde\nu$ and
$\mu=\delta_3$) gives
\begin{align*}
\int_\tau^t e^{-\tilde\nu(t-s)}&\Big(C_7\|\rho_n(s)\|_{L^2(\mathcal{I}_s)}^{2}+2C_{\varepsilon}M^4\|\rho_n(s)\|_{H^1(\mathcal{I}_s)}^{6}\Big)ds\\
&\leq C_{\delta_3}\Big(2C_7\widehat C_{0,\gamma}^2\|\rho_{0,n}\|_{H^1(\mathcal I_\tau)}^2
+64C_\varepsilon\widehat C_{0,\gamma}^6M^4\|\rho_{0,n}\|_{H^1(\mathcal I_\tau)}^6\Big)e^{-\delta_3(t-\tau)}
+\frac{K_1(M)}{\tilde\nu}.
\end{align*}
Consequently, since $e^{-\tilde\nu(t-\tau)}\leq e^{-\delta_3(t-\tau)}$, for all $t>\tau$,
\begin{equation}\label{eq:weighted-integral}
\begin{aligned}
\int_\tau^t & e^{-\tilde\nu(t-s)}\Big[\|u_n(s)\|_{H^1(\mathcal I_s)}^2+\|\rho_n(s)\|_{H^2(\mathcal I_s)}^2\Big]ds \\& \leq\frac{1}{\nu'_{abd}}\Big[\Big(\|u_{0,n}\|_{L^2(\mathcal I_\tau)}^2
+\big(1+2C_7\widehat C_{0,\gamma}^2C_{\delta_3}\big)\|\rho_{0,n}\|_{H^1(\mathcal I_\tau)}^2
\\& +64C_\varepsilon\widehat C_{0,\gamma}^6M^4C_{\delta_3}\|\rho_{0,n}\|_{H^1(\mathcal I_\tau)}^6\Big)e^{-\delta_3(t-\tau)}
+\frac{K_1(M)}{\tilde\nu}\Big].
\end{aligned}
\end{equation}

\noindent $\circ$ \emph{Estimate for the approximations.} Since $C_{ab}\geq\tilde\nu$, we have
$e^{-C_{ab}(t-s)}\leq e^{-\tilde\nu(t-s)}$. Integrating \eqref{eq:step2-n} with respect to
$s\in(\tau,t)$ (note that its left-hand side does not depend on $s$) gives
\begin{align*}
(t-\tau)\Big(\|u_n(t)\|_{H^1(\mathcal I_t)}^2+\|\rho_n(t)\|_{H^2(\mathcal I_t)}^2\Big)
\leq\;&\int_\tau^t e^{-\tilde\nu(t-s)}\Big[\|u_n(s)\|_{H^1(\mathcal I_s)}^2+\|\rho_n(s)\|_{H^2(\mathcal I_s)}^2\Big]ds\\
&+(t-\tau)\,\mathcal{J}\big[\|u_{0,n}\|_{L^2(\mathcal{I}_\tau)},\|\rho_{0,n}\|_{H^1(\mathcal{I}_\tau)}\big]e^{-\delta_2(t-\tau)}
+(t-\tau)K_2(M).
\end{align*}
Dividing by $t-\tau$ and using \eqref{eq:weighted-integral}, we conclude that, for all $t>\tau$ and all $n\in\mathbb N$,
\begin{equation}\label{eq:En-bound}
\begin{aligned}
\|u_n(t)\|_{H^1(\mathcal I_t)}^2+\|\rho_n(t)\|_{H^2(\mathcal I_t)}^2
&\leq\;\frac{1}{\nu'_{abd}(t-\tau)}\Big[\Big(\|u_{0,n}\|_{L^2(\mathcal I_\tau)}^2
+\big(1+2C_7\widehat C_{0,\gamma}^2C_{\delta_3}\big)\|\rho_{0,n}\|_{H^1(\mathcal I_\tau)}^2
\\&+64C_\varepsilon\widehat C_{0,\gamma}^6M^4C_{\delta_3}\|\rho_{0,n}\|_{H^1(\mathcal I_\tau)}^6\Big)e^{-\delta_3(t-\tau)}
+\frac{K_1(M)}{\tilde\nu}\Big]\\
&+\mathcal{J}\big[\|u_{0,n}\|_{L^2(\mathcal{I}_\tau)},\|\rho_{0,n}\|_{H^1(\mathcal{I}_\tau)}\big]e^{-\delta_2(t-\tau)}+K_2(M).
\end{aligned}
\end{equation}

\noindent $\circ$ \emph{Passage to the limit.} Fix $t>\tau$. Since $(u_{0,n},\rho_{0,n})\to(u_0,\rho_0)$ in $L^2(\mathcal I_\tau)\times H^1(\mathcal I_\tau)$,
we have $\|u_{0,n}\|_{L^2(\mathcal I_\tau)}\to\|u_0\|_{L^2(\mathcal I_\tau)}$ and
$\|\rho_{0,n}\|_{H^1(\mathcal I_\tau)}\to\|\rho_0\|_{H^1(\mathcal I_\tau)}$; since
$\mathcal J$ is continuous, the right-hand side of \eqref{eq:En-bound} converges to the same expression with $(u_{0,n},\rho_{0,n})$ replaced by $(u_0,\rho_0)$. In particular, the sequence $(u_n(t),\rho_n(t))$ is bounded in $H^1(\mathcal I_t)\times H^2(\mathcal I_t)$, and therefore there exists a subsequence, which we still denote by $(u_n(t),\rho_n(t))$, that converges weakly in $H^1(\mathcal I_t)\times H^2(\mathcal I_t)$. On the other hand, by the continuous dependence on the initial data, $(u_n(t),\rho_n(t))\to(u(t),\rho(t))=U_1(t,\tau)(u_0,\rho_0)$ in
$L^2(\mathcal I_t)\times H^1(\mathcal I_t)$. By the uniqueness of the weak limit in $L^2(\mathcal I_t)\times H^1(\mathcal I_t)$, the weak limit in
$H^1(\mathcal I_t)\times H^2(\mathcal I_t)$ is $(u(t),\rho(t))$. Since the norm is weakly lower semicontinuous,
\[
\|u(t)\|_{H^1(\mathcal I_t)}^2+\|\rho(t)\|_{H^2(\mathcal I_t)}^2
\leq\liminf_{n\to\infty}\Big(\|u_n(t)\|_{H^1(\mathcal I_t)}^2+\|\rho_n(t)\|_{H^2(\mathcal I_t)}^2\Big).
\]
Passing to the limit in \eqref{eq:En-bound}, we conclude that, for all
$(u_0,\rho_0)\in L^2(\mathcal I_\tau)\times H^1(\mathcal I_\tau)$ and all $t>\tau$,
\begin{equation}\label{eq:decay-L2H1}
\begin{aligned}
\|u(t)\|_{H^1(\mathcal I_t)}^2+\|\rho(t)\|_{H^2(\mathcal I_t)}^2 &
\leq\; \frac{1}{\nu'_{abd}(t-\tau)}\Big[\Big(\|u_{0}\|_{L^2(\mathcal I_\tau)}^2
+\big(1+2C_7\widehat C_{0,\gamma}^2C_{\delta_3}\big)\|\rho_{0}\|_{H^1(\mathcal I_\tau)}^2
\\&+64C_\varepsilon\widehat C_{0,\gamma}^6M^4C_{\delta_3}\|\rho_{0}\|_{H^1(\mathcal I_\tau)}^6\Big)e^{-\delta_3(t-\tau)}
+\frac{K_1(M)}{\tilde\nu}\Big]\\
&+\mathcal{J}\big[\|u_{0}\|_{L^2(\mathcal{I}_\tau)},\|\rho_{0}\|_{H^1(\mathcal{I}_\tau)}\big]e^{-\delta_2(t-\tau)}+K_2(M).
\end{aligned}
\end{equation}

\vspace{0,5cm}

\noindent $\bullet$ {\bf Step 3}: Finally, let $t\in\mathbb R$ and $\widehat D\in\mathcal D_F(\mathbb H_1^M)$, and consider a sequence of initial data $\{(u_{\tau_n},\rho_{\tau_n})\}_{n\in\mathbb N}$
with $(u_{\tau_n},\rho_{\tau_n})\in\widehat D(\tau_n)$ for all $n$ and $\tau_n\to-\infty$. Since $\widehat D$ is a bounded family, there exists $R>0$ such that
\[
\|u_{\tau_n}\|_{L^2(\mathcal I_{\tau_n})}\leq R\quad\text{and}\quad
\|\rho_{\tau_n}\|_{H^1(\mathcal I_{\tau_n})}\leq R,\qquad\forall\,n\in\mathbb N .
\]
Let $n_0\in\mathbb N$ be such that $t-\tau_n\geq1$ for all $n\geq n_0$. Applying \eqref{eq:decay-L2H1} with $\tau=\tau_n$, we obtain, for all $n\geq n_0$,
\begin{align*}
\|U_1(t,\tau_n)(u_{\tau_n},\rho_{\tau_n})\|_{H^1(\mathcal I_t)\times H^2(\mathcal I_t)}^2
\leq\;&\frac{1}{\nu'_{abd}}\Big(\big(2+2C_7\widehat C_{0,\gamma}^2C_{\delta_3}\big)R^2
+64C_\varepsilon\widehat C_{0,\gamma}^6M^4C_{\delta_3}R^6\Big)e^{-\delta_3(t-\tau_n)}\\
&+\sup_{0\leq x,y\leq R}\mathcal J[x,y]\,e^{-\delta_2(t-\tau_n)}
+\frac{K_1(M)}{\nu'_{abd}\tilde\nu}+K_2(M),
\end{align*}
where the supremum is finite because $\mathcal J$ is bounded on bounded sets. Since $\tau_n\to-\infty$, the first two terms on the right-hand side tend to zero, so there exists $m_0\geq n_0$ such that
\[
\|U_1(t,\tau_n)(u_{\tau_n},\rho_{\tau_n})\|_{H^1(\mathcal I_t)\times H^2(\mathcal I_t)}^2
\;\leq\;1+\frac{K_1(M)}{\nu'_{abd}\tilde\nu}+K_2(M),
\qquad\forall\,n\geq m_0 .
\]
In particular, $\{U_1(t,\tau_n)(u_{\tau_n},\rho_{\tau_n})\}_{n\geq m_0}$ is bounded in $H^1(\mathcal I_t)\times H^2(\mathcal I_t)$. Since
$H^1(\mathcal I_t)\times H^2(\mathcal I_t)$ is compactly embedded in $L^2(\mathcal I_t)\times H^1(\mathcal I_t)$, the sequence $\{U_1(t,\tau_n)(u_{\tau_n},\rho_{\tau_n})\}_{n\geq m_0}$ admits a convergent subsequence in $L^2(\mathcal I_t)\times H^1(\mathcal I_t)$. This proves that $U_{1}$ is pullback $\mathcal D_F(\mathbb H_1^M)$-asymptotically compact.
\end{proof}

\begin{theorem}\label{Existence-P-A-in-M2H}
Suppose that conditions (H1)-(H4) hold. Then, there exist the minimal pullback $\mathcal{D}_{F}(\mathbb{H}^M_1)$-attractor $\mathcal{A}_{\mathcal{D}_{F}(\mathbb{H}^M_1)}=\{\mathcal{A}_{\mathcal{D}_{F}(\mathbb{H}^M_1)}(t):t\in\mathbb{R}\}$ for the process
$U_{1}$. The following relationships hold
\begin{equation*}\label{Inc-1-PAM2H}
\begin{array}{cc}
    & \mathcal{A}_{\mathcal{D}_{F}(\mathbb{H}^M_1)}(t)\subset B_{\mathbb{H}_{1,t}^M}[0,R_M]
\;\;\;\forall t \in \mathbb{R}.
\end{array}
\end{equation*}
Moreover, $\mathcal{A}_{\mathcal{D}_{F}(\mathbb{H}^M_1)}\in\mathcal{D}_{F}(\mathbb{H}_1^M)$.
\end{theorem}
\begin{proof} The asymptotic compactness of process $U_{1}$ is given by Lemma \ref{AsympCL} and the existence of an absorbing pullback family follows from Lemma \ref{lemmaR11}. Therefore, the existence of pullback attractors for process $U_{1}$ in the universes $\mathcal{D}_{F}(\mathbb{H}_1^M)$, follows from Theorem \ref{t4}.
\end{proof}

\begin{remark}
The quantity
\[
\int_{\mathcal I_t}u(x)\,dx
\]
represents the total mass of the cell population at time $t$. Therefore, the phase space
$\mathbb{H}_{1,t}^M$ describes all admissible states whose total population does not
exceed the threshold $M$ at time $t$.

If $N>M$, then every state with total mass bounded by $M$ is also admissible for
the larger threshold $N$. Consequently,
\[
\mathcal{D}_F(\mathbb{H}_1^M)
\subset
\mathcal{D}_F(\mathbb{H}_1^N),
\]
and Theorem \ref{compara} yields
\[
\mathcal{A}_{\mathcal{D}_F(\mathbb{H}_1^M)}(t)
\subset
\mathcal{A}_{\mathcal{D}_F(\mathbb{H}_1^N)}(t),
\qquad \forall \ t\in\mathbb{R}.
\]
From a biological viewpoint, this means that allowing a larger total cell mass
enlarges the collection of possible long-term spatial distributions. Every
asymptotic regime attainable by populations with total mass at most $M$
remains attainable when the admissible mass is increased to $N$, while
additional asymptotic behaviors may arise due to the presence of larger
populations.
\end{remark}

\begin{remark}\label{rem:L1-estimate}
Integrating the first equation of \eqref{equ5} over $\mathcal I$ and using the homogeneous Neumann boundary conditions, we obtain
\[
\frac{d}{dt}\int_{\mathcal I}v(y,t)\,dy
=
-\int_{\mathcal I}b(y,t)\,\frac{\partial v}{\partial y}(y,t)\,dy.
\]
Now, integration by parts yields
\[
\frac{d}{dt}\int_{\mathcal I}v(y,t)\,dy
=
-\frac{\gamma'(t)}{\gamma(t)}
\int_{\mathcal I}v(y,t)\,dy
+
\frac{\beta'(t)v(1,t)-\alpha'(t)v(0,t)}{\gamma(t)}.
\]
Consequently,
\[
\frac{d}{dt}
\left(
\gamma(t)\int_{\mathcal I}v(y,t)\,dy
\right)
=
\beta'(t)v(1,t)-\alpha'(t)v(0,t).
\]
By (H2) and since $v\geq 0$, we have $\beta'(t)v(1,t)-\alpha'(t)v(0,t)\le0$, then
\[
\gamma(t)\int_{\mathcal I}v(y,t)\,dy
\le
\gamma(\tau)\int_{\mathcal I}v_\tau(y)\,dy,
\qquad t\ge\tau.
\]
Hence, the weighted mass $\gamma(t)\int_{\mathcal I}v(y,t)\,dy$ is non-increasing along trajectories.
\end{remark}

Now, consider
\begin{equation*}
    \mathbb{H}_{2,t}^M:=\bigg\{ v\in L^2(\mathcal{I}): v\geq 0  \ \text{ and } \ \gamma(t)\int_{\mathcal{I}}v(y)dy \leq M\bigg\}\times \Big\{\phi \in H^1(\mathcal{I}): \phi\geq 0\Big\},
\end{equation*}
with the topology induced by $\mathcal{H} = L^2 (\mathcal{I}) \times H^1(\mathcal{I})$.

By Remark \ref{rem:L1-estimate} and Lemma \ref{lem:positivity}, the weighted mass
$\gamma(t)\int_{\mathcal I}v(y,t)\,dy$ is non-increasing along nonnegative solutions of
\eqref{equ5}; hence we can define  the bi-parametric families of mappings $\{U_2(t,\tau): \mathbb{H}_{2,\tau}^M \to \mathbb{H}_{2,t}^M:(t,\tau)\in \mathbb{R}_{d}^{2}\}$ 
\begin{equation*}
 U_2(t,\tau)(v_0,\phi_0)=\big(v(t;\tau,v_0,\phi_0),\phi(\cdot;\tau,v_0,\phi_0)\big),
\end{equation*}
where $(v,\phi)$ is the unique weak solution of \eqref{equ5} associated to the initial condition $(v_0,\phi_0)\in \mathbb{H}_{2,\tau}^M$ for any $(t,\tau)\in \mathbb{R}_{d}^{2}$.

On the other hand, taking into account Lemma \ref{lem:scaling}, for each $t\in \mathbb{R}$ we can define the following map 
\begin{align}\label{kappa-N} 
\begin{array}{l}
   \displaystyle \Theta_1(t) \hspace{0,3cm}: \hspace{0,3cm} \mathbb{H}_{1,t}^{M} \hspace{0,3cm} \longrightarrow \hspace{0,3cm} \hspace{0,3cm} \mathbb{H}_{2,t}^{M} \\\\ \displaystyle \hspace{1.3cm} (u,\rho) \hspace{0,1cm} \longmapsto  \hspace{0,3cm}\big(u(r^{-1}(\cdot,t)),\rho(r^{-1}(\cdot,t))\big).
   \end{array}
\end{align}
Thus, it is not difficult to see that $\{\Theta_1(t):\mathbb{H}_{1,t}^{M} \to \mathbb{H}_{2,t}^{M}:t\in \mathbb{R}\}$ is a family of homeomorphisms.

\begin{proposition}\label{universo-limitados} Let $\mathcal{D}_{\mathbb{H}_{1}^{M}}$ be a closed inclusion universes on the family $\{ \mathbb{H}_{1,t}^{M} : t \in \mathbb{R}\}$. If the universes of families of bounded sets in $\{\mathbb{H}_{2,t}^{M}: t \in \mathbb{R}\}$, $\mathcal{D}_{F}(\mathbb{H}_{2}^{M})$, and $\mathcal{D}_{\mathbb{H}_{1}^{M}}$ are equivalents. Then 
    \begin{equation*}
      \widehat D \in  \mathcal{D}_{\mathbb{H}_{1}^{M}} \hspace{0,3cm }\text{ if only if } \hspace{0,3cm } \widehat D=\{D(t) \subset B_{\mathbb{H}_{1,t}^{M}}(0,R):t\in\mathbb R\} \text{ for some } R >0.
    \end{equation*}
   Therefore, $\mathcal{D}_F(\mathbb{H}_{1}^{M})=\mathcal{D}_{\mathbb{H}_{1}^{M}}$ and the universes $\mathcal{D}_F(\mathbb{H}_{1}^{M})$ and $\mathcal{D}_F(\mathbb{H}_{2}^{M})$ are equivalents.
\end{proposition}
\begin{proof}
 Let $\{\Theta_1(t): \mathbb{H}_{1,t}^M \to \mathbb{H}_{2,t}^{M} : t\in \mathbb{R}\}$ be the family of homeomorphisms defined in \eqref{kappa-N}. It is not difficult to show that there exist positive constants $\kappa_1 , \kappa_2$ with
\begin{equation}\label{ineq-XY}
   \kappa_{1}\|(v,\psi)\|_{\mathcal{H}}\leq \|\Theta_1^{-1}(t)(v,\psi)\|_{\mathcal{H}_t} \leq \kappa_2 \|(v,\psi)\|_{\mathcal{H}},
\end{equation}
for all $(v,\psi)\in \mathbb{H}_{2,t}^M$. Now, let $\widehat D \in  \mathcal{D}_{\mathbb{H}_{1}^{M}}$. Since the universes $\mathcal{D}_{F}(\mathbb{H}_{2}^M)$ and $\mathcal{D}_{\mathbb{H}_{1}^{M}}$ are equivalents, we have $\Theta_{1}(t)D(t) \subset B(t)$ for all $t \in \mathbb{R}$, where $\widehat{B}=\{B(t)\subset \mathbb{H}_{2,t}^{M}:t \in \mathbb{R}\}$ is a backward bounded family of $\{\mathbb{H}_{2,t}^{M}: t \in \mathbb{R}\}$. Thus, there exists $R>0$ such that $B(t) \subset B_{\mathbb{H}_{2,t}^{M}}(0,R)$ for all $t\in \mathbb{R}$ and by \eqref{ineq-XY} we obtain 
 \begin{equation*}
     \|(u,\rho)\|_{\mathcal{H}_t} \leq \kappa_2 \|\Theta_{1}(t)(u,\rho)\|_{\mathcal{H}}, \text{ for all } (u,\rho) \in B(t), \ t\in \mathbb{R}.
 \end{equation*}
 Therefore $D(t) \subset B_{\mathbb{H}_{1,t}^{M}}(0,\kappa_2 R)$ for all $t \in \mathbb{R}$.

 On the other hand, consider the family $\widehat D=\big\{D(t) \subset B_{\mathbb{H}_{1,t}^{M}}(0,R):t\in\mathbb R \big\} \in \mathcal{D}_F (\mathbb{H}_{1}^{M})$ for some $R>0$. By \eqref{ineq-XY} we obtain
  \begin{equation}\label{dX-00}
    \|\Theta_{1}(t) (u,\rho)\|_{\mathcal{H}} \leq \frac{1}{\kappa_1} \|(u,\rho )\|_{\mathcal{H}_t}, \text{ for all } (u,\rho) \in D(t) \text{ and } t\in \mathbb{R}.
 \end{equation}
 
 Then, the family $\widehat{B}=\{B(t) \subset \mathbb{H}_{2,t}^{M} : t \in \mathbb{R}\}$, where $B(t)=\{ \Theta_{1}(t) (u,\rho) : (u,\rho) \in D(t), \ t\in \mathbb{R}\} $, belong to $\mathcal{D}_{F}(\mathbb{H}_{2}^{M})$, since, by \eqref{dX-00}, $B(t) \subset B_{\mathbb{H}_{2,t}^{M}}(0, R/\kappa_1)$. Finally, by the equivalence of universes $\mathcal{D}_{F}(\mathbb{H}_{2}^{M})$ and $\mathcal{D}_{\mathbb{H}_{1}^{M}}$ and the fact that $\mathcal{D}_{\mathbb{H}_{1}^{M}}$ is closed inclusion, we conclude that $\widehat{D} \in \mathcal{D}_{\mathbb{H}_{1}^{M}}$.
\end{proof}

Using Proposition \ref{universo-limitados}, we establish the following result, which guarantees the existence of a pullback attractor associated with problem \eqref{equ5} on the family of phase space $\{\mathbb{H}_{2,t}^{M}: t\in\mathbb{R}\}.$

\begin{theorem}\label{Existence-P-A-in-IKS2}
Suppose that conditions (H1)-(H4) hold. Then, there exists a minimal pullback $\mathcal{D}_{F}(\mathbb{H}_{2}^{M})$-attractor $\mathcal{A}_{\mathcal{D}_{F}(\mathbb{H}_{2}^{M})}=\{\mathcal{A}_{\mathcal{D}_{F}(\mathbb{H}_{2}^{M})}(t):t\in\mathbb{R}\}$ for the process
$U_2$, given by
\begin{equation*}
    \mathcal{A}_{\mathcal{D}_{F}(\mathbb{H}_{2}^{M})}(t) = \Theta_{1}(t) \mathcal{A}_{\mathcal{D}_{F}(\mathbb{H}_{1}^M)}(t)\;\;\;\forall t \in \mathbb{R},
\end{equation*}
where $\mathcal{A}_{\mathcal{D}_{F}(\mathbb{H}_{1}^{M})}$ is the minimal pullback $\mathcal{D}_{F}(\mathbb{H}_{1}^{M})$-attractor for the process $U_1$, given in Theorem \ref{Existence-P-A-in-M2H}. Moreover, $\mathcal{A}_{\mathcal{D}_{F}(\mathbb{H}_{2}^{M})} \in\mathcal{D}_{F}(\mathbb{H}_{2}^{M}))$.
\end{theorem}
\begin{proof}
    The existence of a minimal pullback $\mathcal{D}_{F}(\mathbb{H}_{2}^{M})$-attractor $\mathcal{A}_{\mathcal{D}_{F}(\mathbb{H}_{2}^{M})}$ for the process
$U_2$, given by
\begin{equation*}
    \mathcal{A}_{\mathcal{D}_{F}(\mathbb{H}_{2}^{M})}(t) = \Theta_{1}(t) \mathcal{A}_{\mathcal{D}_{F}(\mathbb{H}_{1}^M)}(t)\;\;\;\forall t \in \mathbb{R},
\end{equation*}
 follows from Theorem \ref{t5}, since the processes $U_1$ and $U_2$, as well as the universes
$\mathcal{D}_{F}(\mathbb{H}_{1}^{M})$ and
$\mathcal{D}_{F}(\mathbb{H}_{2}^{M})$, are equivalent. Moreover, $\mathcal{A}_{\mathcal{D}_{F}(\mathbb{H}_{2}^{M})} \in\mathcal{D}_{F}(\mathbb{H}_{2}^{M})$, since $\mathcal{A}_{\mathcal{D}_{F}(\mathbb{H}_{1}^{M})} \in\mathcal{D}_{F}(\mathbb{H}_{1}^{M})$.
\end{proof}

\begin{remark}
  In the same way as the pullback attractor for problem \eqref{equ5} was obtained in Theorem \ref{Existence-P-A-in-IKS2}, the time rescaling defined by \eqref{rescaling} allows us to establish the existence of a pullback attractor for problem \eqref{equ4iii}. More precisely, we have the following result.

Consider the family of phase spaces $\{\mathbb H^M_{3,\sigma}:\sigma\in\mathbb{R}\}$,
indexed by the rescaled time $\sigma$,
\[
\mathbb H^M_{3,\sigma}:=\Big\{v\in L^2(\mathcal I): v\ge0,\ \gamma(t(\sigma))\int_{\mathcal I} v\,dy\le M\Big\}
\times\{\phi\in H^1(\mathcal I):\phi\ge0\},
\]
with the topology induced by $\mathcal H=L^2(\mathcal I)\times H^1(\mathcal I)$. Since
\eqref{equ4iii} does not depend on the initial time, we can define
\[
U_3(s,\tau)(w_0,\psi_0)=\big(w(s;\tau,w_0,\psi_0),\psi(s;\tau,w_0,\psi_0)\big),
\qquad s\geq \tau,
\]
where $(w,\psi)$ is the unique weak solution of \eqref{equ4iii} with initial time $\tau$
and initial datum $(w_0,\psi_0)\in\mathbb H^M_{3,\sigma_0}$. For each
$s\in\mathbb{R}$ define
\[
\Theta_2(s):\mathbb H^M_{1,t(s)}\to\mathbb H^M_{3,s},\qquad
(u,\rho)\mapsto\big(u\circ r^{-1}(\cdot,t(s)),\ \rho\circ r^{-1}(\cdot,t(s))\big).
\]
It is not difficult to see that $\{\Theta_2(\sigma)\}_{\sigma\in\mathbb{R}}$ is a family of homeomorphisms satisfying
\[
U_3(s,\tau)=\Theta_2(s)\,U_1\big(t(s),t(\tau)\big)\,\Theta_2(\tau)^{-1}.
\]
Moreover, arguing as in Proposition~\ref{universo-limitados}, the universes
$\mathcal D_F(\mathbb H^M_1)$ and $\mathcal D_F(\mathbb H^M_3)$ are equivalent.
Therefore, by Theorem~\ref{t5}, there exists a minimal pullback
$\mathcal D_F(\mathbb H^M_3)$-attractor for $U_3$, given by
\[
\mathcal A_{\mathcal D_F(\mathbb H^M_3)}(\sigma)
=\Theta_2(\sigma)\,\mathcal A_{\mathcal D_F(\mathbb H^M_1)}(t(\sigma)),\qquad \sigma\in\mathbb{R},
\]
where $\mathcal{A}_{\mathcal{D}_{F}(\mathbb{H}_{1}^{M})}$ is the minimal pullback $\mathcal{D}_{F}(\mathbb{H}_{1}^{M})$-attractor for the process $U_1$, given in Theorem \ref{Existence-P-A-in-M2H}.
\end{remark}

Finally, we summarize the main results of this section by means of the following commutative diagrams.

\begin{figure}[ht]
\centering

\begin{minipage}{0.4\textwidth}
\centering
\[
\begin{tikzcd}[column sep=large,row sep=large]
(\mathbb H_1^M,U_1)
\arrow[r,"\Theta_1"]
\arrow[dr,swap,"\Theta_2"]
&
(\mathbb H_2^M,U_2)
\arrow[d,"\Theta_2\circ\Theta_1^{-1}"]
\\
&
(\mathbb H_3^M,U_3)
\end{tikzcd}
\]
\vspace{-0.2cm}

{\small (a) Equivalent processes}
\end{minipage}
\hfill
\begin{minipage}{0.4\textwidth}
\centering
\[
\begin{tikzcd}[column sep=large,row sep=large]
\mathcal A_{\mathcal D_F(\mathbb H_1^M)}
\arrow[r,"\Theta_1"]
\arrow[dr,swap,"\Theta_2"]
&
\mathcal A_{\mathcal D_F(\mathbb H_2^M)}
\arrow[d,"\Theta_2\circ\Theta_1^{-1}"]
\\
&
\mathcal A_{\mathcal D_F(\mathbb H_3^M)}
\end{tikzcd}
\]
\vspace{-0.2cm}

{\small (b) Corresponding pullback attractors}
\end{minipage}

\caption{Equivalence of the processes and the induced correspondence between their pullback attractors through the families of homeomorphisms.}
\label{fig:diagram}
\end{figure}
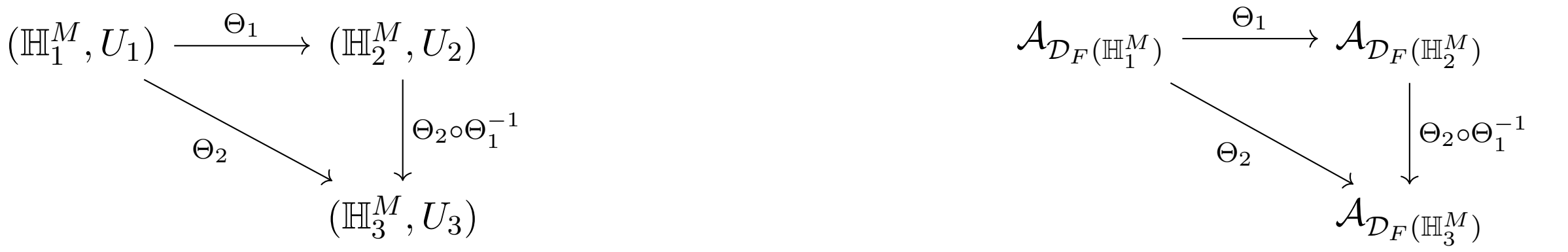

\section{Numerical Examples}\label{sec5}

Our numerical examples follow those carried out in \cite{yahagi} for the one-dimensional Keller--Segel system on a fixed interval. We reproduce analogous scenarios with the essential difference that the domain $\mathcal{I}_t=(\alpha(t),\beta(t))$ now contracts over time according to (H2), which allows us to compare the behavior for different contracting boundaries and to observe the effect of the boundary motion on cell aggregation. In both examples, problem \eqref{equ4} is solved numerically through its equivalent formulation \eqref{equ5} on $\mathcal{I}=[0,1]$, by means of a finite-difference discretization of the spatial variable $y$, with the homogeneous Neumann boundary conditions enforced through ghost nodes, coupled with an implicit backward differentiation formula (BDF) scheme for the time integration, as implemented in SciPy's {\it solve$\_$ivp}. The solutions $(u,\rho)$ on the original moving domain $\mathcal{I}_t$ are then recovered from $(v,\psi)$ through the inverse change of variables $x=\alpha(t)+\gamma(t)y$.

\noindent {\bf Example I.} We consider the chemotaxis problem \eqref{equ4}, where the constants $a,b,c,d,\kappa>0$, for the numerical example, are given as
\begin{equation}\label{eq:constants}
    a=1,\qquad b=1,\qquad \kappa=3,\qquad c=1,\qquad d=1.
\end{equation}
The boundary curves are
\begin{equation}\label{eq:alphabeta}
    \alpha(t) = -10+8\tanh(kt), \qquad \beta(t) = 10-8\tanh(kt), \qquad
    k := \frac{1}{6}\operatorname{arctanh}\!\left(\frac{15}{16}\right),
\end{equation}
so that $\gamma(t) = \beta(t)-\alpha(t) = 20-16\tanh(kt).$ The constant $k$ is chosen so that $\tanh(6k)=15/16$, i.e.\ $\gamma(0)=20$ and $\gamma(6)=5$: the domain $\mathcal{I}_t$ contracts from $[-10,10]$ at $t=0$ to $[-2.5,2.5]$ at $t=6$, and $\gamma(t)\to 4$ as $t\to+\infty$. In particular
\begin{equation}
    \alpha'(t) = \frac{8k}{\cosh^{2}(kt)}>0, \quad \beta'(t) = -\frac{8k}{\cosh^{2}(kt)}<0, \quad
    \gamma'(t) = -\frac{16k}{\cosh^{2}(kt)}<0.
\end{equation}

The initial data are given by: 
\begin{equation}\label{eq:utau-rhotau}
u_0(x) =
\begin{cases}
1-\cos(x), & -2\pi\leq x\leq 2\pi,\\[1mm]
0, & \text{otherwise in } [-10,10],
\end{cases}
\qquad
\rho_0(x) =
\begin{cases}
1+\cos(x), & -\pi\leq x\leq \pi,\\[1mm]
0, & \text{otherwise in } [-10,10].
\end{cases}
\end{equation}

Note that the curves $\alpha,\beta$ describing the moving boundaries satisfy conditions \emph{(H1)} and \emph{(H2)}, and that the initial data are nonnegative, i.e., $u_0\geq0$ and $\rho_0\geq0$. Hence, all the hypotheses of the global existence theorem given in Theorem \ref{teoEx-05-It} are satisfied.

\begin{figure}[htbp]
    \centering
    \begin{subfigure}{0.48\textwidth}
        \centering
        \includegraphics[width=\textwidth]{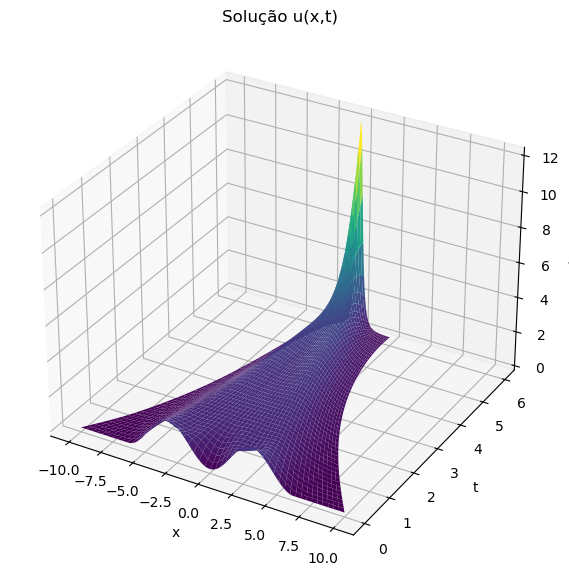}
        \caption{Numerical solution $u(x,t)$.}
        \label{fig:u}
    \end{subfigure}
    \hfill
    \begin{subfigure}{0.48\textwidth}
        \centering
        \includegraphics[width=\textwidth]{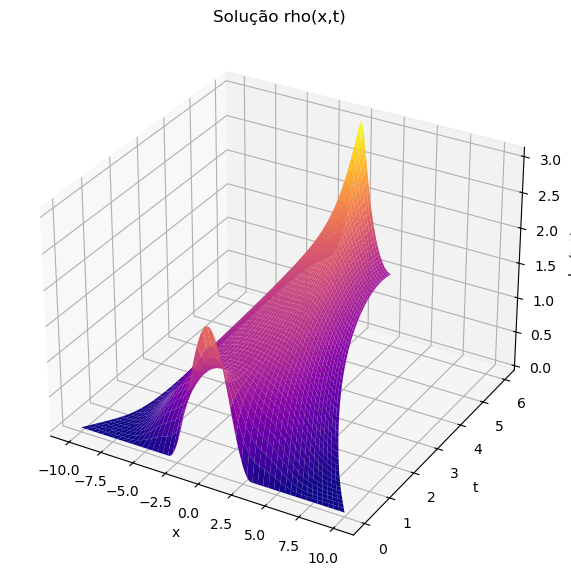}
        \caption{Numerical solution $\rho(x,t)$.}
        \label{fig:rho}
    \end{subfigure}
    \caption{Numerical solutions for $u(x,t)$ and $\rho(x,t)$.}
    \label{fig:solutions}
\end{figure}

Figures \ref{fig:u} and \ref{fig:rho} show the numerical solution of \eqref{equ4}, with the constants \eqref{eq:constants}, the curves \eqref{eq:alphabeta}, and the initial data \eqref{eq:utau-rhotau}, obtained by solving \eqref{equ5} on the fixed domain $\mathcal{I}=[0,1]$ and recovering $u,\rho$ through \eqref{recuperar}.

\vspace{0.5cm}
\noindent {\bf Example II.} For the numerical example associated with the Chemotaxis problem \eqref{equ4}, consider the constants
\begin{equation}\label{eq:constants-0}
    a=1,\qquad b=1,\qquad \kappa=2,\qquad c=1,\qquad d=2.
\end{equation}
The boundary curves are
\begin{equation}\label{eq:alphabeta-0}
    \alpha(t) = -10+\frac{8t}{\sqrt{9+t^{2}}}, \qquad \beta(t) = 10-\frac{8t}{\sqrt{9+t^{2}}},
\end{equation}
so that $\gamma(t) = \beta(t)-\alpha(t) = 20-\frac{16t}{\sqrt{9+t^{2}}}$. Their derivatives are
\begin{equation*}
    \alpha'(t) = \frac{72}{(9+t^{2})^{3/2}}>0, \quad \beta'(t) = -\frac{72}{(9+t^{2})^{3/2}}<0, \quad
    \gamma'(t) = -\frac{144}{(9+t^{2})^{3/2}}<0, \quad \forall t\in\mathbb{R}.
\end{equation*}
In particular $\gamma$ is strictly decreasing, with
\begin{equation*}
    \gamma(0)=20, \qquad \lim_{t\to+\infty}\gamma(t)=4, \qquad \lim_{t\to-\infty}\gamma(t)=36,
\end{equation*}
so $\gamma(t)\in(4,36)$ for all $t\in\mathbb{R}$: the domain $\mathcal{I}_t$ contracts monotonically, from $[-10,10]$ at $t=0$ towards $(-4,4)$ as $t\to+\infty$ (at $t=20$, $\gamma(20)\approx 4.18$).

For $\tau=0$ the initial data are related by
\begin{equation}\label{eq:utau-rhotau-0}
u_0(x) =
\begin{cases}
1, & -10\leq x< -\pi,\\[1mm]
\cos(x)+2, & -\pi\leq x\leq \pi,\\[1mm]
1, & \pi < x \leq 10.
\end{cases} \hspace{0,5cm} \text{ and }\hspace{0,5cm} \rho_0(x) = 1, \, \, -10\leq x\leq 10.
\end{equation}

\begin{figure}[htbp]
    \centering
     \begin{subfigure}{0.48\textwidth}
        \centering
    \includegraphics[width=\textwidth]{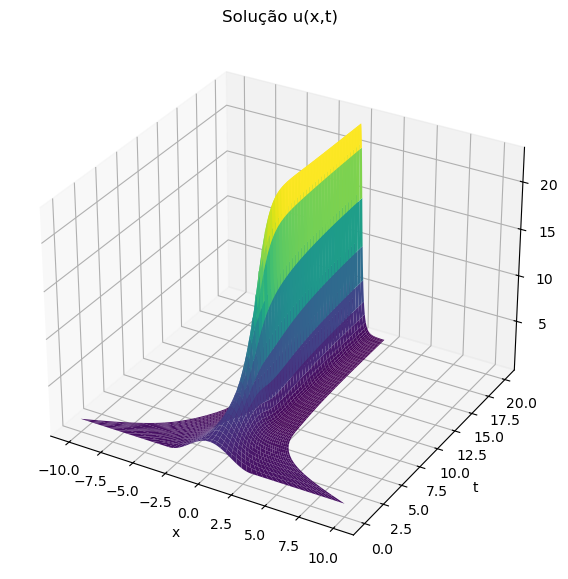}
    \caption{Numerical solution $u(x,t)$.}
    \label{fig:u-II}
\end{subfigure}
\hfill
    \begin{subfigure}{0.48\textwidth}
        \centering
    \includegraphics[width=\textwidth]{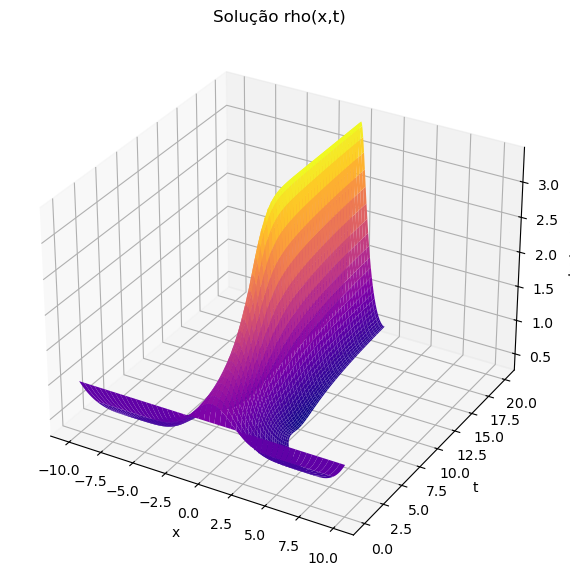}
    \caption{Numerical solution $\rho(x,t)$.}
    \label{fig:rho-II}
    \end{subfigure}
    \caption{Numerical solutions for $u(x,t)$ and $\rho(x,t)$.}
    \label{fig:solutions-II}
\end{figure}

Figures \ref{fig:u-II} and \ref{fig:rho-II} show the numerical solution of \eqref{equ4}, with the constants \eqref{eq:constants-0}, the curves \eqref{eq:alphabeta-0}, and the initial data given in \eqref{eq:utau-rhotau-0}.

\section*{Open questions}

Global existence was obtained under (H1)--(H2) alone, while dissipativity also required (H3) and the small gain condition (H4). On a fixed domain,
(H2)--(H4) hold trivially and we recover \cite{OY01}. Two questions arise
naturally.

\emph{Is (H2) necessary for global existence?} On a moving boundary, the Neumann condition is not a no-flux condition:
$\frac{d}{dt}\int_{\mathcal I_t}u(x,t)\,dx=\beta'(t)u(\beta(t),t)-\alpha'(t)u(\alpha(t),t)$.
Under (H2) the total population is nonincreasing, whereas on expanding domains mass is injected through the boundary. In numerical simulations on expanding domains, which we performed but do not include in this paper, we observed solutions that appear to blow up in finite time, in contrast with the fixed-domain case \cite{OY01,Yagi97}.

\begin{question}\label{conj2}
There exist $\alpha,\beta\in C^1(\mathbb R)$ satisfying (H1), with $\gamma'>0$
on some interval, and nonnegative initial data for which the solution of
\eqref{equ4} blows up in finite time, that is,
$\|u(t)\|_{L^\infty(\mathcal I_t)}\to\infty$ as $t\to T^{*-}$, for some
$T^*<\infty$.
\end{question}

In order to prove the dissipativity of the evolution process, we needed an estimate of singular Gronwall type with an explicit exponential decay rate. This led us to assume the small gain condition (H4), which controls the velocities of the endpoints of the time-dependent interval $\mathcal I_t$. Is this hypothesis really necessary, or is it merely a difficulty arising from the technique used?

\begin{question}\label{conj1}
Under (H1)--(H3), the evolution process associated with \eqref{equ4} is pullback dissipative and possesses a pullback attractor, regardless of the size of $\|\alpha'\|_{L^\infty(\mathbb R)}$ and $\|\beta'\|_{L^\infty(\mathbb R)}$.
\end{question}

\subsection*{Acknowledgments} 

J.C. Huaccha Neyra was partially supported by grants \#2025/25292-1, FAPESP.

\section*{Declarations} 

\medskip \textbf{Conflict of interests.} On behalf of all authors, the corresponding author states that there is no conflict of interest.

\medskip \textbf{Data availability.} We do not analyze or generate any datasets, because our work proceeds within a theoretical and mathematical approach.

\medskip \textbf{AI use.} During the preparation of this work the authors used Claude (Anthropic) to improve the language, readability and clarity of the text. After using this tool, the authors reviewed and edited the content as needed and take full responsibility for the content of the published article.

\end{document}